\documentclass[11pt]{article}
\usepackage{latexsym}
\usepackage{amsfonts,amssymb,amsmath}
\usepackage{fontawesome}

\newtheorem{thm}{Theorem}[section]
\newtheorem{cor}[thm]{Corollary}
\newtheorem{lem}[thm]{Lemma}
\newtheorem{pro}[thm]{Proposition}
\newtheorem{defn}[thm]{Definition}

\newtheorem{notat}[thm]{Notation}

\newcommand{\ov }{\overline }
\newcommand{\para}{\parallel}
\newcommand{\minus}{\smallsetminus}

\title{The word problem of the Brin-Thompson group is {\sf coNP}-complete}

\author{J.C.\ Birget}

\date{\scriptsize{
10.ii.2020}}

\begin{document}
\maketitle

\begin{abstract}
We prove that the word problem of the Brin-Thompson group $n V$ over a 
finite generating set is {\sf coNP}-complete for every $n \ge 2$. 
It is known that $\{n V : n \ge 1\}$ is an infinite family of infinite, 
finitely presented, simple groups. 
We also prove that the word problem of the Thompson group $V$ over a 
certain infinite set of generators, related to boolean circuits, is 
{\sf coNP}-complete. 
\end{abstract}

\section{Introduction} 

The group $nV$ was introduced by Brin \cite{Brin1} as an $n$-dimensional 
generalization of Richard Thompson's group $V$, for any 
positive integer $n$ (with $\,1V = V$). 

Brin proved that $2V$ is finitely generated and simple, that $V$ is not 
isomorphic to $2V$ \cite{Brin1}, that $2V$ is finitely presented 
\cite{Brin2}, and that all $nV$ are simple \cite{Brin3}. 
Hennig and Mattucci \cite{HennigMattucci} show that all $nV$ are finitely 
presented. Bleak and Lanoue \cite{BleakLanoue} show that all $nV$ are 
non-isomorphic.
In short, the groups $nV$ are an infinite family of infinite, finitely 
presented, simple groups.

The word problem of $nV$ is decidable, as is easy to see from the 
definition of $nV$. The main result of the present paper is the following.

\begin{thm} \label{2g21wpCoNP}
 \ The word problem of $nV$ over any finite generating set is
{\sf coNP}-complete, for all $n \ge 2$.
\end{thm}
 
\bigskip

\noindent Remarks on the theorem:

This is only the second example of a finitely presented group with
{\sf coNP}-complete word problem; the first example appeared in 
\cite{BiCoNP}. This is also the first ``naturally occurring'' example of a
finitely presented group with either {\sf NP}-complete or {\sf coNP}-complete
word problem. 
The proof of Theorem \ref{2g21wpCoNP} strengthens the 
connection between acyclic circuits and finite group presentations;  
such a connection already played a crucial role in \cite{BiCoNP}.

The Theorem implies that if ${\sf NP} \ne {\sf coNP}$ then the Dehn
function of $nV$ (for $n \ge 2$) has no polynomial upper bound; more 
strongly, $nV$ cannot be embedded into a finitely presented group with 
polynomially bounded Dehn function (by \cite{SBR, BORS}).

The Theorem implies that if ${\sf P} \ne {\sf NP}$ then $2V$ is not 
embeddable into $V$. It was proved recently \cite[Coroll.\ 11.20]{MatteBon}
that $(n+1)V$ does not embed into $nV$ for any $n \ge 1$.

The groups $nV$ for $n \ge 2$ are the first examples of finitely presented
simple groups whose word problem is harder than {\sf P} (if 
${\sf P} \ne {\sf NP}$).\footnote{      
 \ The Higman-Thompson groups $G_{k,s}$ have their word problem in {\sf P} 
(in fact in {\sf coCFL}, by Lehnert and Schweitzer \cite{LehnSchw}). 
For other currently known finitely presented infinite simple groups 
(Meier \cite{Meier83, Meier85}, 
R\"over \cite{Roever}, 
Burger and  Mozes \cite{BurgMoz}, 
Lodha \cite{YLodha}),
the complexity of the word problem has not been studied, but appears to be 
in {\sf P}.} 
Finitely presented infinite simple groups are related to the Boone-Higman
theorem \cite{BooneHigman}. In \cite{BooneHigman} the authors ask whether 
their theorem can be strengthened as follows: 
{\it Does a finitely generated group $G$ have a decidable word problem iff 
$G$ is embeddable into a finitely presented simple group?} 
In contrast, it was observed in \cite[Section 1]{BiCoNP} that all known 
finitely presented simple groups have a word problem of very low complexity; 
even {\sf coNP} is a low complexity class on the scale of all decidable 
problems. The enormous gap between what is asked, and what has been 
observed so far motivates the following.

\smallskip

\noindent 
{\it Question: Are the computational complexities of the word problems of
all finitely presented simple groups unbounded?} 

\smallskip

\noindent More precisely, the negation of the question is:  
Is there a time-constructible total function $t$ such that the word problems 
of all finitely presented simple groups belong to ${\sf DTime}(t)$? (See 
e.g.\ \cite{HU} for the definitions of ``time-constructible'' and 
``${\sf DTime}(t)$''.) 
In case of a negative answer, the Boone-Higman question also has a negative
answer. If the answer is positive then there is a chance that the  
Boone-Higman question has a positive answer; in that case, the proof of the 
answer to the Question above might be easier than the proof of a
strengthened Boone-Higman theorem, and could be a useful step along the way. 

\bigskip

\noindent {\bf Overview:} In section 2 we define the Higman-Thompson groups 
$G_{k,1}$ and the Brin-Thompson groups $nV$ and $n G_{k,1}$ by (partial) 
actions on finite strings, or $n$-tuples of strings.  For this, the concept 
of prefix code of strings is generalized to the concept of joinless code of 
$n$-tuples of strings.
For the study of the computational complexity of the word problem, the
string-based formalism is more convenient than the geometric approach.
It follows fairly directly that the word problem of $nV$ over a finite 
generating set belongs to {\sf coNP} (section 3). 

The proof of {\sf coNP}-hardness is given in section 4. It goes through 
several steps, following the same strategy as the first half of 
\cite{BiCoNP} (where it was proved that a certain subgroup of $G_{3,1}$, 
over a certain infinite generating set, has a {\sf coNP}-complete word 
problem. 
Based on this we show that the Thompson group $V$, over a certain infinite 
generating set, has a {\sf coNP}-complete word problem. This infinite 
generating set of $V$ consists of a finite generating set, together with 
all the bit-position transpositions $\tau_{i,i+1}$ (where $\, \tau_{i,i+1}:$
$\, x_1 \, \ldots \, x_{i-1} \, x_i \, x_{i+1} \, x_{i+2} \, \ldots \, $ 
$\longmapsto$ 
$\, x_1 \, \ldots \, x_{i-1} \, x_{i+1} \, x_i \, x_{i+2} \, \ldots \ $).
An alternative approach, based on bijective circuits and the work of Jordan
\cite{Jordan}, is described in subsection 4.5. 
Finally, we show that $\tau_{i,i+1}$ can be expressed by $\tau_{1,2}$ and 
the shift $\sigma$.  This reduces the word problem of $V$, over an infinite 
generating set that includes position transpositions, to the word problem of
$2V$ over a finite generating set (subsection 4.6). 

\bigskip

\noindent {\bf Summary of abbreviations and notations:}

\smallskip

\noindent
-- The word {\it function} in this paper means partial function. The domain 
of a function $f: X \to Y$ is denoted by ${\rm Dom}(f) \subseteq X$, and the
image by ${\rm Im}(f) \subseteq Y$. Most often, the sets $X$ and $Y$ will be 
free monoids $A^*$, or Cantor spaces $A^{\omega}$, or their direct powers 
$nA^*$ or $nA^{\omega}$.  

\noindent
-- $A^*$, \ the free monoid freely generated by $A$, a.k.a.\ the set of all   
   strings over $A$;

\noindent
-- $\varepsilon$, \ the empty string;

\noindent
-- $A^+$, \ the free semigroup; $A^+ = A^* \minus \{\varepsilon\}$;

\noindent
-- $|x|$, \ the length of the string $x \in A^*$;

\noindent
-- $x \le_{\rm pref} y$, \ $x$ ($\in A^*$) is a prefix of $y$ 
  ($\in A^* \cup A^{\omega}$);

\noindent
-- $x \para_{\rm pref} y$, \ $x$ is prefix-comparable with $y$;

\noindent
-- $nA^*$, $nA^{\omega}$, \ the $n$-fold cartesian product 
  {\large \sf X}$_{_{i=1}}^{^n} A^*$, respectively  
  {\large \sf X}$_{_{i=1}}^{^n} A^{\omega}$;

\noindent
-- $(\varepsilon)^n$, \ the $n$-tuple of empty strings;

\noindent
-- $A_{\varepsilon,n} = \, \bigcup_{1 \le i \le n}$
$\{\varepsilon\}^{i-1} \times A \times \{\varepsilon\}^{n-i}$, \ the 
unique minimum generating set of the monoid $nA^*$;

\smallskip

\noindent
-- $\ell(x)$, 
 \ \ $\max\{|x_1|, \,\ldots, |x_n|\} \ $ if $\,x = (x_1, \,\ldots, x_n)$ 
     $\in nA^*$;

\noindent
-- $x \le_{\rm init} y$, \ $x$ ($\in nA^*$) is an initial factor of $y$ 
($\in nA^* \cup nA^{\omega}$);

\noindent
-- {\sc dag}, \ directed acyclic graph;

\noindent
-- $f|_M$, \ the restriction of a function $f$ to a set $M$.

\section{Definition of {\boldmath $n V$} based on strings}

The standard definitions in computational complexity require strings as 
inputs.
Brin's original definition of $n V$ uses geometric actions, but for the proof 
of {\sf coNP}-completeness of the word problem of $n V$ we also need a 
(partial) action of $n V$ on $n$-tuples of strings.
The groups $nV$ are generalizations of $V$.  We first look at $V.$

\subsection{Definition of {\boldmath $V$} based on strings}

The group $V$ can be defined in many ways; see e.g.\ \cite{Th0, McKTh, Th, 
Hig74, Scott, CFP}.  We will mostly use two definitions of $V$ from 
\cite{BiThomps} (which are is similar to \cite{Scott}, except that we use 
the terminology of prefix codes, right ideals, and right-ideal morphisms). 

We recall some standard notions. 
An alphabet is any finite set, although we mostly use $\{0,1\}\,$ (the bits),
and $\{0,1, \, \ldots, k-1\} \,$ for any integer $k \ge 2$.
For an alphabet $A$ and $m \in {\mathbb N}$, $A^m$ denotes the 
set of sequences of length $m$ over $A$ (called set of {\it strings} of 
length $m$), and for $x \in A^m$ we say that $|x| = m$ (i.e., the 
{\it length} of $x$ is $m$); $A^{\le m}$ is the set of strings of 
length $\le m$.  
The empty string is denoted by $\varepsilon$, and $|\varepsilon| = 0$.
The set of all strings over $A$ is denoted by $A^*$, and the set of all
infinite strings indexed by the ordinal $\omega$ is denoted by $A^{\omega}$. 
By default a ``string'' is finite; for infinite strings we explicitly say 
``infinite''.
For $x_1, x_2 \in A^*$ the {\it concatenation} is denoted by $x_1 x_2$ or
$x_1 \cdot x_2$; it has length $|x_1| + |x_2|$.
For two subsets $S_1, S_2 \subseteq A^*$, we define the concatenation by 
$\, S_1 \cdot S_2 = \{x_1 \cdot x_2 :$ $x_1 \in S_1$ and $x_2 \in S_2\}$.

For $x, p \in A^*$ we say that $p$ is a {\it prefix} of $x$ iff 
$(\exists u \in A^*) \, x = pu$; this is denoted by $p \le_{\rm pref} x$. 
Two strings $x, y \in A^*$ are called {\it prefix-comparable} (denoted 
by $x \para_{\rm pref} y \,$) iff 
$\, x \le_{\rm pref} y \,$ or $ \, y \le_{\rm pref} x$. 
A {\it prefix code} (a.k.a.\ a prefix-free set) is any subset $P \subset A^*$
such that for all $p_1, p_2 \in P$: $\, p_1 \para_{\rm pref} p_2 \,$ 
implies $p_1 = p_2$.
A {\it right ideal} of $A^*$ is, by definition, any subset $R \subseteq A^*$ 
such that $R = R \cdot A^*$. A subset $C \subseteq R$ is said to generate 
$R$ as a right ideal iff $R = C \cdot A^*$. 
It is easy to prove that every finitely generated right ideal is generated 
by a unique finite prefix code, and this prefix code is the minimum 
generating set of the right ideal (with respect to $\subseteq$). 
By definition, a {\it maximal prefix code} is a prefix code $P \subset A^*$ 
that is not a strict subset of any other prefix code of $A^*$. 
An {\it essential right ideal} is, by definition, a right ideal 
$R \subseteq A^*$ such that all right ideals of $A^*$ intersect $R$ (i.e., 
have a non-$\varnothing$ intersection with $R$). 
It is well known (see e.g.\ \cite[Lemma 8.1]{BiThomps}) that a right ideal 
$R \subseteq A^*$ is essential iff the unique prefix code that generates 
$R$ is maximal. 

A {\it right ideal morphism} of $A^*$ is, by definition, a function 
$f: A^* \to A^*$ such that for all $x \in {\rm Dom}(f)$ and all $w \in A^*$:
 \ $f(xw) = f(x) \ w$.
In that case, ${\rm Dom}(f)$ is a right ideal; one easily proves that 
${\rm Im}(f)$ is also a right ideal. The prefix code that generates
${\rm Dom}(f)$ is denoted by ${\rm domC}(f)$, and is called the {\it domain
code} of $f$; the prefix code that generates ${\rm Im}(f)$ is denoted by 
${\rm imC}(f)$, and is called the {\it image code}.
We are interested in the following monoid:

\medskip
 
${\cal RI}_A^{\sf fin}$ $\, = \,$  $\{ f : f$ is a right ideal morphism of 
    $A^*$ such that $f$ is {\em injective}, and  

\hspace{0.9in}
    ${\rm domC}(f)$ and ${\rm imC}(f)$ are {\em finite maximal} 
    prefix codes\}.

\medskip

\noindent We usually write ${\cal RI}^{\sf fin}$ for 
${\cal RI}_A^{\sf fin}$ since we usually just deal with one alphabet $A$ at 
a time.  It is proved in \cite[Prop.\ 2.1]{BiThomps} that every 
$f \in {\cal RI}^{\sf fin}$ is contained in a unique $\subseteq$-maximum 
right ideal morphism in ${\cal RI}^{\sf fin}$; this is called the 
{\it maximum extension} of $f$. 
The Higman-Thompson group $G_{k,1}$ (where $k = |A|$) is a homomorphic 
image of ${\cal RI}^{\sf fin}$:

\begin{defn} \label{ThompsV} {\bf (Thompson group $V$ and Higman-Thompson 
groups $G_{k,1}$).}
 \ The Thompson group $V$, as a set, consists of the 
right ideal morphisms $f \in {\cal RI}_{\{0,1\}}^{\sf fin}$ that are maximum
extensions in ${\cal RI}_{\{0,1\}}^{\sf fin}$. 
The multiplication in $V$ consists of composition, followed by maximum
extension.

\smallskip

The same definition for ${\cal RI}_A^{\sf fin}$ with 
$A = \{0,1, \, \ldots, k-1\}$ yields the Higman-Thompson group $G_{k,1}\,$ 
for every $k \ge 2$; $\, V = G_{2,1}$. 
\end{defn}

Every element $f \in {\cal RI}^{\sf fin}$ (and in particular, every 
$f \in G_{k,1}$) is determined by the restriction of $f$ to ${\rm domC}(f)$. 
This restriction $\,f_{{\rm domC}(f)}: {\rm domC}(f) \to {\rm imC}(f)\,$ is 
a finite bijection, called the {\it table} of $f$ \cite{Hig74}. 
Obviously, $f$ ($\in {\cal RI}^{\sf fin}$) determines ${\rm domC}(f)$ and 
hence a unique table.  
When we use tables we do not always assume that $f$ is a maximum extension.
The well known tree representation of $G_{k,1}$ is obtained by using the 
{\em prefix trees} of ${\rm domC}(f)$ and ${\rm imC}(f)$.

\begin{lem} \label{extStep} 
 \ Let $P, Q \subset A^*$ be finite maximal prefix codes. 
The right ideal morphism $f \in {\cal RI}^{\sf fin}$ determined by a 
table $F$: $P \to Q$ can be extended \ iff \ there exist $s, t \in A^*$ such 
that for every $\alpha \in A$:
 \ $s \alpha \in P$, \ $t \alpha \in Q$, and $\, F(s \alpha) = t \alpha$.

In that case, $f$ can be extended by defining $\, f(s) = t$. The table for
this extension is obtained be replacing 
$\,P\,$ by $\,(P \minus s A) \cup \{s\}$, 
$\,Q\,$ by $\,(Q \minus qA) \cup \{q\}$, and  
$\, \{(s \alpha, t \alpha) : \alpha \in A\}\,$ by $\,\{(s, t)\}$. 

This is called an {\em extension step} of the table $F$.  
\end{lem}
{\sc Proof.} See \cite[Lemma 2.2]{BiThomps} and \cite{Hig74}. 
 \ \ \  \ \ \  $\Box$

\medskip

Since in an extension step the cardinality of ${\rm domC}(f)$ decreases,
only finitely steps are needed to reach the maximum extension of $f$; the
number of steps is $\, < |{\rm domC}(f)|$.

Based on the representation of the elements of $V$ (and of $G_{k,1}$) by 
tables, one can show easily that the word problem of these groups is in
{\sf P}. A much stronger result is that the word problem is in {\sf coCFL} 
(the set of languages whose complement is context-free) \cite{LehnSchw}; 
{\sf coCFL} is a strict subclass of the parallel complexity class 
${\sf AC}^1$, which is a subclass of {\sf P} (see e.g., \cite{HandbookA}).

\bigskip

\noindent {\bf The $A^{\omega}$ definition of {\boldmath $G_{k,1}$:}}
Maximality of finite prefix codes has the following characterization in
terms of $A^{\omega}$. A finite prefix code $P \subset A^*$ is maximal iff 
$\, P A^{\omega} = A^{\omega}$. 
(This is not true for infinite prefix codes; a counter example is 
$\, 0^*1$.)

It follows that every element $f \in G_{k,1}$ determines a permutation of
$A^{\omega}$.
Conversely, let $P \subset A^*$ be a finite maximal prefix code. Then for
every $w \in A^{\omega}$ there exists a unique $p \in P$ and 
$v \in A^{\omega}$ such that $w = p v$. Let $f$ be a permutation of 
$A^{\omega}$ for which there exists a table $F$: $P \to Q$ such $f$ is 
defined by $\, f(p v) = F(p) \ v$ \ (for every $p \in P$ and 
$v \in A^{\omega}$).  Then $f \in G_{k,1}$.

Thus, $G_{k,1}$ can be defined as a certain group of permutations of
$A^{\omega}$.

\begin{lem} \label{VequivCommRestr}
 \ Let $F_1$: $P_1 \to Q_1$ and $F_2$: $P_2 \to Q_2$ be two tables that 
determine, respectively, the right ideal morphisms $f_1, f_2 \in$
${\cal RI}^{\sf fin}$. Then the following are equivalent:    

\smallskip

\noindent
{\rm (1)} \ \ $F_1$ and $F_2$ determine the same element of $G_{k,1}$ (by
maximum extension);  \\  
{\rm (2)} \ \ $f_1$ and $f_2$ have the same maximum extension in 
 ${\cal RI}^{\sf fin}$; \\   
{\rm (3)} \ \ $f_1$ and $f_2$ have a common restriction in 
 ${\cal RI}^{\sf fin}$; \\  
{\rm (4)} \ \ $f_1$ and $f_2$ have a common restriction to an essential 
 right ideal of $A^*$; \\   
{\rm (5)} \ \ $F_1$ and $F_2$ determine the same function on 
    $A^{\omega}$; \\      
{\rm (6)} \ \ $f_1$ and $f_2$ determine the same function on $A^{\omega}$.  
\end{lem}
{\sc Proof.} (1) and (2) are equivalent by the definition of $G_{k,1}$. 
(2) implies (3) (which implies (4)): 
The intersection $f_1 \cap f_2$ is a common restriction; 
by \cite[Lemma 8.3]{BiThomps}, ${\rm Dom}(f_1) \cap {\rm Dom}(f_2)$ is an
essential right ideal. Moreover,
${\rm domC}(f_1 \cap f_2) \subset {\rm domC}(f_1) \cup {\rm domC}(f_2)$;
hence ${\rm domC}(f_1 \cap f_2)$ is finite.
(4) implies (2) by uniqueness of maximum extensions in 
${\cal RI}^{\sf fin}$ (see \cite[Lemma 2.1]{BiThomps}, which does not 
require finiteness of prefix codes). (3) implies (5) in an obvious way.
And (5) implies (1), based on finiteness and uniqueness of maximum 
extension. (5) and (6) are obviously equivalent. 
 \ \ \ $\Box$

\bigskip  

\noindent {\bf The piecewise linear definition of $V$:} 
Brin's definition of $nV$ extends the definition of $V$ as given in 
\cite{CFP}; the latter is based on piecewise linear actions on the 
interval $[0,1]$ $\subset$ ${\mathbb R}$. We use half-open intervals, so 
neighboring intervals do not intersect; however, when the right boundary 
is 1, we use ``$1]$''.  
The boundary-points of the subintervals that appear are binary rational 
numbers (i.e., the denominator is a power of 2). 
A string $s = s_1 \ldots s_m \in \{0,1\}^*$ with $m = |s|$ determines the 
half-open subinterval 
$\,[0.s, \ 0.s + 2^{-|s|}[\,$; \ but if $s + 2^{-|s|} = 1$ then we take 
$\,[0.s, \, 1]$, i.e., in that case we close the interval.
Here, $0.s$ is a rational number written in fractional {\em binary 
representation}; i.e., $0.s = \sum_{i=1}^m s_i \, 2^{-i}$.
E.g., 01100 (of length 5) determines the subinterval
$\, [0.011, \ 0.011 + 2^{-5}[$ $ \, = \,$ 
$[0.011, \ 0.01101[ \,$.

\subsection{Right ideals of {\boldmath $nA^*$}}

Here we completely develop the string description of $nV$, which is briefly
alluded to in \cite[subsection 4.3]{Brin1}. A hybrid string-geometric
description was used in \cite{BleakLanoue} (where some crucial concepts 
appear only in geometric form).  Our description is entirely based on 
strings, but the correspondence with geometric concepts is often pointed out.
The present subsection focuses on finitely generated right ideals of $nA^*$; 
in the next subsection, $n V$ will be defined based on right-ideal morphisms 
of $nA^*$. 

\smallskip

As before, let $A$ be a finite alphabet of cardinality $k \ge 1$, usually
denoted by $\{0, \, \ldots, k-1\}$ or $\{a_0, \,\ldots, a_{k-1} \}$.
The $n$-fold cartesian product \ {\large \sf X}$_{_{i=1}}^{^n} A^* \, $ will 
be denoted by $nA^*$; we choose this notation in analogy with the notation 
$nV$, and also in order to avoid confusion with $n$-fold concatenation (of 
the form  
$\, S^n = \{s_1 \cdot \ldots \cdot s_n :\, s_1, \ldots, s_n \in S\}$
$\subseteq A^*$).
Similarly, $nA^{\omega}$ denotes the $n$-fold cartesian product
 \ {\large \sf X}$_{_{i=1}}^{^n} A^{\omega}$.  
Multiplication in $\, nA^*$ is done coordinatewise, i.e., $nA^*$ is the 
direct product of $n$ copies of the free monoid $A^*$. For $u \in nA^*$ we 
denote the coordinates of $u$ by $u_i \in A^*$, for $1 \le i \le n$; i.e., 
$u = (u_1, \, \ldots, u_n)$.

\smallskip

{\bf Geometrically:}  $x = (x_1, \ldots, x_n) \in n \{0,1\}^*$ represents 
the hyperrectangle 
$\,${\large \sf X}$_{_{i=1}}^{^n} [0.x_i, \ 0.x_i + 2^{-|x_i|}[ \,$ 
(except that ``$0.x_i + 2^{-|x_i|}[$'' is replaced by ``$1]$'' if 
$\,0.x_i + 2^{-|x_i|} = 1$).
The measure of this hyperrectangle is
 \ $2^{-(|x_1| \, + \ \ldots \ + \, |x_n|)}$.
In particular, $(\varepsilon)^n$ represents $[0,1]^n$ and has measure 1. 

\smallskip

The concept of prefix is similar to the one in $A^*$, but in order to avoid 
confusion we will use the phrase ``initial factor''.
So the {\em initial factor order} on $nA^*$ is defined for $u, v \in nA^*$
by $u \le_{\rm init} v$ iff there exists $x \in nA^*$ such that $u x = v$.
In a similar way we have the concepts of comparability (denoted by 
$\|_{\rm init}$), right ideal, generating set of a right ideal, and essential
right ideal. 
It is easy to prove that $\, u \le_{\rm init} v$ in $nA^*$ iff 
$u_i \le_{\rm pref} v_i$ for all $i= 1, \, \ldots, n$.
For any $u, v \in nA^*$ there exists a unique $\le_{\rm init}$-maximum 
common initial factor, denoted by $u \wedge v$. 
In terms of coordinates, $(u \wedge v)_i = u_i \wedge_{\rm pref} v_i$, where 
$u_i \wedge_{\rm pref} v_i$ is the longest common prefix of the strings $u_i$ 
and $v_i$.

An {\em initial factor code} is a set $S \subset nA^*$ such
that no two different elements of $S$ are $\le_{\rm init}$-comparable.

As we shall see, a crucial way in which $nA^*$ with $n \ge 2$ differs from 
$A^*$ concerns the {\em join} operation with respect to $\,\le_{\rm init}$. 
For all $n$, the join of $u,v \in nA^*$ is defined by $ \ u \vee v$ $\,=\,$  
$\min_{\le_{\rm init}}\{z \in nA^* :$
$ \, u \le_{\rm init} z$ and $v \le_{\rm init} z\}$. 
Of course, $u \vee v$ does not always exist.

\begin{defn} \label{DEFjoinless}
 \ A set $S \subset nA^*$ is {\em joinless} iff no two elements of $S$ 
have a join with respect to $\, \le_{\rm init}$ in $nA^*$.
Joinless sets will be called {\em joinless codes}, since they are 
necessarily initial factor codes.

A set $S \subset nA^*$ is a {\em maximal} joinless code iff 
$\,S$ is $\subseteq$-maximal among the joinless codes of $nA^*$. 
(In other words, adding a new element to a maximal joinless code $S$ 
results in a set, some of whose elements have joins.)

A right ideal $R \subseteq nA^*$ is called {\em joinless generated}
iff $R$ is generated, as a right ideal, by a joinless code. 
\end{defn}
(About the grammar: ``Joinless{\bf ly} generated'' would not make sense 
since it is not the generating process that is joinless.)

\medskip

\noindent
{\bf Examples:} Not every initial factor code is joinless;  e.g.,
$\{(\varepsilon, 0), \ (0, \varepsilon)\}$ is an initial factor code where
$\, (\varepsilon, 0) \vee (0, \varepsilon) = (0,0)$.
An example of a maximal joinless code is 
$\, \{(\varepsilon, 0), \ (0, 1), \ (1, 1)\}$. 
A maximal joinless code is usually not maximal as an initial factor code; 
for example, in $\, \{(\varepsilon, 0), \ (0, 1), \ (1, 1)\} \,$ one could 
add $(00, \varepsilon)$; the result would be a initial factor code (that is 
not joinless).
The only maximal joinless code that is also maximal as an initial factor 
code is $\{(\varepsilon, \varepsilon)\}$. 

From here on, a joinless code will be called maximal if it is maximal as a 
joinless code.

\medskip

\noindent {\bf Connection with the geometric description:}
For $u, v \in nA^*$ we have $v \le_{\rm init} u$ iff
the hyperrectangle $u$ is contained in the hyperrectangle $v\,$ (i.e., 
$\le_{\rm init}$ corresponds to $\supseteq$); note that ``shorter'' 
$n$-tuples correspond to ``larger'' hyperrectangles. The join $u \vee v$ 
represents the hyperrectangle obtained by intersecting the 
hyperrectangles $u$ and $v$ (so $\vee$ corresponds to $\cap$). Note that 
$u \vee v$ does not exist iff the intersection is the empty set (since the 
empty set is not a hyperrectangle). 
The meet $u \wedge v$ (which always exists) does not represent the union, 
nor the smallest hyperrectangle that contains $u$ and $v$, but the smallest hyperrectangle representable by an $n$-tuple in $nA^*$ that contains $u$ 
and $v$.
Joinlessness of a code means that any two hyperrectangles in the chosen 
subdivision of $[0,1]^n$ are {\em disjoint} as sets. A joinless code is 
maximal iff its hyperrectangles form a {\em tiling} of $[0,1]^n$.
In an initial factor code, $\le_{\rm init}$-incomparability means that no 
hyperrectangle in the code is contained in another one.

\smallskip

\noindent {\bf Examples} (for the correspondence between strings and 
geometry): Fig.\ 1 shows a few elements of $\,2\, \{0,1\}^*$. 
The large square $[0,1] \times [0,1]$ is represented by
$(\varepsilon, \varepsilon)$.
The numbers use fractional binary representation; e.g.,
$0.1101 = \frac{1}{2} + \frac{1}{4} + \frac{1}{16}$.

$(0, 00) \in 2\, \{0,1\}^*$ represents the rectangle 
$\,[0, \,0.1[ \,\times\, [0, \,0.01[$ \ (horizontally hashed); and 
$(010, 0)$ represents $\,[0.01, \,0.011[ \,\times\, [0, \,0.1[$ 
 \ (vertically hashed).  The join
$(010, 00) = (0, \, 00) \vee (010, \,0)\,$ represents 
$\,[0.01, \,0.011[ \,\times\, [0, \,0.01[$ \ (doubly hashed).
$(0,0) = (0, 00) \,\wedge\, (010, 0)\,$ represents 
$\,[0, 0.1[ \,\times\, [0, 0.1[$. 

The rectangle $[0.1, \, 0.11[ \,\times\, [0.1101, \, 0.111[$ is represented 
by $(10, 1101)$ (horizontally hashed). 
And $[0.1, \, 1] \times [0.111, \, 0.1111]$ is represented by 
$(1, 1110)$ (vertically hashed). 
Here, $(10, 1101) \vee (1, 1110)$ does not exist, and the meet
$(10, 1101) \wedge (1, 1110) = (1, 11)$ represents 
$[0.1, \, 1] \times [0.11,\, 1]$.

\unitlength=1.0mm       
\begin{picture}(110,60)

\put(5,10){\framebox(40,40)}    
\put(5,10){\framebox(20,20)}    

\put(-2,0){\makebox(0,0)[cc]{\sf Fig. 1}}

\put(5,10){\framebox(20,10)}    
\put(15,10){\framebox(5,20)}    

\put(25,42.5){\framebox(10,2.5)}  
\put(25,45){\framebox(20,2.5)}  

\put(25,40){\framebox(20,10)}

\linethickness{0.05mm}

\multiput(5,10)(0,0.5){20}{\line(1,0){20}}  
\multiput(15,10)(0.5,0){10}{\line(0,1){20}}   

\multiput(25,42.5)(0,0.5){5}{\line(1,0){10}}    
\multiput(25,45)(0.5,0){40}{\line(0,1){2.5}}   

\end{picture}

\bigskip

\bigskip

For $u,v \in A^*$, $\, u \vee_{\rm pref} v$ exists in $A^*$ iff $u$ and $v$ 
have a common upper bound for $\le_{\rm pref}$. This holds iff 
$u \para_{\rm pref} v$; in that case,
$u \vee_{\rm pref} v = u$ if $v \le_{\rm pref} u$, and 
$u \vee_{\rm pref} v = v$ if $u \le_{\rm pref} v$. Hence in $A^*$, prefix 
codes are the same thing as joinless codes.
This is not the case for $nA^*$ with $n \ge 2$; here, joinless codes are a 
special case of initial factor codes, and the join is characterized as 
follows:

\begin{lem} \label{joinAn} {\bf (join for {\boldmath $\le_{\rm init}$} in 
{\boldmath $nA^*$}).}
 \ For all $\,u = (u_1, \, \ldots, u_n), \ v = (v_1, \, \ldots, v_n)$
$\in nA^*$, the following are equivalent:  

\smallskip

\noindent
{\rm (1)} \ \ the join $\,u \vee v \,$ (with respect to $\, \le_{\rm init}$) 
exists; \\   
{\rm (2)} \ \ $u$ and $v$ have a common upper bound for $\, \le_{\rm init}$,
 \ i.e., $ \, (\exists z) \, [\, u \le_{\rm init} z$  $ \,{\rm and}\, $
$v \le_{\rm init} z \,]$; \\  
{\rm (3)} \ \ for all $i = 1, \, \ldots, n$: 
 \ $u_i \para_{\rm pref} v_i \ $ in $A^*$. 

\medskip

\noindent Moreover, if 
$\, u \vee v = ((u \vee v)_i : i =   1, \, \ldots, n) \,$ exists, 
then 

\smallskip

\hspace{.5in}  {\rm
$ (u \vee v)_i \ = \ \left\{ \begin{array}{ll} 
    u_i & \ \ \ \mbox{if $\, v_i \le_{\rm pref} u_i \,$ (in $A^*$),} \\   
    v_i & \ \ \ \mbox{if $\, u_i \le_{\rm pref} v_i \,$ (in $A^*$).}
\end{array} \right.    $ 
}

\smallskip

\noindent In other words, if $u \vee v$ exists then 
$\,(u \vee v)_i \, = \, \max_{\le_{\rm pref}}\{u_i, v_i\}$, and
$\, |(u \vee v)_i| \, = \, \max\{|u_i|, |v_i|\}$.
\end{lem}
So, in $nA^*$ the relation $\, \|_{\rm init} \,$ is not equivalent to
coordinatewise $\parallel_{\rm pref}$; the latter is equivalent to the
existence of a join; $\, \|_{\rm init} \,$ implies (but is not
equivalent to) existence of a join.

\medskip

\noindent {\sc Proof.} \ [(1) $\Rightarrow$ (2)] is obvious. 
[(2) $\Rightarrow$ (3)] is straightforward: If $u \le_{\rm init} z$ and
$v \le_{\rm init} z$ for some $z \in nA^*$ then $u r = v t = z$ for some
$s,t,z \in nA^*$. Hence, $u_i s_i = v_i t_i = z_i$, so
$\, u_i \parallel_{\rm pref} v_i$ in $A^*$.

\noindent [(3) $\Rightarrow$ (1)] \ Suppose $u_i \para_{\rm pref} v_i$ 
for all $i$.  Then $u_i \le_{\rm pref} v_i$ for some $i$, and 
$v_i \le_{\rm pref} u_i$ for the other $i$.
Hence, $(u \vee v)_i = u_i$ if $v_i \le_{\rm pref} u_i$ (in $A^*$), and 
$(u \vee v)_i = v_i$ otherwise; so $u \vee v$ exists.
 \ \ \ $\Box$

\begin{notat} \label{GenA}
 \ Let \ \ $A_{\varepsilon,n} \ = \ $  $\bigcup_{1 \le i \le n}$
$\{\varepsilon\}^{i-1} \times A \times \{\varepsilon\}^{n-i}$ . 
\end{notat}
Note that $A_{\varepsilon,n}$ is the unique minimum generating set of $nA^*$
as a monoid; the cardinality is $\, |A_{\varepsilon,n}| = n \, |A|$.

\begin{lem} \label{UniquePrefC} \hspace{-.07in}.

\noindent {\rm (1)} \ Every right ideal $R \subseteq nA^*$ is generated, 
as a right ideal, by a unique initial factor code.  
(Finiteness of generating sets is not assumed here.)

\noindent {\rm (2)} \ If a right ideal $R \subseteq nA^*$ is generated by a 
{\em joinless} code then the unique initial factor code that generates $R$ 
is joinless.
\end{lem}
{\sc Proof.} \ Let 
 \ $P \ = \ R \ \minus \ R \cdot A_{\varepsilon, n}$.
We claim that $P$ is an initial factor code that generates $R$, and that 
$P$ is the unique such initial factor code. 
(We closely follow the proof of \cite[Lemma 8.1(1)]{BiThompsMon}.)

Let us show that $P$ generates $R$. Obviously, since $P  \subset R$, we have 
$P \, (nA^*) \subseteq R \, (nA^*) = R$.
Conversely, to show that $R \subseteq P \, (nA^*)$, consider any $r \in R$. 
In $nA^*$, $\, r$ has only finitely many 
initial factors, hence there exists a (not necessarily unique) $p \in R$ 
which is an initial factor of $r$ and is $\le_{\rm init}$-minimal in $R$. 
So $r = px$ for some $x \in nA^*$.  And 
$p \not\in R \, A_{\varepsilon, n}$, otherwise there would exist
$p = r' a$ for some $r' \in R, a \in A_{\varepsilon, n}$, which 
would contradict that $p$ is $\le_{\rm init}$-minimal in $R$.  Hence
$p \in R \minus R \, A_{\varepsilon, n}$. 

To show that $P$ is an initial factor code, let $p, p' \in P$ and suppose 
$p = p' x$ for some $x \in nA^*$. If $x \ne (\varepsilon)^n$ then 
$p \in R A_{\varepsilon, n}$, contradicting the assumption that $p \in P$
($= R \minus R A_{\varepsilon, n}$). So, $p = p'$. 

To prove uniqueness of the initial factor code that generates $R$, we 
generalize the proof of \cite[Lemma 8.1(1')]{BiThompsMon}.
If $P_1 \, (nA^*) = P_2 \, (nA^*)$ for two initial factor codes $P_1, P_2$, 
then for every $p_1 \in P_1$ there exists $p_2 \in P_2$ such that 
$p_1 = p_2 x$ (for some $x \in nA^*$). Also, there is $p_1' \in P_1$ such 
that $p_2 = p_1' y$ (for some $y \in nA^*$). Hence $p_1 = p_1' xy$, which 
implies $x = y = (\varepsilon)^n$, since $P_1$ is an initial factor code. 
Thus, $p_1 = p_2 \in P_2$. Therefore, $P_1 \subseteq P_2$. Similarly we have 
$P_2 \subseteq P_1$, so $P_1 = P_2$.

\smallskip

Part (2) follows immediately from the uniqueness of the initial factor code 
that generates $R$.
 \ \ \ $\Box$

\begin{lem} \label{OmegaJoilessPref}
 \ Let $P \subset nA^*$ be a finite maximal joinless code. Then every 
$w \in nA^{\omega}$ has a unique initial factor in $P$; \ i.e., 
$\, (\forall w \in nA^{\omega})(\exists! \,p \in P,$ 
$u \in nA^{\omega}) \, [\, w = p u \,]$.
\end{lem}
{\sc Proof.} If there were two different initial factors $p, q$ of $w$ in 
$P$ then $p$ and $q$ would be initial factors of a finite initial factor of
$w$; hence $p$ and $q$ would have a join, contradicting that $P$ is joinless.
This shows uniqueness.

Let us show existence.
Since $P$ is a maximal joinless code, every initial factor $v$ of $w$ has 
a join with some element of $P$. Let us pick $v$ so that its coordinates 
(in $A^*$) are longer than all the coordinates of the elements of $P$. Then 
the element of $P$ that has a join with $v$ is an initial factor of $v$.
 \ \ \ $\Box$

\begin{lem} \label{JoinlessEssential}
 \ Let $P \subset nA^*$ be any finite {\em joinless} code, and let 
$\,R = P \cdot (nA^*)\,$ be the right ideal generated. (Recall that by Lemma
\ref{UniquePrefC}, $\, P$ is uniquely determined by $R$.) 
Then the following are equivalent:

\smallskip

\noindent {\rm (1)} \ $R$ is an {\em essential} right ideal;

\noindent {\rm (2)} \ $P$ is {\em maximal} as a joinless code;

\noindent {\rm (3)} \ $P \cdot (nA^{\omega}) \, = \, nA^{\omega}$;

\noindent {\rm (4)} \ $R \cdot (nA^{\omega}) \, = \, nA^{\omega}$.
\end{lem}
{\sc Proof.} \ $[(1) \Leftrightarrow (2)]$ \ Suppose $P$ is a finite joinless
code. Then $P$ is maximal joinless iff
every $v \in nA^*$ has a join with some element of $P$ (as follows
directly from the definition of maximality). This is equivalent to the
property that every monogenic right-ideal of $nA^*$ intersects
$P \, (nA^*)$; i.e., $P \, (nA^*)$ is essential.

$[(3) \Rightarrow (1)]$ \ If $\, P \, (nA^{\omega}) = nA^{\omega}$, then
every $w \in nA^{\omega}$ has an initial factor in $P$. It follows that for 
every right ideal $R \subset nA^*$,
 \ $R \, (nA^{\omega}) \subseteq P \, (nA^{\omega})$.
Hence $R$ intersects $P \, (nA^*)$.
So, $P \, (nA^*)$ is essential.

$[(2) \Rightarrow (3)]$ \ Suppose $P$ is a finite maximal joinless code.
Let $w \in nA^{\omega}$, and for any
$(i_1, \, \ldots, i_n) \in {\mathbb N}^n$, let $w^{(i_1, \, \ldots, i_n)}$ 
be the initial factor of $w$ in 
$A^{i_1} \times \, \, \ldots \, \times A^{i_n}$.
Then $w^{(i_1, \, \ldots, i_n)}$ has a join with some $p \in P$.
Since $P$ is finite, $p$ is an initial factor of $w^{(i_1, \, \ldots, i_n)}$ 
if each of $i_1, \, \ldots, i_n$ is larger than
$\max\{|p_i| : p \in P, \ i \in \{1, \, \ldots,n\} \}$.
Hence, $p$ is an initial factor of $w$, so $w \in P \, (nA^{\omega})$.
Since for every $w \in nA^{\omega}$ such a $p \in P$ exists (by Lemma
\ref{OmegaJoilessPref}), we conclude that 
$nA^{\omega} \subseteq P \ (nA^{\omega})$.

The equivalence of (3) and (4) is obvious since
$\, nA^* \cdot nA^{\omega} = nA^{\omega}$, so $\, R \cdot (nA^{\omega})$ 
$=$  $P \cdot (nA^*) \cdot (nA^{\omega}) = P \cdot (nA^{\omega})$.
 \ \ \ $\Box$

\bigskip

\noindent {\bf Remark.} Lemma \ref{JoinlessEssential} only talks about 
joinless generated right ideals. 
Indeed, an essential finitely generated right ideal in $nA^*$ is not
necessarily joinless generated. An example for $\,A = \{0,1\}\,$ is

\smallskip

 \ \ \ $R$ \ $=$
 \ $\{ (\varepsilon, 0), (0, \varepsilon), (1,1)\} \cdot (2A^*)$.

\smallskip

\noindent It is easy to prove that $R$ is essential, and that 
$\{(\varepsilon, 0),$ $(0, \varepsilon),$  $(1,1)\}$ is an initial-factor 
code that is not joinless (since 
$(\varepsilon, 0) \vee (0, \varepsilon) = (0,0)$ exists).
By Lemma \ref{UniquePrefC}, this initial factor code is unique, i.e., $R$ is 
not generated by any other initial-factor code; hence $R$ is not joinless 
generated.

Section 5 of version 1 of \cite{BinG} gives a detailed proof (independently 
of Lemma \ref{UniquePrefC}) that $R$ is essential in $2 A^*$, and that $R$ 
is not generated (as a right ideal) by any finite joinless code.

\bigskip

\noindent {\bf DAGs and} {\boldmath $nA^*$}{\bf :} 
The following generalizes the well known concepts of the tree of $A^*$ and 
the tree of a prefix code.  
We abbreviate {\em directed acyclic graph} by {\sc dag}.  A few definitions:
The {\em leaves} of a {\sc dag} are the vertices of out-degree 0; all the
other vertices are {\em interior vertices}. For a {\sc dag} $D$, the 
sub-{\sc dag} spanned by the interior vertices of $D$ is called the 
{\em interior} {\sc dag} of $D$.   The {\em sources} of a {\sc dag} are the 
vertices of in-degree 0; if there is only one source, and all vertices are 
reachable from this source, this source is called the {\em root}, and the 
{\sc dag} is then called {\em rooted}.  The 
{\em depth} of a vertex $v$ in a rooted {\sc dag} is defined to be the length 
of the shortest path from the root to $v$; by ``path'' we will always mean a
directed path.

\medskip

\noindent $\bullet$ \ {\em The {\sc dag} of $nA^*$} is the infinite rooted 
{\sc dag} with vertex set $nA^*$ and root $(\varepsilon)^n$; the edges are
the ordered pairs $(s,t) \in (nA^*) \times (nA^*)$ such that there exists 
$i \in \{1, \, \ldots, n\}$ and $a \in A$ with
$t = (s_1, \, \ldots, s_{i-1}, \, s_i a, \, s_{i+1}, \, \ldots, s_n) \,$ 
(where $\,s = (s_1,\, \ldots, s_{i-1},\, s_i,\, s_{i+1},\, \ldots, s_n)$).
Hence every vertex has $|A_{\varepsilon,n}|$ ($ = n \cdot |A|$) children;
see Notation \ref{GenA}.  And $u \le_{\rm init} v$ iff there exists a 
directed path from $u$ to $v$ in the {\sc dag}.
It is easy to show that the depth of a vertex 
$\,v = (v_1,\, \ldots\, , v_n)$ in the {\sc dag} of $nA^*$ is 
$ \ \sum_{i=1}^n |v_i|$.

\smallskip

The {\sc dag} of $nA^*$ is the {\em right Cayley graph} of the monoid 
$nA^*$ over the generating set $A_{\varepsilon, n}$.

\medskip

\noindent $\bullet$ \ For any finite subset $P \subset nA^*$ we define 
the {\em initial factor {\sc dag} of $P$} (also called the $P$-{\sc dag}): 
This is a finite rooted sub{\sc dag} of the {\sc dag} of $nA^*$;
the root of the $P$-{\sc dag} is the root of the {\sc dag} of $nA^*$;
the vertices and edges are those vertices, respectively edges, of the 
{\sc dag} of $nA^*$ that appear on any path from the root to any vertex 
in $P$. Hence the vertices of the $P$-{\sc dag} are all the initial factors
of the elements of $P$ (so the $P$-{\sc dag} is uniquely determined by $P$).
The set of leaves of the $P$-{\sc dag} is $P$ iff $P$ is an initial factor 
code.

\medskip

Note that the trees and {\sc dag}s considered here are not ordered trees or
{\sc dag}s; i.e., the children of a vertex are defined as a set, not a 
sequence; similarly, the leaves form a set, not a sequence.

\begin{lem} \label{SaturatedLeaf}
 \ Let $P \subset nA^*$ be a finite maximal joinless code such that
$P \ne \{\varepsilon\}^n$.  
Let $v = (v_1, \, \ldots, v_n)$ be any leaf of the interior {\sc dag} of the 
{\sc dag} of $P$, and let $v_+$ be the set of children of $v$ in the 
$P$-{\sc dag}; so \ $v_+ = v \cdot A_{\varepsilon,n} \,\cap\, P$, and  
$v_+$ is non-empty (since $v$ is an interior vertex).

\smallskip

\noindent {\rm (0)} \ Then $v_+$ satisfies

\medskip

 \ \ \  \ \ \ $v_+ \ \subseteq \ $
$\{(v_1, \, \ldots, v_{i-1}, \, v_i a, \, v_{i+1}, \, \ldots, v_n)$
$ \, : \, $ $a \in A\}$   $ \ = \ $
$v \cdot (\{\varepsilon\}^{i-1} \times A \times \{\varepsilon\}^{n-i})\,$, 

\medskip

\noindent for some $\, i \in \{1, \, \ldots, n\}$; and $i$ is unique 
(for a given $v$ and $P$). 

\smallskip

\noindent {\rm (1)} For $n = 1$, part {\rm (0)} holds with equality for 
every leaf $v$ of the interior {\sc dag}:
 \ $v_+ \ = \ \{v a : a \in A\}$.

\smallskip

\noindent {\rm (2) (Lawson and Vdovina \cite[Thm.\ 12.11]{LawsonVdovina}, 
but with a different formalism.)}  
 \ For $n = 2$ and $|A| = 2$, part {\rm (0)} holds with 
equality for some maximum-depth leaf $v$ of the interior {\sc dag} of $P$: 

\smallskip

 \ \ \  \ \ \ $v_+ \ = \ \{(v_1 a, \, v_2) : a \in A\}$ \ \ \ or 
 \ \ \ $v_+ \ = \ \{(v_1, \, v_2 a) : a \in A\}$.   

\smallskip

\noindent However, equality does not necessarily hold for every interior 
leaf, not even for every interior leaf of maximum depth.  

\smallskip

\noindent {\rm (3) (Lawson and Vdovina \cite[Ex.\ 12.8]{LawsonVdovina})} 
 \ For $n \ge 3$ there exist finite maximal joinless codes 
$P \subset n \, \{0,1\}^*$ for which the inclusion in part {\rm (0)} is 
strict. I.e., for every leaf $v$ of the interior {\sc dag} and for every 
$i \in \{1, \, \ldots, n\}$:

\smallskip

 \ \ \  \ \ \ $v_+ \ \ne \ $
$\{(v_1, \, \ldots, v_{i-1}, \, v_i a, \, v_{i+1}, \, \ldots, v_n)$
$ \, : \, $ $a \in A\}$.
\end{lem}
{\sc Proof.}
{\bf (0)} Since $v$ is interior without having interior children, it 
contains a least one child in $P$, of the form 
$(v_1, \, \ldots, v_{i-1}, \, v_i a, \, v_{i+1}, \, \ldots, v_n)$, for some
$a \in A$, $i \in \{1, \, \ldots, n\}$.  

Any possible child of $v$ belongs to $v \cdot A_{\varepsilon, n}$.  If, in 
addition to 
$(v_1, \, \ldots, v_{i-1}, \, v_i a, \, v_{i+1}, \, \ldots, v_n)$, $v$ had 
an additional child of the form 
$(v_1,$  $\ldots,$  $v_{j-1},$  $\, v_j b,$   $\, v_{j+1},$   $\ldots,$ 
$v_n)$ with $i \ne j$
(for any $b \in A$), then $P$ would not be joinless. Indeed, these two 
children have the join 
$(v_1, \, \ldots, v_j b, \, \ldots, v_i a, \, \ldots, v_n)$ (if $j < i$), or
$(v_1, \, \ldots, v_i a, \, \ldots, v_j b, \, \ldots,  v_n)$ (if $i < j$).
This shows that all children of $v$ belong to 
$\, \{(v_1, \, \ldots, v_{i-1}, \, v_i a, \, v_{i+1}, \, \ldots, v_n) :$
$a \in A\} \, $   for one particular $i$ (depending on $v$).  

\smallskip

\noindent {\bf (1)} For $n=1$ the Lemma is folklore knowledge. 

\smallskip

\noindent {\bf (2)} (This result is equivalent to 
\cite[Thm.\ 12.11]{LawsonVdovina}, but the proof given here is rather 
different.) 
 
Here $A = \{0,1\}$.
Let $v = (v_1,v_2)$ be a maximum-depth leaf of the interior {\sc dag} 
of $P$. Since $v$ is an interior leaf, at least one of its children is in 
$P$. Hence either $(v_1 a, v_2) \in P$ or $(v_1, v_2 a) \in P$, for some 
$a \in A$. 

Let us assume that $a = 0$ and that $(v_1 0, v_2) \in P$; the other cases are 
very similar. Since $(v_1,v_2)$ has maximum depth in the interior {\sc dag},
$(v_1 0, v_2)$ has maximum depth in $P$.  

If it is also the case that $(v_1 1, v_2) \in P$, then 
$\{(v_1 0, v_2), (v_1 1, v_2)\} \subseteq P$, and the Lemma holds. 
Therefore, from here on we only consider the situation where 
$(v_1 1, v_2) \not\in P$ (but $(v_1 0, v_2) \in P$). 
Then there exists $(u_1,u_2) \in P \minus \{(v_1 0, v_2)\}$ with
$(u_1,u_2) \ne (v_1 1, v_2)$, 
such that $(u_1,u_2)$ has a join with $(v_1 1, v_2)$.
By Prop.\ \ref{joinAn}, this is equivalent to
$\, u_1 \,\|_{\rm pref}\, v_1 1 \,$ and $\, u_2 \,\|_{\rm pref}\, v_2$.
 
This leads to four cases.

\medskip

\noindent Case 1: \ $v_1 1 \le_{\rm pref} u_1$ and $v_2 \le_{\rm pref} u_2$.

\smallskip

Then $v_1 <_{\rm pref} u_1$ and $v_2 \le_{\rm pref} u_2$. Since $(v_1,v_2)$
is a leaf of the interior {\sc dag} of $P$, and $(u_1,u_2) \in P$, it 
follows that $(u_1,u_2)$ is a child of $(v_1,v_2)$. Since $(u_1, u_2)$
$\not\in$ $\{(v_1 0, v_2), (v_1 1, v_2)\}$, it follows that $(u_1, u_2)$ is 
of the form $(v_1, v_2 c)$ for some $c \in A$. But then $(u_1, u_2)$ 
($ = (v_1, v_2 c)$) has a join with $(v_1 0, v_2) \in P$, contradicting the 
fact that $P$ is joinless.  So, case 1 is ruled out.

\medskip

\noindent Case 2: \ $v_1 1 \ge_{\rm pref} u_1$ and $v_2 \ge_{\rm pref} u_2$;
 \ since $(u_1,u_2) \ne (v_1 1, v_2)$, at least one of these
$\ge_{\rm pref}$ is strict.

\smallskip

\noindent Case 2.1: \ $v_1 1 >_{\rm pref} u_1$ and $v_2 \ge_{\rm pref} u_2$:
 
Then $(u_1, u_2)$ is interior, since $(v_1, v_2)$ is interior. But 
$(u_1, u_2)$ being an interior vertex contradicts the assumption that 
$(u_1, u_2) \in P$. So case 2.1 is ruled out. 

\smallskip

\noindent Case 2.2: \ $v_1 1 = u_1$ and $v_2 >_{\rm pref} u_2$:

Then $u_2 = v_2 c z$, for some $c \in A$ and $z \in A^*$.
But now $(u_1, u_2) = (v_1 1, v_2 cz)$ has greater depth than $(v_1 0, v_2)$,
which has maximum depth in $P$. So case 2.2 is ruled out.

\smallskip

\noindent Case 3: \ $v_1 1 \ge_{\rm pref} u_1$ and $v_2 \le_{\rm pref} u_2$;
 \ since $(u_1,u_2) \ne (v_1 1, v_2)$, at least one of $\ge_{\rm pref}$ or
$\le_{\rm pref}$ is strict.

\smallskip

\noindent Case 3.1: $v_1 1 > _{\rm pref} u_1$ and $v_2 \le_{\rm pref} u_2$.

Then $v_1 1 = u_1 x 1$, and $u_2 = v_2 y$ for some $x, y \in A^*$; so 
$v_1 = u_1 x$.  But then 
$\, (u_1, u_2) \vee (v_1 0, v_2)$ $=$ $(u_1, v_2 y) \vee (u_1 x 0, v_2)$ $=$  
$(u_1 x 0, v_2 y) \,$ exists, contradicting the fact that 
$\{(u_1, u_2), (v_1 0, v_2)\}$  $\subseteq P$.
So case 3.1 is ruled out.

\smallskip

\noindent Case 3.2: $v_1 1 =  u_1$ and $v_2 <_{\rm pref} u_2$.

Then $(u_1, u_2)$ has greater depth than $(v_1 0, v_2)$, contradicting the 
fact that $(v_1 0, v_2)$ has maximum depth in $P$. 
So case 3.2 is ruled out.
 
\smallskip

\noindent Case 4: \ $v_1 1 \le_{\rm pref} u_1$ and $v_2 \ge_{\rm pref} u_2$.

Then $u_1 = v_1 1 x$ and $v_2 = u_2 y$ for some $x, y \in A^*$. 
Since $(v_1 0, v_2)$ has maximum depth in $P$ we have 
$|u_1| + |u_2| \le |v_1 0| + |v_2|$, hence
$|v_1| + 1 + |x| + |u_2| \le |v_1| + 1 + |u_2| + |y|$, hence
$|x| \le |y|$. Moreover, $y \ne \varepsilon$, otherwise $|x| = 0$, hence
$x = \varepsilon$, hence $(u_1,u_2) = (v_1 1, v_2)$, which would imply 
$(v_1 1, v_2) \in P$.
In summary this proves: 

\smallskip

 \ \ \  \ \ \ $|x| \le |y| \ne 0$ \ and \ $v_2 >_{\rm pref} u_2$.

\smallskip

\noindent Notation (used in the remainder of the proof):  
For any $z \in \{0,1\}^+$, let $z^-$ denote the bitstring obtained by 
complementing the right-most bit of $z$. And $z \{0,1\}^{-1}$ denotes the 
bitstring obtained by removing the right-most bit of $z$.

\smallskip

Note that since $(v_1 0, v_2) \in P$, if we prove that 
$(v_1 0, v_2^-) \in P\,$ then the Lemma holds for the interior vertex 
$(v_1 0, \, v_2 \{0,1\}^{-1})$.

\smallskip

\noindent Claim: \ $(v_1 0, v_2^-) \in P$.
 
\smallskip

\noindent Proof of the Claim: Assume by contradiction that there exists 
$(w_1,w_2) \in P$ such that $(w_1,w_2) \ne (v_1 0, v_2^-)$, and $(w_1,w_2)$ 
has a join with $(v_1 0, v_2^-)$.
The existence of this join is equivalent to $w_1 \,\|_{\rm pref}\, v_1 0$
and $w_2 \,\|_{\rm pref}\, v_2^-$. 

This leads to four cases. 

\smallskip

\noindent Case 4.1: 
$w_1 \le_{\rm pref} v_1 0$ and $w_2 \le_{\rm pref} v_2^-$. 
At least one of the $\le_{\rm pref}$ is strict.

\noindent Case 4.1.1: $w_1 \le_{\rm pref} v_1 0$ and 
$w_2 \le_{\rm pref} v_2^- A^{-1} = v_2 A^{-1}$.

Then the join $(w_1,w_2) \vee (v_1 0, v_2) = (v_1 0, v_2)$ exists, 
contradicting the fact that $(w_1,w_2)$ and $(v_1 0, v_2)$ belong to $P$.
So case 4.1.1 is ruled out.

\noindent Case 4.1.2: $w_1 \le_{\rm pref} v_1$ and
$w_2 \le_{\rm pref} v_2^-$.

Then $v_1 = w_1 \alpha$ and $v_2^- = u_2 y^- = w_2 \beta$ for some
$\alpha, \beta \in A^*$. The latter equality implies that 
$w_2 \,\|_{\rm pref}\, u_2$. 
Recall that in case 4, $u_1 = v_1 1 x$; this and $v_1 = w_1 \alpha$ imply 
that $u_1 = w_1 \alpha 1 x$, hence $u_1 \,\|_{\rm pref}\, w_1$.
Now, since $u_1 \,\|_{\rm pref}\, w_1$ and $w_2 \,\|_{\rm pref}\, u_2$, the
join $(w_1,w_2) \vee (u_1,u_2)$ exists, which contradicts the fact that 
$(w_1,w_2)$ and $(u_1,u_2)$ belong to $P$.
So case 4.1.2 is ruled out.

\smallskip

\noindent Case 4.2:
$w_1 \ge_{\rm pref} v_1 0$ and $w_2 \ge_{\rm pref} v_2^-$.

Since $(v_1 0,v_2)$ has maximum depth in $P$, and $(v_1 0,v_2^-)$ has the
same depth, it follows that $(w_1,w_2) = (v_1 0,v_2^-)$. This contradicts 
the assumption $(w_1,w_2) \ne (v_1 0,v_2^-)$. So case 4.2 is ruled out.

\smallskip

\noindent Case 4.3: 
$w_1 \le_{\rm pref} v_1 0$ and $w_2 \ge_{\rm pref} v_2^-$. 

\noindent Case 4.3.1: 
$w_1 = v_1 0$ and $w_2 >_{\rm pref} v_2^-$ \ (since 
$(w_1, w_2) \ne (v_1 0, v_2^-)$, equality in the first coordinate implies 
strictness in the second).

Then $|w_1| + |w_2| > |v_1 0| + |v_2^-| = |v_1 0| + |v_2|$, i.e., 
$(w_1,w_2)$ has greater depth than $(v_1 0,v_2)$, which contradicts the fact
that $(v_1 0,v_2)$ has maximum depth in $P$. So case 4.3.1 is ruled out.

\noindent Case 4.3.2:
$w_1 <_{\rm pref} v_1 0$ and $w_2 \ge_{\rm pref} v_2^-$.

Then $w_1 \le_{\rm pref} v_1 = w_1 s$, and $w_2 = v_2^- t = u_2 y^- t$, for
some $s, t \in A^*$. Recall that $y \ne \varepsilon$ in case 4. 
Then $(w_1,w_2) \vee (u_1,u_2) = (w_1, u_2 y^- t) \vee (w_1 s, u_2) = $
$(w_1 s, u_2 y^- t)$ exists. 
This contradicts the fact that $(w_1,w_2)$ and $(u_1,u_2)$ belong to $P$.
So case 4.3.2 is ruled out.

\smallskip

\noindent Case 4.4: 
$w_1 \ge_{\rm pref} v_1 0$ and $w_2 \le_{\rm pref} v_2^-$; \ since 
$(w_1,w_2) \ne (v_1 0, v_2^-)$, $\le_{\rm pref}$ or $\ge_{\rm pref}$ is 
strict.  

\noindent Case 4.4.1: 
$w_1 >_{\rm pref} v_1 0$ and $w_2 = v_2^-$.

Then $|w_1| + |w_2| > |v_1 0| + |v_2^-| = |v_1 0| + |v_2|$, hence 
$(w_1,w_2)$ has greater depth than $(v_1 0,v_2)$, which contradicts the fact
that $(v_1 0,v_2)$ has maximum depth in $P$. So case 4.4.1 is ruled out.  

\noindent Case 4.4.2: 
$w_1 \ge_{\rm pref} v_1 0$ and $w_2 <_{\rm pref} v_2^-$.

Then $w_1 = v_1 0 s$; also, $w_2 <_{\rm pref} v_2$ (since $v_2$ and $v_2^-$
only differ in the right-most bit), so $v_2 = w_2 t$, for some $s, t \in A^*$.
Now, $(w_1,w_2) \vee (v_1 0, v_2) = (v_1 0 s, w_2) \vee (v_1 0, w_2 t) =$
$(v_1 0 s, w_2 t)$ exists. This contradicts the fact that $(w_1,w_2)$ and
$(v_1 0, v_2)$ belong to $P$. So case 4.4.2 is ruled out.

Since we now ruled out all sub-cases of case 4, this completes the proof (by 
contradiction) of the Claim.

\medskip

Summary of the proof so far: We have $(v_1 0, v_2) \in P$ for some 
maximum-depth vertex $(v_1, v_2)$ in the interior of the $P$-{\sc dag}. 
(The cases where, instead, we have $(v_1 1, v_2)$ or $(v_1, v_2 0)$, or 
$(v_1,v_2 1)$ in $P$, are similar.) 

If we also have $(v_1 1, v_2) \in P$ then the Lemma holds.

If $(v_1 1, v_2) \not\in P$ then there exists $(u_1,u_2) \in P$ that has a 
join with $(v_1 1, v_2)$. Four cases are possible, of which cases 1, 2, and 
3 were ruled out. In case 4 we showed that $(v_1 0, v_2^-) \in P$; hence   
in case 4, $(v_1 0, v_2)$ and $(v_1 0, v_2^-)$ belong to $P$, i.e.,
the Lemma holds for the interior vertex $(v_1 0, v_2 \{0,1\}^{-1})$.

\medskip

The following is an example where not every maximum-depth interior leaf has 
two children in $P$. Consider the maximal joinless code $P = $
$\{(0,0), (0,1), (1,\varepsilon)\}$. Here the interior leaf 
$v = (\varepsilon,0)$ has maximum depth, and has only one child  in $P$
(namely $(0,0)$).  
Nevertheless, there is another maximum-depth interior leaf, namely
$(0,\varepsilon)$, that has two children in $P$ (namely $(0,0)$ and $(0,1)$).  

\medskip

\noindent {\bf (3)} Example (from \cite[Ex.\ 12.8]{LawsonVdovina}): Let 
 \ $P = $ 
$\{(0,0,\varepsilon), \, (1,\varepsilon,0), \, (\varepsilon,1,1),$ 
$(0,1,0), \, (1,0,1)\} \, \subset \, 3 \, \{0,1\}^*$.  
It is easy to verify that $P$ is a finite maximal joinless code, and that 
no leaf of the interior {\sc dag} has two children in $P$.  
 \ \ \ $\Box$

\bigskip

\noindent {\bf Remark about Lemma \ref{SaturatedLeaf}:} 
Version 1 of this paper (see \cite{BinG}) stated incorrectly that ``for 
every $n \ge 1$ and  every leaf $v$ of the interior {\sc dag} of $P$: 
$ \ v_+ =$ 
$v \cdot (\{\varepsilon\}^i \times A \times \{\varepsilon\}^{n-i-1})\,$
(for some $i$, $0 \le i < n$)''. This statement had to be modified for $n=2$ 
(from ``for every leaf'' to ``there exists a leaf''), and dropped for 
$n \ge 3$.  The above counter-example for $n \ge 3$ was given in 
\cite{Lawson2019} and \cite[Ex.\ 12.8]{LawsonVdovina}.

\begin{lem} \label{PjoinlessPA}
 \ Let $P \subset nA^*$ be a finite set. For any
$\, p = (p_1, \, \ldots, p_n) \in P$ and $i \in \{1, \, \ldots, n\}$, let

\medskip

 \ \ \  \ \ \ $P_{p,i}' \ = \ (P \minus \{p\})$  $ \ \cup \ $
$\{(p_1, \, \ldots, p_{i-1}, \, p_i a, \, p_{i+1}, \, \ldots, p_n) \, : \, $
$a \in A\}$.

\medskip

\noindent Then we have:

\smallskip

\noindent {\bf (1)} \ $P$ is joinless iff $\, P_{p,i}' \,$ is joinless.

\smallskip

\noindent {\bf (2)} \ $P$ is a maximal joinless code iff $\, P_{p,i}' \,$ is
a maximal joinless code.

\medskip

\noindent The set $P_{p,i}'$ is called a {\em one-step restriction} of $P$
(``restriction'' because 
$\, P_{p,i}' \cdot (nA^*) \,\subsetneqq\, P \cdot (nA^*)$); 
and $P$ is called a {\em one-step extension} of $P_{p,i}'$.
Clearly, $ \ |P_{p,i}'| = |P| - (|A|-1)$. 
\end{lem}
{\sc Proof.} {\bf (1)} $[\Rightarrow]$  Let us assume that $P$ is joinless.  
For any $a, a' \in A$ with $a \ne a'$, the join of 
$\,(p_1,$  $\ldots,$  $p_{i-1},$ $p_i a,$ $p_{i+1},$  $\ldots,$ $p_n)\,$
and $\,(p_1,$ $\ldots,$ $p_{i-1},$ $p_i a',$  $p_{i+1},$ $\ldots$, $p_n)\,$
does not exist, since $p_i a$ and $p_i a'$ are not prefix-comparable.

If $q \in P \minus \{p\}$ and 
$(p_1, \, \ldots, p_{i-1}, p_i a, p_{i+1}, \, \ldots, p_n)$ were both
initial factors of some $z \in nA^*$, then $q$ and $p$ would also both be
initial factors of $z$, contradicting the assumption that $P$ is joinless.

Finally, all pairs $q_1, q_2 \in P \minus \{p\}$ 
($\subset P_{p,i}'$) are joinless since $P$ is joinless.
Thus $P_{p,i}'$ is joinless.

\smallskip

\noindent $[\Leftarrow]$ Let us assume that $P_{p,i}'$ is joinless.  
Then every pair $q_1, q_2 \in P \minus \{p\}$ ($\subset P_{p,i}'$) 
is joinless.

If $q \in P \minus \{p\}$ and $p$ had a join $z$, then both $p$ and 
$q$ would be initial factors of $z$. 
By Lemma \ref{joinAn}, $z_j = \max\{q_j, p_j\}$ for 
all $j \in \{1, \, \ldots, n\}$. We have two cases.

\smallskip

\noindent Case 1: $z_i = p_i$ \ (for the $i$ used in $P_{p,i}'$).

This is equivalent to $q_i$ being a prefix of $p_i$. Then $q_i$ is 
a prefix of $p_i a$ as well (for every $a \in A$), hence 
$q \in P \minus \{p\}$ and 
$(p_1, \, \ldots, p_{i-1}, p_i a, p_{i+1}, \, \ldots, p_n)$ have a join. 
But this contradicts the assumption that $P_{p,i}'$ is joinless.

\smallskip

\noindent Case 2: $z_i \ne p_i$, and $z_i = q_i$ \ (for the $i$ used in 
$P_{p,i}'$).

Then $p_i$ is a strict prefix of $q_i$ ($= z_i$), hence $p_i a$ is a prefix
of $q_i$ for some $a \in A$. It follows that $z$ has $q$ 
and $(p_1, \, \ldots, p_{i-1}, p_i a, p_{i+1}, \, \ldots, p_n)$ as initial
factors; this contradicts the assumption that $P_{p,i}'$ is joinless.

\smallskip

\noindent {\bf (2)} $[\Rightarrow]$ Suppose $P$ is a maximal joinless code. 
Hence, every $x \in nA^*$ has a join with some $q \in P$ (otherwise $x$ 
could be added to $P$, which would contradict that $P$ is maximal joinless).
We want to show that $x$ also has a join with some element of $P_{p,i}'$. 

If $q \ne p$ then $q \in P_{p,i}'$, hence $x$ also has a join with some 
$q \in P_{p,i}'$.

If $q = p$, i.e., $z = x \vee p$, then $z_j = \max\{x_j,\, p_j\}$ for all
$j \in \{1, \, \ldots, n\}$. 
We have two cases.

\smallskip

\noindent Case 1: $z_i = p_i$ \ (for the $i$ used in $P_{p,i}'$).

This is equivalent to $x_i$ being a prefix of $p_i$. Then $x_i$ is a prefix
of $p_i a$ too (for every $a \in A$), hence $x$ and 
$(p_1, \, \ldots, p_{i-1}, p_i a, p_{i+1}, \, \ldots, p_n)$ have a join.
So, $x$ has a join with some element of $P_{p,i}'$.

\smallskip

\noindent Case 2: $z_i \ne p_i$, and $z_i = x_i$ \ (for the $i$ used in
$P_{p,i}'$).

Then $p_i$ is a strict prefix of $x_i$ ($= z_i$), hence $p_i a$ is a
prefix of $x_i$ for some $a \in A$. It follows that $z$ has $x$ and 
$(p_1, \, \ldots, p_{i-1}, p_i a, p_{i+1}, \, \ldots, p_n)$ as initial
factors; this implies that $x$ has a join with 
$(p_1, \, \ldots, p_{i-1}, p_i a, p_{i+1}, \, \ldots, p_n)$ $\in P_{p,i}'$
(for this particular $a \in A$).
 
\smallskip

\noindent $[\Leftarrow]$ Suppose that $P_{p,i}'$ is maximal joinless. Then
every $x \in nA^*$ has a join with some $q \in P_{p,i}'$.
We want to show that $x$ also has a join with some element of $P$.

If $q \ne (p_1, \, \ldots, p_{i-1}, p_i a, p_{i+1}, \, \ldots, p_n)$ for all
$a \in A$, then $q \in P$ so $x$ also has a join with $q \in P$.

If $q = (p_1, \, \ldots, p_{i-1}, p_i a, p_{i+1}, \, \ldots, p_n)$ for some 
$a \in A$, then let $z$ be the join of $x$ and
$(p_1,$  $\ldots,$  $p_{i-1},$  $p_i a,$  $p_{i+1},$  $\ldots, p_n)$.
Then $z$ has $x$ and 
$(p_1, \, \ldots, p_{i-1}, p_i a, p_{i+1}, \, \ldots, p_n)$ 
as initial factors, hence $p$ is an initial factor of $z$. Hence $x \vee p$ 
exists, so $x$ has a join with an element of $P$.
 \ \ \ $\Box$

\bigskip

The properties of joinless codes given in Lemma \ref{PjoinlessPA} do not 
hold for initial factor codes in general.
For example, for $A = \{0,1\}$ consider the initial factor code 
$\, P = \{(\varepsilon, 0), \, (0, \varepsilon)\}$.
Then for $p = (0, \varepsilon)$ and $i = 2$ we obtain
$\, P_{p,i}' = \{(\varepsilon, 0), \, (0, 0), \, (0,1)\}$, which is not an
initial factor code.

\medskip

The process of one-step restriction or extension can be iterated, which 
inspires the following definition and the algorithm.

\begin{defn} \label{DEFtreeofP} {\bf (parse trees).} 
 \ Let $P \subset nA^*$ be a finite joinless code. A {\em parse tree of} $P$
is any subtree $T$ of the {\sc dag} of $P$ with the following properties:

\smallskip

\noindent {\rm (1)} The root of $T$ is $(\varepsilon)^n$ (i.e., the root 
of the {\sc dag} of $P$); and the set of leaves of $T$ is $P$ (i.e., the 
leaves of the {\sc dag} of $P$).

\noindent {\rm (2)} For every interior vertex $v$ of $T$ the set of children 
in $T$ is 
$ \ v \cdot (\{\varepsilon\}^{i-1} \times A \times \{\varepsilon\}^{n-i})$,
for a unique $i \in \{1, \, \ldots, n\}$.  \ So $v$ has exactly $|A|$ 
children in $T$.
\end{defn}

Given the {\sc dag} of $P$ and a subtree $T$, it is easy to check whether $T$
is a parse tree of $P$; one just needs to check that $(\varepsilon)^n$ 
occurs in $T$, and that every vertex in $T$ is reachable from 
$(\varepsilon)^n$; moreover, for each vertex $v$ of $T$ one checks whether 
it is in $P$, or whether its set of children is of the form 
$\, v \cdot (\{\varepsilon\}^{i-1} \times A \times \{\varepsilon\}^{n-i})$.
Recall the {\sc dag}s and trees are not oriented (children and leaves are 
not ordered).

A maximal joinless code $P$ can have more than one parse tree. E.g., the 
joinless set $ \ \{(0,0),$  $(0,1),$  $(1,0),$   $(1,1)\} \ $ has the 
following two parse trees:

\bigskip

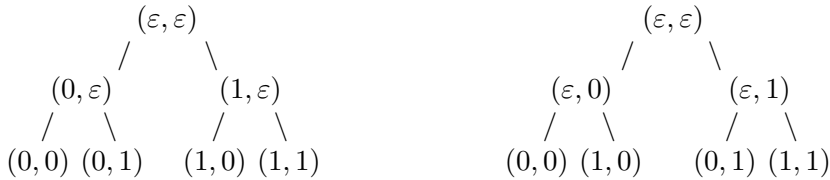
\begin{minipage}{\textwidth}

\hspace{1.05in} $(\varepsilon,\varepsilon)$ \hspace{2.22in}
$(\varepsilon,\varepsilon)$

\hspace{.95in} $/$ \hspace{0.28in} $\setminus$
\hspace{2.in} $/$ \hspace{0.28in} $\setminus$

\hspace{.6in} $(0,\varepsilon)$ \hspace{.45in} $(1,\varepsilon)$ 
\hspace{1.3in} $(\varepsilon,0)$ \hspace{.5in} $(\varepsilon,1)$

\hspace{.55in} $/$ \hspace{0.15in} $\setminus$
\hspace{.39in} $/$ \hspace{0.15in} $\setminus$
\hspace{1.23in} $/$ \hspace{0.13in} $\setminus$
\hspace{.47in} $/$ \hspace{0.13in} $\setminus$

\hspace{.37in} $(0,0)$ $(0,1)$ \hspace{.1in} $(1,0)$ $(1,1)$ 
\hspace{.86in} $(0,0)$ $(1,0)$ \hspace{.16in} $(0,1)$ $(1,1)$
\end{minipage}

\bigskip

\bigskip

\noindent Burillo and Cleary \cite{BurCleary} give a similar tree 
description of tilings of $[0,1]^2$, and point out that the tree is not 
unique.

\smallskip

If $P$ is not maximal (as a joinless code) then it has no parse
tree (according to our definition of parse tree).

By Lemma \ref{SaturatedLeaf}(2), every maximal joinless code in 
$\,2\,\{0,1\}^*$ has at least one parse tree. But in $nA^*$ with $n \ge 3$
there are maximal joinless codes that have no parse tree, by Lemma 
\ref{SaturatedLeaf}(3); geometrically, codes in $3\,\{0,1\}^*$ without parse
tree correspond to tilings of the cube that cannot be obtained by successive 
bipartitions of cuboids (perpendicularly to an axis).  
This motivates the following.

\noindent {\bf Questions:} Is there a simple geometric or combinatorial
characterization of the finite maximal joinless codes in $nA^*$ (for 
$n \ge 3$) that have no parse tree?  
Is the non-existence of a parse tree equivalent to the presence of one of 
certain joinless subsets (``forbidden patterns'')? An example of such a 
forbidden pattern is the subset 
$\,\{(0,0,\varepsilon), \,(1,\varepsilon,0), \,(\varepsilon,1,1)\}\,$
of Lawson and Vdovina \cite{LawsonVdovina}, used in \ref{SaturatedLeaf}(3).

\bigskip

The following algorithm nondeterministically constructs any parse tree of 
$P$, if a parse tree exists. If $P$ has no parse tree the algorithm will 
discover this for some (but not all) nondeterministic choices. 
For a finite joinless code $P \subset \,2\,\{0,1\}^*$, the deterministic 
version of the algorithm decides whether $P$ is maximal (as a joinless code).

{\it Outline of the algorithm: } 
Initially, the algorithm puts $P$ into $T$ (as its leaf set), and makes a 
working copy $P_0$ of $P$. 
The algorithm keeps looking for an initial factor $v$ of an element of 
$P_0$ such that 
$\, v \cdot (\{\varepsilon\}^{i-1} \times A \times \{\varepsilon\}^{n-i}) \,$
$\subseteq$ $P_0$ \ (for some $i \in \{1, \, \ldots, n\}$). 
When such a $v$ is found, it is added to $T$ and to $P_0$; and 
$\, v \cdot (\{\varepsilon\}^{i-1} \times A \times \{\varepsilon\}^{n-i}) \,$ 
is removed from the working copy $P_0$.
If $(\varepsilon)^n$ is reached, and put into $T$, the construction of $T$ 
is complete and the algorithm concludes that $P$ is maximal (as a joinless 
code), and that it has a parse tree. 

The algorithm can be made deterministic by picking a total order for $nA^*$ 
(e.g., the lexicographic dictionary order), and always picking the first 
$v$ that works. 

Notation: $\, {\sf init}(P_0) \,$ denotes the set of {\em strict} 
initial factors of the elements of $P_0$; because of strictness (and since
$P_0$ is joinless), $P_0 \, \cap \, {\sf init}(P_0) = \varnothing$.

\bigskip

\begin{minipage}{\textwidth}
\noindent {\bf Algorithm}

\smallskip

\noindent {\sc Input:} A finite set $P \subset nA^*$, 
 given by a list of $n$-tuples of strings in $A^*$.

\noindent {\sc Precondition:} $P \ne \{(\varepsilon)^n\}$, and $P$ is 
joinless.  (This can easily be checked, by Lemma \ref{joinAn}.)

\noindent {\sc Output:} A set of vertices $V(T)$ and edges $E(T)$ of a parse
tree of $P$, if $P$ has a parse tree; 

\medskip

$P_0 := P$; \hspace{.2in} 
            {\small \# $P_0$ is a a working copy of $P$}

$V(T) := P$; \ \ \ $E(T) := \varnothing$;

{\tt while} $\, (\exists v \in {\sf init}(P_0))$ 
$(\exists i \in \{1, \, \ldots, n\})$ 
$[\, v \cdot (\{\varepsilon\}^{i-1} \times A \times \{\varepsilon\}^{n-i})$
$\, \subseteq \,$ $P_0 \,]$: 

\hspace{.4in} choose any $v$ that satisfies the while-condition;

\hspace{.4in} {\small \# for a deterministic algorithm, pick the first $v$ 
    that works (in a fixed total order) }

\hspace{.4in} $V(T) := V(T) \cup \{v\}$; 
 
\hspace{.4in} $E(T) := E(T)$ $\cup$ \ set of all edges from $v$ to the 
 elements of
$\, v \cdot (\{\varepsilon\}^{i-1} \times A \times \{\varepsilon\}^{n-i})$;
 
\hspace{.4in} $P_0 := (P_0$  $ \,\minus \, $  
$\, v \cdot (\{\varepsilon\}^{i-1} \times A \times \{\varepsilon\}^{n-i}))$
$ \ \cup \ $ $\{v\}$; \ \ \ \# Hence $P_0$ remains joinless.

\smallskip

{\tt if} $\, (\varepsilon)^n \in V(T)$: 

\hspace{.4in} {\tt then} output $(V(T), E(T))$ and conclude that $P$ is 
maximal;

\hspace{.4in}  {\tt else} (in case $n = 2$ and $A = \{0,1\}$) conclude that 
$P$ is not maximal
 
\hspace{.8in} (and hence has no parse tree).

\noindent  $\Box$
\end{minipage}

\medskip

\begin{pro} \label{PROP_Ptree}
 \ Let $P$ be any finite joinless code in $\,2\,\{0,1\}^*$.

Then $P$ has a parse tree iff $P$ is maximal as a joinless code.  

The Algorithm (deterministic version) decides maximality of $P$ and finds a
parse tree in polynomial time, when $P$ is given as a list of pairs of 
bitstrings.
\end{pro}
{\sc Proof.} The Algorithm uses one-step extensions of maximal joinless
codes; by Lemma \ref{PjoinlessPA}, each one-step extension or restriction 
preserves joinlessness and maximality.
Since $\{\varepsilon\}^n$ is a maximal joinless code, it follows that 
$P$ is maximal if the root $(\varepsilon)^n$ is reached.  It follows also 
that if the root is reached, a parse tree of $P$ exists (and the Algorithm
returns such a tree).

Conversely (for $n=2$ and $A = \{0,1\}$), if $P$ (or, at any later stage, 
$P_0$) is maximal, then by Lemma \ref{SaturatedLeaf}(2) there exists $v$ in 
the interior {\sc dag} such that 
$\, v \cdot (\{\varepsilon\} \times \{0,1\}$    $\,\cup\,$ 
$\{0,1\} \times \{\varepsilon\})$ $\subseteq$ $P$ (or $\subseteq$ $P_0$).
And this process does not stop until $P_0 = \{\varepsilon\}^n$. 
 \ \ \ $\Box$

\begin{cor} \label{OneStepAll} {\bf (cardinality of joinless codes).} 

\smallskip

\noindent Let $n$ be any positive integer and $A$ any finite alphabet.

\smallskip

\noindent {\rm (0.1)} \ For every $k_1, \ldots, k_n \in {\mathbb N}${\rm :}
 \ \ {\large \sf X}$_{_{i=1}}^{^n} A^{k_i} \ $
is a maximal joinless code that has a parse tree. 

\smallskip

\noindent {\rm (0.2)} \ 
For any finite joinless code $P \subset nA^*$: \ $P$ is maximal \ iff  
 \ $P$ can be transformed into \ {\large \sf X}$_{_{i=1}}^{^n} A^{k_i}\,$ 
by a finite sequence of restriction steps, where $k_i = $
$\max \{|v_i| : (v_1, \ldots, v_n) \in P\}\,$ for $1 \le i \le n$.

\medskip

\noindent {\rm (1)} \ For every finite maximal joinless code $P$
$\subseteq$  $nA^*$ there exists $N \in {\mathbb N}$ such that

\medskip

 \ \ \  \ \ \  \ \ \ $|P| \,=\, 1 + (|A| - 1) \cdot N$.  

\medskip

\noindent {\rm (1.1)} If $P$ has a parse tree then $P$ can be obtained 
from $\{\varepsilon\}^n$ by a finite sequence of one-step restrictions.
The number of one-step restrictions used is equal to the number of interior
vertices of every parse tree of $P$, and is equal to 
$\, N = (|P| - 1)/(|A| - 1)$.

\smallskip

\noindent {\rm (1.2)} If $P$ has no parse tree, then $P$ can be obtained from
$\{\varepsilon\}^n$ by a finite sequence of one-step restrictions, followed 
by a finite sequence of one-step extensions. 

Even when $P$ has no parse tree, $N = (|P| - 1)/(|A| - 1)$ 
is still the number of interior vertices in any parse tree of any maximal 
joinless code that has a parse tree and that has the same cardinality as 
$P\,$ (e.g., of the form $P_1 \times \{\varepsilon\}^{n-1}$ where $P_1$ is a 
prefix code in $A^*$).

\smallskip

\noindent {\rm (2)} \ Conversely, for all $N \in {\mathbb N}$ there are 
maximal joinless codes in $nA^*$ of cardinality $\, 1 + (|A| - 1) \cdot N$.
In particular, when $|A| = 2$ every positive integer is the cardinality of
some maximal joinless code.
\end{cor}
{\sc Proof.} (0.1) Let $\,C(k_1, \ldots, k_i, \ldots, k_n)$  $\,=\,$
{\large \sf X}$_{_{i=1}}^{^n} A^{k_i}$.
Let us prove by induction on $\,\sum_{i=1}^n k_i$ that 
$C(k_1, \ldots, k_i, \ldots, k_n)$ has a parse tree.  For
$C(0\ldots, 0, \ldots, 0) = \{\varepsilon\}^n$, the parse tree consists 
of one vertex. Inductively, 

\smallskip

$C(k_1, \ldots, k_{i-1}, k_i + 1, k_{i+1}, \ldots, k_n)$ 
 \ $=$ 
 \ {\large \sf X}$_{_{j=1}}^{^{i-1}} A^{k_j}$
$\,\times\, A^{k_i +1} \times\,$
{\large \sf X}$_{_{j=i+1}}^{^n} A^{k_j}$

\smallskip

$= \ $
({\large \sf X}$_{_{j=1}}^{^n} A^{k_j}$)
$\cdot$
$(\{\varepsilon\}^{i-1} \times A \times \{\varepsilon\}^{n-i-1})$

\smallskip

$= \ $
$C(k_1, \ldots, k_i, \ldots, k_n)$
$\cdot$
$(\{\varepsilon\}^{i-1} \times A \times \{\varepsilon\}^{n-i-1})$

\smallskip

$= \ $
$\bigcup_{v \in A^{k_i}}$
$\big($({\large \sf X}$_{_{j=1}}^{^{i-1}} A^{k_j}$ 
$\times \{v\} \times$
{\large \sf X}$_{_{j=i+1}}^{^n} A^{k_j}$)   
$\cdot$
$(\{\varepsilon\}^{i-1} \times A \times \{\varepsilon\}^{n-i-1}) \big)$.

\smallskip

\noindent So, $C(k_1, \ldots, k_{i-1}, k_i + 1, k_{i+1}, \ldots, k_n)$ is 
obtained form $C(k_1, \ldots, k_i, \ldots, k_n)$ by $|A|^{k_i}$ one-step 
restrictions (one one-step restriction for every $v \in A^{k_i}$). 
It follows that if $C(k_1, \ldots, k_i, \ldots, k_n)$ has a parse tree then
$C(k_1, \ldots, k_{i-1}, k_i + 1, k_{i+1}, \ldots, k_n)$ has a parse tree.
Moreover, any joinless code that has a parse tree is maximal.

\smallskip

\noindent (0.2) Let 

\smallskip

 \ \ \  \ \ \ $\ell_i(P)$  $\,=\,$ 
$\max \{|v_i| : (v_1, \ldots, v_n) \in P\}$, for $1 \le i \le n$; \ and 

\smallskip

 \ \ \  \ \ \ $\nu(P)$  $\,=\,$
$\,\prod_{i=1}^n |A|^{\ell_i(P)} \ - \ \sum_{u \in P} \sum_{i=1}^n |u_i|\,$.

\smallskip

\noindent The fact that $P$ can be restricted to 
$C(\ell_1(P), \,\ldots\, , \ell_n(P))$ follows by induction on $\nu(P)$:

If $\nu(P) = 0$ then $P = C(\ell_1(P), \,\ldots\, , \ell_n(P))$.

If $\nu(P) > 0$, and $(u_1, \ldots, u_n) \in P$ is such that 
$|u_i| < \ell_i(P)$ for some $i$, then a one-step restriction decreases 
$\nu(P)$, as $(u_1, \ldots, u_n)$ is replaced by 
$\,(u_1, \ldots, u_n) \cdot$
$(\{\varepsilon\}^{i-1} \times A \times \{\varepsilon\}^{n-i-1})$. 

\smallskip

\noindent (1.1) We prove the equivalent statement that from $P$ one can 
reach $\{\varepsilon\}^n$ by $N = (|P| - 1)/(|A| - 1)$ ones-step extensions.
We use induction on $|P|$.  When $|P| = 1$ then $P = \{\varepsilon\}^n$, and
the formula holds. 
For $|P| > 1$, an extension step can be applied to some leaf of the interior
of a parse tree of $P$, by Lemmas \ref{SaturatedLeaf} and \ref{PjoinlessPA}.   
In this extension step, a new maximal joinless code $Q$ is obtained; one 
leaf of the interior the parse tree of $P$ becomes a leaf of the parse tree 
of $Q$, so this parse tree of $Q$ has $N-1$ interior vertices; and 
$|Q| = |P| - (|A| - 1)$.
By induction, $N-1 = (|Q| - 1)/(|A| -1)$; and the latter is equal to 
$(|P| - (|A| - 1) - 1)/(|A| -1)$  $=$  $(|P| - 1)/(|A| - 1) - 1$. Hence 
$N = (|P| - 1)/(|A| - 1)$. 

\smallskip

\noindent (1.2) By applying one-step restrictions as in part (0.2), from any
maximal joinless code $P$ one can reach 
{\large \sf X}$_{_{i=1}}^{^n} A^{k_i}$, where $k_i$ is as in part (0.2). And
from {\large \sf X}$_{_{i=1}}^{^n} A^{k_i}$ one can reach 
$\{\varepsilon\}^n$ by one-step extensions by (0.1). In any one-step
restriction or extension the cardinality of the maximal joinless code
increases or decreases by $|A|-1$.
So, to reach $P$ from $\{\varepsilon\}^n$ we can apply restrictions to reach
{\large \sf X}$_{_{i=1}}^{^n} A^{k_i}$, then apply extensions to obtain $P$.

\smallskip

\noindent (1) The formula follows from (1.1) and (1.2).

\smallskip

\noindent
(2) For the existence of codes of the given cardinality, take for example 
$\,Q \times \{\varepsilon\}^{n-1}$, where $Q$ is any maximal prefix code in 
$A^*$, and apply the corresponding result for maximal prefix codes (which is 
folklore; see e.g.\ \cite[Lemma 9.9(0)]{BiCoNP}).
 \ \ \ $\Box$

\begin{pro} \label{PolyAlgos}
 \ There exist {\em polynomial-time algorithms} that on input 
$P \subset nA^*$ (a finite set, given by an explicit list of $n$-tuples 
of strings) decide whether $P$ has the following properties: \\  
{\rm (1)} \ $P$ is {\em joinless}; \\  
{\rm (2)} \ $P$ is {\em maximal} as a joinless code.  
\end{pro}
{\sc Proof.} The input to the algorithms is $P$, given as a list of 
$n$-tuples of strings, so the input size is
$\, \sum_{p \in P} \sum_{i=1}^n |p_i|$.
 
\smallskip

\noindent (1) Lemma \ref{joinAn}, applied to every two elements $u, v \in P$
with $u \ne v$, will decide in quadratic time whether $P$ is joinless. 

\noindent (2) For finite joinless codes in $\,2\,\{0,1\}^*$,  
the Algorithm given after Def.\ \ref{DEFtreeofP} has polynomial time 
complexity, in view of Coroll.\ \ref{OneStepAll} which proves that every 
parse tree of $P$ has size that is linearly bounded in terms of $|P|$.

For $nA^*$ in general, the Algorithm that follows Prop.\ \ref{PROnKraft}, 
based on the generalized Kraft equality, decides in polynomial time whether 
a joinless code is maximal.
 \ \ \ $\Box$

\medskip

The corresponding questions about initial factor codes are also decidable in 
polynomial time. Suppose $P \subset nA^*$ is finite and given by an 
explicit list of $n$-tuples of strings. It is easy to decide whether $P$ 
is an initial factor code; it is sufficient to check for every two elements 
$u, v \in P$ with $u \ne v$, whether $u \, \|_{\rm init} \, v$.
A code $P$ is maximal as an initial factor code iff for every interior vertex 
$v$ of the $P$-{\sc dag}, the $n \, |A|$ children of $v$ in the 
$nA^*$-{\sc dag} are also children of $v$ in the $P$-{\sc dag}.

\bigskip

An algorithm for testing maximality of a joinless code can be derived from 
the following generalization of the Kraft (in)equality to higher dimensions. 
We mentioned earlier that in the geometric description of the
Brin-Thompson groups, a word $x = (x_1, \ldots, x_n) \in n\, \{0,1\}^*$
represents the hyperrectangle
 \ {\large \sf X}$_{_{i=1}}^{^n} [0.x_i \, , \ 0.x_i + 2^{-|x_i|}[ \ $ 
(where we close the intervals whose right-bound is 1).
The measure of this hyperrectangle is
 \ $2^{-(|x_1| \, + \ \ldots \ + \, |x_n|)}$.
More generally, we have the following.

\begin{defn} \label{DEFmeasure}
 \ Let $A$ be an alphabet of cardinality $|A| = k \ge 2$. For every
$\, x = (x_1, \ldots, x_n) \in nA^*\,$ we define the measure

\smallskip

 \ \ \  \ \ \ $\mu(x) \,=\, k^{-(|x_1| + \, \ldots \, + |x_n|)}$.

\smallskip

\noindent For every joinless code $P \subset nA^*$ (not necessarily 
finite) we define the measure 

\smallskip

 \ \ \  \ \ \ $\mu(P) \,= \,\sum_{x \in P} \,\mu(x)$.
\end{defn}

\begin{pro} \label{PROnKraft}  {\bf ({\boldmath $n$}-dimensional Kraft 
(in)equality).}
 \ Let $P \subset nA^*$ be a finite joinless code, where 
$|A| \ge 2$ and $n \ge 1$. Then we have:

\smallskip

{\rm (1)} \ \ $\mu(P) \le 1$.

\smallskip

{\rm (2)} \ \ $P$ is maximal (as a joinless code) \ iff \ $\mu(P) = 1$.
\end{pro}
{\sc Proof.} This follows from the geometric picture. For a joinless code 
$P$, all the words in $P$ represent non-overlapping hyperrectangles in 
$[0,1]^n$, so their total measure is at most the measure of $[0,1]^n$, 
which is 1.

And $P$ is maximal iff the corresponding hyperrectangles tile $[0,1]^n$,
which is iff the sum of the measures of the hyperrectangles is 1.
 \ \ \ $\Box$

\medskip

\noindent Prop.\ \ref{PROnKraft} probably holds for infinite joinless codes
too; but since we don't need it in that case, we'll that question open.

\smallskip

\noindent Prop.\ \ref{PROnKraft} leads to the following algorithm.

\bigskip

\begin{minipage}{\textwidth}

\noindent {\bf Algorithm (maximality of a finite joinless code)}

\noindent {\sc Input:} A finite set $P \subset nA^*$, given as an explicit 
list of words.

\noindent {\sc Precondition:} $P$ is joinless. (This is easily checked, by 
Prop.\ \ref{PolyAlgos}(1).)

\noindent {\sc Question:} Is $P$ maximal?

\medskip

Compute $ \ \mu(P) = \sum_{x \in P} k^{- \sum_{i=1}^n |x_i|}$ \ in fractional
base-$k$ representation;

if $\, \mu(P) = 1$, output ``yes'';

else, output ``no''. \ \ \  \ \ \ $\Box$ 

\end{minipage}

\bigskip

\noindent This algorithm runs in polynomial time, in terms of the total input
length $\, \sum_{x \in P} \sum_{i=1}^n |x_i|$. In fractional
base-$k$ representation the sum $\mu(P)$ is easy to compute.

\bigskip

\noindent We will need the intersection of joinless generated right ideals, 
and the elementwise join of joinless codes.

\begin{pro} \label{INTERSjoinless} 
 \ Let $P, Q \subset nA^*$ be joinless codes.

\smallskip

\noindent {\rm (1)} \ \ The {\em elementwise join} $P \vee Q$, defined by

\smallskip

\hspace{.5in} $P \vee Q \, = \, \{p \vee q \, : \, p \in P, \ q \in Q\}$,  

\smallskip

 \ \ is a joinless code.
(Here, $p \vee q$ ranges over the joins that exist.)

\smallskip

 \ \   Hence, \ $|P \vee Q| \, \le \, |P| \cdot |Q|$. 

\smallskip

\noindent {\rm (2)} \ \ $P$ and $Q$ are both maximal (as joinless codes) 
 \ iff \ $P \vee Q$ is maximal.

\smallskip

\noindent {\rm (3)} \ \ $(P \vee Q) \cdot (nA^*) \ = \ $
    $P \cdot (nA^*) \ \cap \ Q \cdot (nA^*)$. 

\smallskip

 \ \ Hence, if $P \, (nA^*)$ and $Q \, (nA^*)$ are joinless
generated then so is $\, P \, (nA^*) \, \cap \, Q \, (nA^*)$.
\end{pro}
{\sc Proof.} (1) Suppose $p, p' \in P$, $q, q' \in Q$, and $p \ne p'$ or  
$q \ne q'$. Then \, $(p \vee q) \vee (p' \vee q') \,$ does not exist, 
because $(p \vee q) \vee (p' \vee q')$ would have $p$, $p'$, $q$, 
and $q'$ as prefixes. But either $p$ and $p'$ (if different) or $q$ and
$q'$ (if different) do not have a join.  

\smallskip

\noindent  
(2) $[\Leftarrow]$ If $P \vee Q$ is maximal then every $x \in nA^*$ has 
a join with some $p \vee q \in P \vee Q$, i.e., $x$ and $p \vee q$ are 
initial factors of some $z \in nA^*$. Then $p$ and $q$ are also initial
factors of $z$, so $x \vee p$ and $x \vee q$ exist. Hence, every $x \in nA^*$ 
has a join with some $p \in P$ and some $q \in Q$, thus $P$ and $Q$ are 
maximal. \\   
$[\Rightarrow]$ If $P$ is maximal then every $x \in nA^*$ has a join with
some $p \in P$; and if $Q$ is maximal, $x \vee p$ has a join with some 
$q \in Q$. Hence, $x$, $p$, and $q$, are all initial factors of some word $z$,
hence $z \vee p \vee q$ exists. So, every $x \in nA^*$ has a join
with some $p \vee q$, so $P \vee Q$ is maximal. 

\smallskip

\noindent
(3) $[\supseteq]$ Every $w \in P \, (nA^*) \, \cap \, Q \, (nA^*)$ 
satisfies $w = p u = q v$ for some $p \in P$, $q \in Q$, and 
$u, v \in nA^*$. 
This implies that $p$ and $q$ are initial factors of $w$, so $p \vee q$ 
exists, and is an initial factor of $w$. Hence, 
$w \in (P \vee Q) \cdot (nA^*)$. 

\noindent
$[\subseteq]$ If $p \vee q$ exists then it has $p$ and $q$ as initial
factors, hence $p \vee q \in P \, (nA^*) \, \cap \, Q \, (nA^*)$.  
 \ \ \ $\Box$

\subsection{Right ideal morphisms of {\boldmath $nA^*$}, and string-based 
definition of {\boldmath $n G_{k,1}$} and {\boldmath $n V$} }

Just as for $A^*,$ one defines the concepts of {\it right ideal morphism},
domain code, and image code in $nA^*$.
We only consider domain and image codes that are {\it joinless}. 
Indeed, if $P \subset nA^*$ is not joinless, some definitions of right 
ideal morphisms on $P$ will be inconsistent.
E.g., let $P = \{(0, \varepsilon), \, (\varepsilon, 0)\}$, 
so $(0, \varepsilon) \vee (\varepsilon, 0) = (0,0)$; and let
$f(0, \varepsilon) = (0,0)$ and $f(\varepsilon, 0) = (1,1)$;
then $\, f(0,0) = f((0, \varepsilon) \cdot (\varepsilon, 0)) = $
$(0,0) \cdot (\varepsilon, 0) = (0,00) \ne $
$(10,1) = (1,1) \cdot (0, \varepsilon) =$ 
$f((\varepsilon, 0) \cdot (0, \varepsilon)) = f(0,0)$; so $f(0,0)$ receives
two different values.

\bigskip

\noindent Before we get to $n G_{k,1}$ we define the following monoid:

\begin{defn} \label{DEFnRIfin}\hspace{-0.1in}.

\smallskip

$n {\cal RI}_A^{\sf fin}$  $ \ = \ $ 
$\{f : \, f$ is a right ideal morphism of $nA^*$ such that $f$ is injective, 

\hspace{1.06in} and ${\rm domC}(f)$ and ${\rm imC}(f)$ 
are {\em finite, maximal, joinless} codes\} .
\end{defn}
``Maximal'' means maximal as a joinless code.  
Usually we just write $n {\cal RI}^{\sf fin}$ when a fixed alphabet $A$ 
is used.

\begin{lem} \label{InjNormal}
 \ For every $f \in n {\cal RI}^{\sf fin}$: 
 \ \ $f({\rm domC}(f)) = {\rm imC}(f)$.

\smallskip

Hence, if $f \in n {\cal RI}^{\sf fin}$ then 
$f^{-1} \in n {\cal RI}^{\sf fin}$, and ${\rm domC}(f^{-1}) = {\rm imC}(f)$, 
 \ ${\rm imC}(f^{-1}) = {\rm domC}(f)$.
\end{lem}
{\sc Proof.} For every $p_1 \in {\rm domC}(f)$:
$\, f(p_1) = q_1 u \in {\rm Im}(f)$, for some $q_1 \in {\rm imC}(f)$ and 
$u \in nA^*$. 
Since $q_1 \in {\rm Im}(f)$, $q_1 = f(p_2 v)$ for some 
$p_2 \in {\rm domC}(f)$ and $v \in nA^*$.
Hence, $q_1 u = f(p_2 v) \ u = f(p_2 v u)$.  Thus, 
$f(p_1) = q_1 u = f(p_2 v u)$.
Since $f$ is injective, this implies that $p_1 = p_2 v u$. Since 
$p_1, p_2 \in {\rm domC}(f)$, which is an initial factor code, $p_1 = p_2$ 
and $u = v = (\varepsilon)^n$. Hence, 
$f(p_1) = q_1 u = q_1 \in {\rm imC}(f)$. 
So $f({\rm domC}(f)) \subseteq {\rm imC}(f)$. 

Conversely, if $q \in {\rm imC}(f)$, then $q = f(p) \, v$ for some 
$p \in {\rm domC}(f)$ and $v \in nA^*$. Since $f(p) \in {\rm Im}(f)$
and $q \in {\rm imC}(f)$ (which is the initial factor code that generates 
${\rm Im}(f)$), we conclude that $q = f(p)$ and $v = (\varepsilon)^n$.
Hence, $q \in f({\rm domC}(f)$. 
So, ${\rm imC}(f) \subseteq f({\rm domC}(f))$.   

Now $f^{-1}$ satisfies the following: For $q \in {\rm  imC}(f)$, 
 \ $f^{-1}(q) = p$ iff $p \in {\rm domC}(f)$ and $f(p) = q$. 
Hence $f^{-1} \in n {\cal RI}^{\sf fin}$, and 
${\rm domC}(f^{-1}) = {\rm imC}(f)$, and ${\rm imC}(f^{-1}) = {\rm domC}(f)$.
 \ \ \ $\Box$

\begin{lem} \label{LEM_f_of_P}
 \ Let $f \in n {\cal RI}^{\sf fin}$ and let $P \subset nA^*$ be a 
finite set.

\smallskip

\noindent {\rm (1.1)} \ If $P \subset {\rm Dom}(f)$ we have: 
 \ \ $f(P)$ is joinless \ iff \ $P$ is joinless.

\smallskip

\noindent {\rm (1.2)} \ If $P \subset {\rm Dom}(f)$ and $P$ is joinless, 
we have:  \ \ $P$ is maximal \ iff \ $f(P)$ is maximal.

\smallskip

\noindent {\rm (2.1)} \ In general (not assuming $P \subset {\rm Dom}(f)$),
we have: 

\smallskip

\hspace{0.5in}  
$f(P \,\vee\, {\rm domC}(f))$ is joinless \ iff \ $P$ is joinless.  

\smallskip

\noindent {\rm (2.2)} \ In general, if $P$ is joinless then the following
are equivalent:

\smallskip

\hspace{0.5in} $P$ is maximal,  

\hspace{0.5in} $P \vee {\rm domC}(f) \ $ is maximal, 

\hspace{0.5in} $f(P \vee {\rm domC}(f)) \ $ is maximal.
\end{lem}
{\sc Proof.} (1.1) $[\Leftarrow]$ \ Let $p, q \in P$, and assume by 
contradiction that there exists $z \in nA^*$ such that $f(p)$ and $f(q)$ 
are initial factors of $z$.  Then $z = f(p) \, u = f(q) \, v$ for some 
$u,v \in nA^*$. Hence, $f^{-1}(z) = f^{-1}(f(p) \, u)$ $=$ 
$f^{-1}(f(p)) \ u$; the latter holds since 
$f^{-1} \in n {\cal RI}^{\sf fin}$, and $f(p) \in {\rm Dom}(f^{-1})$ 
$=$  ${\rm Im}(f)$ \ (by Lemma \ref{InjNormal}).
Hence, $f^{-1}(z) = p u$. Similarly, $f^{-1}(z) = qv$. So, $p u = q v$,
but that contradicts the assumption that $P$ is joinless.

\noindent (1.1) $[\Rightarrow]$ \ Conversely, if some $p, q \in P$ have a 
join $z$ then $z = pu = qv$ for some $u, v \in \in nA^*$. Then 
$\, f(z) = f(p) \, u = f(q) \, v$, so $f(p) \vee f(q)$ exists, hence $f(P)$ 
is not joinless. 

\noindent (1.2) $[\Rightarrow]$ \ Suppose $P$ is maximal, and assume by 
contradiction that $f(P)$ is not maximal. Then there exists $x \in nA^*$ 
such that $\{x\} \cup f(P)$ is a joinless code. Since 
$f^{-1} \in n {\cal RI}^{\sf fin}$, 
$f^{-1}(\{x\} \cup f(P))$ is joinless (by what was proved in the previous 
paragraph). So, $f^{-1}(\{x\} \cup f(P)) = P \cup \{f^{-1}(x)\}$ is 
joinless, which contradicts $P$ the assumption that $P$ is maximal.
Thus, if $P$ is maximal then $f(P)$ is maximal. 

\noindent (1.2) $[\Leftarrow]$ \ Similarly, if $f(P)$ is maximal then 
$f^{-1}f(P)$ is maximal (since $f^{-1} \in n {\cal RI}^{\sf fin}$). Hence 
if $f(P)$ is maximal, $P$ is maximal.

\noindent (2.1) If $P$ is joinless iff $P \vee {\rm domC}(f)$ is joinless, by
Lemma \ref{INTERSjoinless}(1), since ${\rm domC}(f)$ is joinless for all
$f \in n {\cal RI}^{\sf fin}$. And $P \vee {\rm domC}(f)$ is joinless iff 
$f(P \vee {\rm domC}(f))$ is joinless, by (1.1).

\noindent (2.2) If $P$ is maximal then $P \vee {\rm domC}(f)$ is maximal
by Lemma \ref{INTERSjoinless}(2), since ${\rm domC}(f)$ is maximal for 
$f \in n {\cal RI}^{\sf fin}$. This implies that $f(P \vee {\rm domC}(f))$ 
is maximal, by (1.2). And if $f(P \vee {\rm domC}(f))$ is maximal then 
$P \vee {\rm domC}(f)$ is maximal, again by (1.2). Moreover, maximality of
$P \vee {\rm domC}(f)$ implies maximality of $P$ (and of ${\rm domC}(f)$),
by Lemma \ref{INTERSjoinless}(2).  
 \ \ \ $\Box$

\bigskip

Every right ideal morphism $f \in n {\cal RI}^{\sf fin}$ is uniquely 
determined by its restriction to ${\rm domC}(f)$; this is an obvious 
consequence of the fact that $f$ is a right-ideal morphism and 
${\rm domC}(f)$ is a joinless code.
So $f$ is determined by the finite function 
$\, f$: ${\rm domC}(f) \to {\rm imC}(f)$. 

Conversely, let $P, Q \subset nA^*$ be two finite maximal joinless codes 
with the same cardinality, and let $F$: $P \to Q$ by any bijection from $P$ 
onto $Q$. Then $F$ determines a right ideal morphism $f$ of $nA^*$, such 
that $F$ is the restriction of $f$ to its domain code; $f$ is defined in a 
unique way by $f(p v) = F(p) \ v$ for all $p \in P$, $\, v \in nA^*$. 
Since $P$ is joinless, $f$ is well defined.

\begin{defn} \label{DEFtable} {\bf (table).} 
 \ A bijection $F$: $P \to Q$ between finite maximal joinless codes 
$P, Q \subset nA^* \,$ is called a {\em table}. 
\end{defn}
Tables and right ideal morphisms in $n {\cal RI}^{\sf fin}$ determine each 
other bijectively, and can the treated as ``the same thing''. 

\bigskip

Every function $f \in n {\cal RI}^{\sf fin}$ determines a permutation 
$f^{(\omega)}$ of $nA^{\omega}$, as follows.
For any $w \in nA^{\omega}$ there exists a unique $p \in {\rm domC}(f)$ 
such that $w = pu$ for some $u \in nA^{\omega}$, by Lemma
\ref{OmegaJoilessPref}.  
Then we define $f^{(\omega)}$ by

\smallskip

\hspace{1.in}   $\, f^{(\omega)}(w) = f(p) \ u$. 

\medskip

\noindent The converse does not hold; i.e., $f \in n {\cal RI}^{\sf fin}$ 
is not determined by $f^{(\omega)}$, as will be seen in Lemma
\ref{LEMendequiv}.

\begin{defn} \label{DEFomegaequiv} {\bf (end-equivalence).}
\ Two right ideal morphisms $f, g \in n {\cal RI}^{\sf fin}$ are
{\em end-equivalent} \ iff \ $f$ and $g$ agree on
$\,{\rm Dom}(f) \,\cap\, {\rm Dom}(g)$. 
This will be denoted by $f \equiv_{\rm end} g$. 
\end{defn}
By Prop.\ \ref{INTERSjoinless}, $\, {\rm Dom}(f) \cap {\rm Dom}(g)$ is
generated by a joinless code, namely
$\, {\rm domC}(f) \vee {\rm domC}(g)$. 

In \cite{BiCongr} the congruence $\equiv_{\rm end}$ is defined in much 
greater generality, and other congruences are introduced.

\begin{lem} \label{LEMendequiv}
 \ For all $f, g \in n {\cal RI}^{\sf fin}$:
  \ $f \equiv_{\rm end} g$ \ iff \ $f^{(\omega)} = g^{(\omega)}$.
\end{lem}
{\sc Proof.} For every $f \in n {\cal RI}^{\sf fin}$, ${\rm domC}(f)$ and 
${\rm imC}(f)$ are maximal joinless codes. Therefore (by Lemma
\ref{JoinlessEssential}): \ ${\rm domC}(f) \cdot (nA^{\omega})$  $=$ 
$nA^{\omega}$  $=$  ${\rm imC}(f) \cdot (nA^{\omega})$.
And by Lemma \ref{INTERSjoinless}, 
$\, {\rm domC}(f) \vee {\rm domC}(g)$ is also a maximal joinless code.

Let $R = {\rm Dom}(f) \cap {\rm Dom}(g)$, and let $f|_R$ and $g|_R$ be the 
restrictions of $f$ or $g$ to $R$. Then $f \equiv_{\rm end} g$ is equivalent
to $f|_R = g|_R$. 

\smallskip

$[\Rightarrow]$ Suppose $f \equiv_{\rm end} g$, i.e., $f|_R = g|_R$, where
$R = {\rm Dom}(f) \cap {\rm Dom}(g)$.
For every $w \in nA^{\omega}$, let $z \in nA^*$ be an initial factor of $w$
such that in all coordinates, $z$ is longer than the longest coordinate of
any element of $P = {\rm domC}(f) \vee {\rm domC}(g)$. And $w = z u$ for 
some $u \in nA^{\omega}$.  
Since $P$ is a maximal joinless code, $z$ has a join with an element of 
$P$; by the chosen length of $z$, $z$ has an initial factor in $P$, hence 
$z \in R$. 
Now $f^{(\omega)}(z u) =  f(z) \ u$, since $z \in R \subseteq {\rm Dom}(f)$;
and $g^{(\omega)}(z u) = g(z) \ u$, since $z \in R \subseteq {\rm Dom}(g)$.
Since $f(z) = g(z)$ (because $f|_R = g|_R$), if follows that 
$f^{(\omega)}(z u) = g^{(\omega)}(z u)$.

$[\Leftarrow]$ Suppose $f^{(\omega)} = g^{(\omega)}$.
For every $r \in R$ and every $u \in nA^{\omega}$, 
$f^{(\omega)}(r u) = g^{(\omega)}(r u)$. And since $r \in R = $
${\rm Dom}(f) \cap {\rm Dom}(g)$, $f^{(\omega)}(r u) = f(r) \ u$, and 
$g^{(\omega)}(r u) = g(r) \ u$. From $f(r) \ u = g(r) \ u$ it follows that
$f(r) = g(r)$. Hence, $f|_R = g|_R$, i.e., $f \equiv_{\rm end} g$. 
 \ \ \ $\Box$

\begin{lem} \label{EquivFinCongr}
 \ For all $f_1, f_2 \in n {\cal RI}^{\sf fin}$: 
 \ \ $(f_2 \circ f_1)^{(\omega)} = f_2^{(\omega)} \circ f_1^{(\omega)} $.

\smallskip

The relation $\equiv_{\rm end}$ is a {\em congruence} on
$n {\cal RI}^{\sf fin}$.
\end{lem}
{\sc Proof.} For every $w \in nA^{\omega}$ there exist
$r \in {\rm Dom}(f_2 \circ f_1)$ and $u \in nA^{\omega}$ such that
$w = r u$; this follows from Lemma \ref{OmegaJoilessPref}. Then
$r \in {\rm Dom}(f_1)$ and $f_1(r) \in {\rm Dom}(f_2)$.
Now by the definition of $f^{(\omega)}(w)$,
 (\ $(f_2 \circ f_1)^{(\omega)}(w) = (f_2 \circ f_1)(r) \ u$
$= f_2(f_1(r)) \ u$.
And $f_2^{(\omega)}(f_1^{(\omega)}(r u)) = f_2^{(\omega)}(f_1(r) \, u)$
$= f_2(f_1(r)) \ u$; the latter holds since $f_1(r) \in {\rm Dom}(f_2)$.
This proves that $(f_2 \circ f_1)^{(\omega)}(w) = $
$f_2^{(\omega)}(f_1^{(\omega)}(w)$.

It follows immediately that $\equiv_{\rm end}$ is a congruence on
$n {\cal RI}^{\sf fin}$ (by Lemma \ref{LEMendequiv}).
 \ \ \ $\Box$

\bigskip

Next we develop criteria about extensions and restrictions of functions in 
$n {\cal RI}^{\sf fin}$ that enable us to decide efficiently whether two 
tables determine end-equivalent functions.
The Remark below applies to finite maximal joinless codes in 
$\,2\,\{0,1\}^*$ and is similar to a criterion for end-equivalence of finite 
maximal prefix codes in $A^*$.
But because of Lemma \ref{SaturatedLeaf}(3) it does not apply for $n \ge 3$.  
For $nA^*$ in general, Prop.\ \ref{FinExtVSsteps} gives an efficient 
algorithm for deciding whether two tables determine end-equivalent functions.

\medskip

\noindent {\bf Remark (extension-restriction criterion):}
Let $P, Q$ be finite maximal joinless codes in $\,2\,\{0,1\}^*$, let 
$F$: $P \to Q$ be a table, and let 
$f \in 2\,{\cal RI}^{\sf fin}$ be the corresponding right ideal morphism of
$\,2\,\{0,1\}^*$. Then: 
$f$ is extendable in $2\,{\cal RI}^{\sf fin}$ iff there exist 
$ \ p = (p_1, p_2), \ q = (q_1, q_2) \in 2\,\{0,1\}^*$ 
such that for every $a \in \{0,1\}$: 

 \ \ \ {\rm (1)} \ \ \ $\{(p_1, p_2 a): a \in \{0,1\}\} \,\subseteq\, P$ 
 \ (or \ $\{(p_1 a, p_2): a \in \{0,1\}\} \,\subseteq\, P$), \ and 

 \ \ \ {\rm (2)} \ \ \ $\{(q_1, q_2 a): a \in \{0,1\}\} \,\subseteq\, Q$ 
 \ (or \ $\{(q_1 a, q_2): a \in \{0,1\}\} \,\subseteq\, Q$), \ and 

 \ \ \ {\rm (3)} \ \ \ $F(p_1, p_2 a) \,=\, (q_1, q_2 a)$ 
 \ \ (or \ $F(p_1 a, p_2) \,=\, (q_1 a, q_2)$). 

\noindent In that case, let 

 \ \ \ $P' \, = \, $   $(P$ $\minus$
$\{(p_1, p_2 a) : a \in \{0,1\}\}) \ \cup \ \{p\}$, 
(or \ $(P \minus \{(p_1 a, p_2) : a \in \{0,1\}\})  \ \cup \ \{p\}$),
 \ and 

 \ \ \ $Q' \, = \, $   $(Q$ $\minus$
$\{(q_1, q_2 a) : a \in \{0,1\}\}) \ \cup \ \{q\}$,        
(or \ $(Q \minus \{(q_1 a, q_2) : a \in \{0,1\}\})  \ \cup \ \{q\}$).  

\noindent Then $P'$ and $Q'$ are finite maximal joinless codes in
$\,2\,\{0,1\}^*$, and $f$ can be extended to a function 
$f' \in 2\,{\cal RI}^{\sf fin}$ with table $F'$: $P' \to Q'\,$ defined by

 \ \ \ $F'(p) = q$, and 

 \ \ \ $F'(p') = F(p')\,$ for all $\, p' \in P$  $\minus$
$\{(p_1, p_2 a): a \in \{0,1\}\}$ \ (or \ $P$ $\minus$
$\{(p_1 a, p_2): a \in \{0,1\}\}$).

\noindent The passage from $f$ to $f'$ is called a {\em one-step extension},
and $f$ is called a {\em one-step restriction} of $f'$.

\smallskip

\noindent The Remark follows from Lemma \ref{SaturatedLeaf}(2), in the 
same way as for prefix codes in $A^*$ (see \cite[Lemma 2.2]{BiThomps} 
and \cite{Hig74}). The Remark is not always applicable when $n \ge 3$, by 
Lemma \ref{SaturatedLeaf}(3).

\begin{pro} \label{FinExtVSsteps} 
{\bf (restrictions, and deciding $\, \equiv_{\rm end}$).} 

\smallskip

\noindent {\rm (1)} \ Let $F: P \to Q$ be a table, and let 
$f \in n {\cal RI}_A^{\sf fin}$ be the right-ideal morphism given by this 
table. Suppose $P' \cdot (nA^*) \subseteq P \cdot (nA^*)$, where 
$P' \subset nA^*$ is a finite maximal joinless code.
Then the {\em restriction} $f' = f|_{P' \cdot (nA^*)}\,$ of $f$ to 
$\,P' \cdot (nA^*)\,$ is an element of $\,n {\cal RI}_A^{\sf fin}$ with table 
$\, F': P' \to f(P')$.

\smallskip

Moreover, $ \ f \equiv_{\rm end} f|_{P' \cdot (nA^*)}$.

\smallskip

\noindent {\rm (2)} \ Let $F^{(j)}: P^{(j)} \to Q^{(j)}$ be tables (for 
$j = 1,2$), and let $\,f^{(j)} \in n {\cal RI}_A^{\sf fin}$ be the 
right-ideal morphisms given by these tables.  Then:

\smallskip

 \ \ \ \ \ $f^{(1)} \equiv_{\rm end} f^{(2)}$ \ \ \ iff  
 \ \ \ $F^{(1)}|_{P^{(1)} \vee P^{(2)}}$   $\,=\,$ 
    $F^{(2)}|_{P^{(1)} \vee P^{(2)}}$ , 

\smallskip

\noindent where $F^{(j)}|_{P^{(1)} \vee P^{(2)}}$ is the table of the 
restriction of $f^{(j)}$ to the finite maximal joinless code 
$P^{(1)} \vee P^{(2)}$.

\smallskip

Hence there is a polynomial-time algorithm that decides whether the tables 
$F^{(1)}$ and $F^{(2)}$ represent $\equiv_{\rm end}$-equivalent elements of 
$n {\cal RI}_A^{\sf fin}$.
\end{pro}
{\sc Proof.} (1) The restricted table $F': P' \to f(P')$ is defined as 
follows. For every $p' \in P'$ there exists $p = (p_1, \ldots, p_n) \in P$ 
and $w = (w_1, \ldots, w_n) \in nA^*$ such that $p' = pw$. We define  
 \ $F'(p') = F(p) \ w$.

In order to verify that $F'$ is a well defined function, suppose that
$p' = p^{(1)} \,u = p^{(2)}\, v$ for some $p^{(1)}, p^{(2)} \in P$ and 
$u, v \in nA^*$. Since $P$ is joinless, it follows that $p^{(1)} = p^{(2)}$;
let $p^{(1)} = p^{(2)} = p$. Since multiplication in $nA^*$ is cancelative,
$p u = p v$ implies $u = v$.  
So, $p' \in P'$ determines a unique $p \in P$ and $w = u = v \in nA^*$
such that $p' = pw$. Hence $F'(p') = F(p) \ w$ defines $F'(p')$ in a
unique way.

\smallskip

\noindent (2) The largest common restriction of $f^{(1)}$ and $f^{(2)}$ is 
$f^{(1)} \cap f^{(2)}$, which has domain 
${\rm Dom}(f^{(1)}) \cap {\rm Dom}(f^{(2)}$, and domain code 
${\rm domC}(f^{(1)} \cap f^{(2)}) = P^{(1)} \vee P^{(2)}$ (by Lemma 
\ref{INTERSjoinless}). So, $f^{(1)} \equiv_{\rm end} f^{(2)}\,$ iff the 
tables of $f^{(1)}$ and $f^{(2)}$, restricted to $P^{(1)} \vee P^{(2)}$, 
are the same.
 
We have $\, |P^{(1)} \vee P^{(2)}| \le |P^{(1)}| \cdot |P^{(2)}|$. 
And for all $p^{(1)} \in P^{(1)}$ and $p^{(2)} \in P^{(2)}$,
$|p^{(1)} \vee p^{(2)}|_{\rm max} \le $
$\max\{|p^{(1)}|_{\rm max}, \, |p^{(1)}|_{\rm max}\}$.
(Recall the notation $|x|_{\rm max} = \max \{|x_i| : 1 \le i \le n\}$.)
Hence one can check in polynomial time whether  
$f^{(1)} \equiv_{\rm end} f^{(2)}$.  
 \ \ \ $\Box$

\begin{lem} \label{NONmaxExt} {\bf (non-uniqueness of maximal extensions).}
 \ There exists a right ideal morphism $f \in 2 \, {\cal RI}^{\sf fin}_2$ 
such that $f$ has {\bf two} {\em maximal extensions} in 
$2 \, {\cal RI}^{\sf fin}_2$.
\end{lem}
As a consequence, $nV$ with $n \ge 2$ cannot be defined by maximum 
extended morphisms (unlike $V$).

\medskip

\noindent {\sc Proof.} Let $f$ be defined by 
$ \, {\rm domC}(f) = {\rm imC}(f) = $
$\{(0,0), (0,1), (1,0), (1,10), (1,11)\}$, and the table

\bigskip

\begin{tabular} {c||l|l|l|l|l|}
$x$    & (0,0) & (0,1) & (1,0) & (1,10) & (1,11) \\ \hline
$f(x)$ & (0,0) & (0,1) & (1,0) & (1,11) & (1,10) 
\end{tabular} .

\bigskip

\noindent The geometric representation of $f$ is given in Fig.\ 2 (with 
mapping-by-number as in \cite{Brin1}):

\unitlength=.8mm
\begin{picture}(110,60)

\put(5,10){\framebox(40,20)}
\put(5,10){\framebox(20,40)}
\put(25,30){\framebox(20,10)}
\put(25,40){\framebox(20,10)}
\put(15,20){\makebox(0,0)[cc]{\sf 1}}
\put(35,20){\makebox(0,0)[cc]{\sf 2}}
\put(15,40){\makebox(0,0)[cc]{\sf 3}}
\put(35,35){\makebox(0,0)[cc]{\sf 4}}
\put(35,45){\makebox(0,0)[cc]{\sf 5}}

\put(60,27){\vector(1,0){20}}
\put(70,30){\makebox(0,0)[cc]{$f$}}

\put(95,10){\framebox(40,20)}
\put(95,10){\framebox(20,40)}
\put(115,30){\framebox(20,10)}
\put(115,40){\framebox(20,10)}
\put(105,20){\makebox(0,0)[cc]{\sf 1}}
\put(125,20){\makebox(0,0)[cc]{\sf 2}}
\put(105,40){\makebox(0,0)[cc]{\sf 3}}
\put(125,35){\makebox(0,0)[cc]{\sf 5}}
\put(125,45){\makebox(0,0)[cc]{\sf 4}}

\put(-2,0){\makebox(0,0)[cc]{\sf Fig. 2}}
\end{picture}
 
\bigskip

\bigskip

\noindent In Fig.\ 2, the squares labeled ``1'' and ``2'' could be merged
into one binary rectangle; alternatively, the squares labeled ``1'' and ``3'' 
could be merged into one binary rectangle.  After either step, no further 
extension is possible.  
Thus, $f$ has the following two maximal extensions $F_1$ and $F_2$:

\smallskip

\noindent (1) \ \ ${\rm domC}(F_1) = {\rm imC}(F_1) = $
$\{(\varepsilon, 0), (0,1), (1,10), (1,11)\}$, \ and

 \ \ $F_1 \, = \, \{((\varepsilon, 0), (\varepsilon, 0)), \ ((0,1), (0,1)),$
$             \ ((1,10), (1,11)), \ ((1,11), (1,10))\}$.

\smallskip

\noindent (2) \ \ ${\rm domC}(F_2) = {\rm imC}(F_2) = $
$\{(0,\varepsilon), (1,0), (1,10), (1,11)\}$, \ and

 \ \ $F_2 \, = \, \{((0,\varepsilon),(0,\varepsilon)), \ ((1,0),(1,0)),$
$             \ ((1,10), (1,11)), \ ((1,11), (1,10))\}$.
 \ \ \  \ \ \ $\Box$

\bigskip

\noindent We now give the definition of $n G_{k,1}$ and $nV$ based on 
strings.

\begin{defn} \label{BrinThompnV} {\bf (Brin-Thompson groups 
{\boldmath $nV$} and {\boldmath $n G_{k,1}$}).} 
 \ Let $A = \{0, \, \ldots, k-1\}$ and $n \ge 2$. 
The Brin-Thompson group $n G_{k,1}$ is 
$ \ n \, {\cal RI}^{\sf fin}_A/\!\! \equiv_{\rm end}$.  
 \ Equivalently, $n G_{k,1}$ is the group determined by the action of
$n \, {\cal RI}^{\sf fin}_A$ on $nA^{\omega}$.
 \ When $k = 2$ we obtain $n V$.
\end{defn}

\noindent Every element of $n G_{k,1}$ can be represented (in infinitely 
many ways) by a table of the form $F$: $P \to Q$, which is a bijection 
between two finite maximal joinless codes.

\begin{lem} \label{LEMcompositionnV} {\bf (composition in 
{\boldmath $n G_{k,1}$} based on tables).}
 \ Let $F_j:P_j \to Q_j$ be a table representing 
$f_j \in n \, {\cal RI}^{\sf fin}$, which in turn determines 
$f_j^{(\omega)} \in n G_{k,1}$ (for $j = 1, 2$).
Then the {\em composite} $\, f_2 \circ f_1$, and hence also 
$f_2^{(\omega)} \circ f_1^{(\omega)}$, is represented by the table 

\bigskip

 \ \ \  \ \ \ $(f_2 \circ f_1)|_P: \, P \to Q$, \ \  where

\medskip

 \ \ \  \ \ \ $P = f_1^{-1}(P_2 \vee Q_1)$, 

\medskip

 \ \ \  \ \ \ $Q = f_2(P_2 \vee Q_1)$.
\end{lem}
{\sc Proof.} It is a general fact about partial functions $f_2, f_1$, that 
$\, {\rm Dom}(f_2 \circ f_1) = f_1^{-1}({\rm Dom}(f_2) \cap {\rm Im}(f_1))$, 
and $\, {\rm Im}(f_2 \circ f_1) = f_2({\rm Dom}(f_2) \cap {\rm Im}(f_1))$.
Obviously, $f_2 \circ f_1 = (f_2 \circ f_1)|_{{\rm Dom}(f_2 \circ f_1)}$.

For $f_2, f_1 \in n \, {\cal RI}^{\sf fin}$, given by tables,
${\rm Dom}(f_2) \cap {\rm Im}(f_1) = (P_2 \vee Q_1) \cdot (nA^*) \,$ (by 
Lemma \ref{INTERSjoinless}). And $f_1^{-1}(P_2 \vee Q_1)$ and 
$f_2(P_2 \vee Q_1)$ are maximal joinless codes (by Lemmas \ref{InjNormal} 
and \ref{LEM_f_of_P}).
Moreover, $f_1^{-1}({\rm Dom}(f_2) \cap {\rm Im}(f_1))$ $=$ 
$f_1^{-1}(P_2 \vee Q_1) \cdot (nA^*)$, and
$\, f_2({\rm Dom}(f_2) \cap {\rm Im}(f_1))$ $=$
$f_2(P_2 \vee Q_1) \cdot (nA^*)$.
Hence, $(f_2 \circ f_1)|_{{\rm domC}(f_2 \circ f_1)}$ is given 
by the table described in the Lemma. 
 \ \ \ $\Box$

\section{The word problem of {\boldmath $n V$} is in {\sf coNP}}

For a fixed group $G$ with a fixed finite generating set $\Gamma$ the 
{\em word problem} is the following decision problem. 

{\sc Input:} A string $\,w \in (\Gamma^{\pm 1})^*$.

{\sc Question:} Does $w$ represent the identity element of $G$ ?

\medskip 

\noindent We mentioned in the Introduction that $nV$ is finitely generated 
for all $n \ge 1$. The groups $n G_{k,1}$, for $k > 2$, are presumably 
finitely generated too, but this has not been proved, so we will only 
address the word problem of $n V$ here.

\medskip

\noindent {\bf Notation.} We mostly use the alphabet 
$A = \{0, 1, \, \ldots \, , k-1\}$, usually with $k = 2$.  

For any integer $j \ge 0$, let
  \ $nA^{\le j} \, = \, \{(x_1, \ldots, x_n) \in nA^*:$
$\, |x_i| \le j$ for $i = 1, \ldots, n \}$;
for a string $w \in A^*$, $\,|w|$ denotes the length of $w$. 

\bigskip

\noindent {\bf Definition of {\sf coNP} and {\sf NP}:}  
 \ We use the following logic-based definitions of {\sf coNP} and {\sf NP} 
(see e.g., \cite{HandbookA}). Let $\Gamma$ be a finite alphabet.
A set $S \subseteq \Gamma^*$ is in {\sf coNP} \ iff \ there exists $m \ge 1$,
a two-variable predicate $\,R(.,.) \subseteq mA^* \times \Gamma^*$, and a 
polynomial $p(.)$, such that 

\smallskip

(1) \ \ $R \in {\sf P}$ \ \ (i.e., the membership problem of $R$ is in 
{\sf P}); 

\smallskip

(2) \ \ $S \, = \, \{w \in \Gamma^* \,: \ \ $
$(\forall x \in mA^{\le p(|w|)}) \, R(x, w) \, \}$.

\smallskip

\noindent Similarly, $S$ is in {\sf NP} \ iff \ for some $m \ge 1$, some
$\, R(.,.) \subseteq mA^* \times \Gamma^*$ in {\sf P}, and some polynomial
$p(.)$, 
 
\smallskip

 \ \ \ $S \, = \, \{w \in \Gamma^* \,: \ \ $
$(\exists x \in mA^{\le p(|w|)}) \, R(x, w) \, \}$.

\begin{defn} \label{DEFlength} {\bf (max length).}

\smallskip

\noindent For every $z = (z_1, \, \ldots, z_n) \in nA^*$: 
 \ \ $\ell(z) \, = \, \max\{|z_1|, \, \ldots, |z_n|\}$. 

\smallskip

\noindent For every finite set $P \subset nA^*$:
 \ \ $\ell(P)$  $\, = \,$
$\max\{\ell(z) \, : \, z \in P\}$.

\smallskip

\noindent For every $f \in n \, {\cal RI}^{\sf fin}$: 
 \ \  $\ell(f)$  $\, = \,$ 
$\max\{\ell(z) \,:\, z \in {\rm domC}(f) \,\cup\, {\rm imC}(f)\}$.
%
%
\end{defn}

\begin{pro} \label{ellAdd} {\bf (length formula).} 
 \ For all $f_2, f_1 \in n {\cal RI}^{\sf fin}$:
 \ \ $\ell(f_2 \circ f_1) \, \le \, \ell(f_2) + \ell(f_1)$.
\end{pro}
{\sc Proof.} Let $F_i: P_i \to Q_i$ be a table for $f_i$ ($i = 1, 2$). 
Recall the table for $f_2 \circ f_1$, given in Lemma \ref{LEMcompositionnV}.
We have: 

\smallskip

\noindent (L1) \hspace{1.3in} 
$\ell(P_2 \vee Q_1) = \max\{ \ell(P_2), \ \ell(Q_1)\}$. 

\medskip

\noindent Indeed, for every $p = (p_1, \, \ldots, p_n) \in P_2$, 
 \ $q = (q_1, \, \ldots, q_n) \in Q_1$, and $i \in \{1, \, \ldots, n\}$ we 
have:
 \ $|(p \vee q)_i| = \max\{|p_i|, \, |q_i|\} \,$ (by Lemma \ref{joinAn}). 
 
\noindent We have:

\smallskip

\noindent (L2) \hspace{1.3in} $\ell(f_2(P_2 \vee Q_1))$ $\, \le \,$
$\ell(Q_2) + \ell(Q_1)$ $\, \le \,$ $\ell(f_2) + \ell(f_1)$.

\medskip

\noindent Indeed, $(p \vee q)_i = \max_{\le_{\rm pref}}\{p_i, q_i\}$, 
for every $p \in P_2$, $q \in Q_1$, and $i \in \{1, \, \ldots, n\}$.
By Prop.\ \ref{INTERSjoinless}(3), $p \vee q \in {\rm Dom}(f_2)$.
Since $p$ is an initial factor of $p \vee q$ there exists $u \in nA^*$ such 
that $p u = p \vee q$. 
Since $(p \vee q)_i = \max_{\le_{\rm pref}}\{p_i, q_i\}$, the following 
holds: $\, u_i = \varepsilon$ when $(p \vee q)_i = p_i$, and 
$u_i$ is a suffix of $q_i$ when $(p \vee q)_i = q_i$.
Hence, $\ell(u) \le \ell(q)$.
Now, $f_2(p \vee q) = f_2(p) \, u$, where $f_2(p) \in Q_2$ (since 
$p \in P_2$). And $q \in Q_1$. Hence 
$\ell(f_2(p \vee q)) \le \ell(f_2(p)) + \ell(u)$ 
$\le$ $\ell(Q_2) + \ell(Q_1)$.

\noindent We also have:

\smallskip

\noindent (L3) \hspace{1.3in} $\ell(f_1^{-1}(P_2 \vee Q_1))$ $\, \le \,$ 
$\ell(P_1) + \ell(P_2)$ $\, \le \,$ $\ell(f_2) + \ell(f_1)$.

\medskip

\noindent Indeed, $f_1^{-1}$ is given by the table $f_1^{-1}|_{Q_1}$:
$Q_1 \to P_1$. Consider any $p \vee q$ for $p \in P_2$, $q \in Q_1$. 
Since $q$ is an initial factor of $p \vee q$ there exists $v \in nA^*$ such 
that $q v = p \vee q$.
Since $(p \vee q)_i = \max_{\le_{\rm pref}}\{p_i, q_i\}$, the following
holds: $\, v_i = \varepsilon$ when $(p \vee q)_i = q_i$, and
$v_i$ is a suffix of $p_i$ when $(p \vee q)_i = p_i$.
Hence, $\ell(v) \le \ell(p)$.
Now, $f_1^{-1}(p \vee q) = f_1^{-1}(q) \ v$, where $f_1^{-1}(q) \in P_1$ 
(since $q \in Q_1$). And $p \in P_2$. Hence 
$\ell(f_1^{-1}(p \vee q))$  $\le$
$\ell(f_1^{-1}(q)) + \ell(v)$ $\le$ $\ell(P_1) + \ell(P_2)$.

\smallskip

\noindent Finally, since $\, \ell(f_2 \circ f_1)$ $\, = \,$ 
$\max\{\ell(f_1^{-1}(P_2 \vee Q_1))$, $\ell(f_2(P_2 \vee Q_1))\}$, we obtain:

\medskip

\noindent (L4) \hspace{.5in}
$\, \ell(f_2 \circ f_1) \, \le \, \max\{\ell(Q_2) + \ell(Q_1), \ $ 
$\ell(P_1) + \ell(P_2)\}$ $\, \le \,$ $\ell(f_2) + \ell(f_1)$.

\smallskip

\noindent $\Box$

\begin{cor} \label{LengthPROD}
 \ Let $f_t, \, \ldots, f_1 \in n {\cal RI}^{\sf fin}$, and let 
$\lambda \in {\mathbb N}$ be such that $\ell(f_j) \le \lambda$ for 
$j = 1, \, \ldots, t$.
Then \ $\ell(f_t \circ \, \ldots \, \circ f_1) \le \lambda \, t$. 
 \ \ \  \ \ \ $\Box$
\end{cor}

\begin{lem} \label{nAcInP}
 \ For any $f \in n {\cal RI}^{\sf fin}$ and $\lambda = \ell({\rm domC}(f))$ 
we have: \ \ \ $n A^{\lambda} \subset {\rm Dom}(f)$.

\smallskip

Hence, $f|_{n A^{\lambda}}$ determines $f^{(\omega)}$.
In particular, \ $f|_{n A^{\lambda}} = {\sf id}|_{n A^{\lambda}} \,$ \ iff 
 \ $\, f^{(\omega)} = {\bf 1} \,$ in $ \, n G_{k,1}$. 
\end{lem}
{\sc Proof.} By Coroll.\ \ref{OneStepAll}(0), $n A^{\lambda}$ is a maximal 
joinless code, hence every element $p \in {\rm domC}(f)$ has a join with some
element $u \in n A^{\lambda}$. Since $\lambda = \ell({\rm domC}(f))$, $p$ is 
actually an initial factor of $u$. Hence, $u \in P \cdot (nA^*)$ 
($= {\rm Dom}(f)$). This proves that $n A^{\lambda} \subset {\rm Dom}(f)$.

For any finite maximal joinless code $P \subset {\rm Dom}(f)$, the 
restriction $f|_P$: $P \to f(P)$ is a table for $f^{(\omega)}$, hence it 
determines $f^{(\omega)}$.  
Since $n A^{\lambda}$ is a finite maximal joinless code contained in 
${\rm Dom}(f)$, the result follows.
 \ \ \ $\Box$

\begin{lem} \label{nVinCoNP}
 \ The word problem of $n V$ over any finite generating set belongs to
{\sf coNP}.
\end{lem}
{\sc Proof.} Let $\Gamma$ be any finite generating set of $n V$. 
To simplify the notation we assume that $\Gamma$ is closed under inverse, 
i.e., $\Gamma = \Gamma^{\pm 1}$.  
Every $\gamma \in \Gamma$ is represented by a table 
$F_{\gamma}$: $P_{\gamma} \to Q_{\gamma}$.
For any $w \in \Gamma^*$, let $f_w \in n {\cal RI}^{\sf fin}$ be the function 
obtained by composing the generators in $w$ (given by tables). 
Let $F_w$: $P \to Q$ be the table of $f_w$. 
By Prop.\ \ref{ellAdd} and Coroll.\ \ref{LengthPROD}: 
$\, \ell(f_w) \le c_{_{\Gamma}} \, |w|$, where
$\, c_{_{\Gamma}} = \max\{\ell(\gamma) : \gamma \in \Gamma\}$, and $|w|$
denotes the length of $w$ as a word over $\Gamma$.  So $c_{_{\Gamma}}$ is a 
known constant, determined by the finite generating set $\Gamma$.

For the word problem we have: \ $w = {\bf 1}$ in $n V$ \ iff 
 \ $f_w^{(\omega)} = {\sf id} \,$ (the identity function on $A^{\omega}$) 
 \ iff \ $\, f_w = {\sf id}|_{{\rm Dom}(f_w)}$ in $n {\cal RI}^{\sf fin}$. 
Since ${\rm domC}(f_w) = P \,$ (in the table $F_w$: $P \to Q$), we have: 
$\, f_w = {\sf id}|_{{\rm Dom}(f_w)} \,$ iff $\, P = Q \,$ and 
$\, F_w = {\sf id}|_P$.
By Prop.\ \ref{ellAdd}, $\, P \cup Q \subset$ 
$n A^{\le \, c_{_{\Gamma}} |w|}$. 
We can further restrict $f_w$ to 
$n A^{c_{_{\Gamma}} |w|} \cdot (n A^*)$; then by Lemma \ref{nAcInP}, 
we obtain the following {\sf coNP}-formula for the word problem:

\medskip

\hspace{0.4in}  $w = {\bf 1}$ in $n V$ \ \ \ iff 
 \ \ \ $(\forall x \in n A^{c_{_{\Gamma}} |w|})$ $[ \, f_w(x) = x \,]$. 

\medskip

\noindent We still need to show that the predicate $R(x,w)$, defined by 

\smallskip

 \ \ \ $R(x,w)$ \ \ $\Leftrightarrow$ 
 \ \ $[\,(\forall i \in \{1,\ldots,n\})[|x_i| = c_{_{\Gamma}} |w|]$ 
$ \ \Rightarrow \ $  $f_w(x) = x \,]$,

\smallskip

\noindent belongs to {\sf P}. 
I.e., we want a deterministic polynomial-time algorithm that on input 
$w \in \Gamma^*$ and $x \in n A^{c_{_{\Gamma}} |w|}$, checks whether 
$f_w(x) = x$.
To do this we apply, to $x \in n A^{c_{_{\Gamma}} |w|}$, the tables of the 
generators $\gamma_j \in \Gamma$ that appear in 
$w = \gamma_t \, \,\ldots \, \gamma_1$.  We compute 
$\, x \longmapsto \gamma_1(x) = y^{(1)}$  $\longmapsto$  
$\gamma_2(y^{(1)}) = y^{(2)}$  $\longmapsto$  $ \ \ \ldots \ \ $
$\longmapsto$  $\gamma_t(y^{(t-1)}) = y^{(t)} = f_w(x)$.
Since $x \in n A^{c_{_{\Gamma}} |w|}$ $\subseteq$ ${\rm Dom}(f_w)$, every 
$y^{(j)}$ is defined.  Moreover, 
$x = p u$ for some $p \in P$ and $u \in n A^{\le c_{_{\Gamma}} |w|}$.
By Prop.\ \ref{ellAdd}, 
$\, |y^{(j)}| \le \, n \, c_{_{\Gamma}} \, j + |u| \,$  $\le$ 
$\, 2n \, c_{_{\Gamma}} \, |w|$. After computing $y^{(t)}$ we check whether 
$y^{(t)} = x$.

The application of the table of $\gamma_j$ to $y^{(j-1)}$ takes time 
proportional to $\,|y^{(j-1)}| \,$ (for $j = 1, \, \ldots, t$). 
So, the time complexity of verifying whether $x$ and $w$ satisfy the 
predicate is (up to a constant multiple) $ \ \le \ $  
$|x| \,+\, \sum_{j=1}^t |y^{(j)}|$ 
$\,\le\,$ $n\,c_{_{\Gamma}}\, |w| + |w| \cdot 2\,n\, c_{_{\Gamma}}\, |w|$. 
Hence the time-complexity of the predicate is quadratic in $|w|$.
 \ \ \  \ \ \ $\Box$

\section{{\sf coNP}-completeness of the word problem of {\boldmath $nV$} }

In this section we prove that the word problem of $nV$ with $n \ge 2$, over 
any finite generating set, is {\sf coNP}-hard with respect to polynomial-time
many-one reduction. The result for all $nV, n \ge 2$, follows quickly from 
the result for $2V$.
We proved already in Lemma \ref{nVinCoNP} that the word problem of $nV$ 
belongs to {\sf coNP}.

For $2V$, {\sf coNP}-hardness of the word problem follows fairly directly
from the {\sf coNP}-hardness of the word problem of $V$ over the infinite 
generating set $\Gamma_{\! V} \cup \tau$, by making use of the shift 
$\sigma\,$ (subsection 4.6). 
Here, $\Gamma_{\! V}$ is any finite generating set of $V$ and $\tau$ is the 
set of position transpositions $\{\tau_{i,i+1} : i \ge 1\}$. 

At the end of subsection 4.1 we show that the word problem of $V$ over
$\Gamma_{\! V} \cup \tau$ belongs to {\sf coNP}. 
The main difficulty is to prove that the word problem of $V$ over 
$\Gamma_{\! V} \cup \tau$ is {\sf coNP}-hard; this is proved in subsections 
4.2 - 4.4, by constructing a binary conjunctive polynomial-time reduction 
of the circuit equivalence problem to this word problem.
An alternative proof, that gives a polynomial-time many-one reduction, appears
in subsection 4.5.

\subsection{Preliminaries on the word problem and complexity}

We give some definitions and facts about complexity and the word problem of 
a group, especially when an infinite generating set is used. 
Here we use finite and infinite alphabets (but we always point out when an
alphabet is infinite).

\begin{defn} \label{DEFpolyconjred}
 \ Let $\Sigma_1, \Sigma_2$ be two finite alphabets, and let $m$ be a 
positive integer. 
A polynomial-time {\em conjunctive reduction} of arity $m$ from
$L_1 \subseteq \Sigma_1^*$ to $L_2 \subseteq \Sigma_2^*$ is a 
polynomial-time computable total function 
$\, \rho: \Sigma_1^* \, \to \, m \Sigma_2^*$ 
such that for all $x \in \Sigma_1^*$:

\medskip

 \ \ \  \ \ \    $x \in L_1$ \ \ iff 
 \ \ $\rho(x) \, \in \ {\large \sf X}_{_{j=1}}^{^{m}} L_2$
  \ ($= m L_2$).

\medskip

\noindent 
Equivalently, $L_1 = \rho^{-1}({\large \sf X}_{_{j=1}}^{^{m}} L_2)$.
In other words, $\rho$ reduces the problem $L_1$ to $m$ instances of the 
problem $L_2$, and the $m$ answers are combined by ``and''. 

A polynomial-time conjunctive reduction of arity 1 is called a 
{\em many-one} reduction.
\end{defn}
The reductions in Def.\ \ref{DEFpolyconjred} are a very special case of 
polynomial-time truth-table reductions; see e.g.\ \cite[Def.\ 7.18]{DuKo}. 
It is straightforward to show that each of {\sf P}, {\sf NP}, and 
{\sf coNP}, is closed under downward polynomial-time conjunctive reduction 
of bounded arity. 

\smallskip

In this paper we use the following definition of {\sf coNP}-hardness and
{\sf coNP}-completeness.

\begin{defn} \label{DEFhardness} 
 \ Let $\Sigma_0$ be a finite alphabet. 
A problem $L_0 \subseteq \Sigma_0^*$ is {\sf coNP}{\em -hard} iff for 
every finite alphabet $\Sigma$ and every problem $L \subseteq \Sigma^*$
there exists a polynomial-time conjunctive reduction $\rho$ of bounded
arity that reduces $L$ to $L_0$.

Moreover, $L_0 \subseteq \Sigma_0^*$ is {\sf coNP}{\em -complete} iff $L_0$
is {\sf coNP}-hard and $L_0$ belongs to {\sf coNP}.
\end{defn}
There are many well-known {\sf coNP}-complete problems, e.g., the tautology
problem for boolean formulas, the integer linear programming equivalence
problem, the 4-coloring problem, the connectivity lower-bound problem (see 
e.g.\ \cite{GaryJohns}, \cite[Introduction]{BiCoNP}).  
We will use the equivalence problem for acyclic boolean circuits (defined 
in Section 4.2).

\begin{defn} \label{DETequalinG}
 \ Let $G$ be a group, and let $\Gamma$ ($\,\subseteq G$) be a (possibly
infinite) generating set for $G$. For words
$w_1, w_2 \in (\Gamma^{\pm 1})^*$ we say that {\em ``$w_1 = w_2$ in $G$''}
iff the generator sequences $w_1$ and $w_2$ have the same value when their
elements are multiplied in $G$.  In a similar way, for
$w \in (\Gamma^{\pm 1})^*$ and $g \in G$, we say {\em ``$w = g$ in $G$''}
iff $g$ is the value obtained when the elements of $w$ are multiplied in
$G$. We also use the notation $w_1 =_G w_2$ or $w =_G g$ for this. 
\end{defn}
To simplify the notation, we will from now on take group generating sets
$\Gamma$ that are closed under inverse, i.e., $\Gamma = \Gamma^{\pm 1}$.

\begin{lem} \label{ComplexSub} {\bf (folklore).} 
 \ Let $G_2$ be a finitely generated subgroup of a finitely generated 
group $G_1$, and let $\Gamma_i$ be a finite generating set of $G_i$ for
$i = 1,2$.

\smallskip

\noindent {\rm (1)} \ If the word problem of $G_1$ over $\Gamma_1$ is 
decidable in deterministic (or nondeterministic, or co-nondeterministic) 
time $\le t_1(.)$, then the word problem of $G_2$ over $\Gamma_2$ is 
decidable in deterministic (respectively in nondeterministic, or
co-nondeterministic) time $\le t_1(O(.))$.

\smallskip

\noindent
{\rm (2)} \ If a problem $L \subseteq \Sigma^*$ (where $\Sigma$ is finite) 
is reducible to the word problem of $G_2$ over $\Gamma_2$ by a 
polynomial-time conjunctive reduction of arity $m$, then $L$ is also 
reducible to the word problem of $G_1$ over $\Gamma_1$ by a polynomial-time conjunctive reduction of arity $m$.
\end{lem}
{\sc Proof.}  
To simplify the notation, let us assume that $\Gamma_1$ and $\Gamma_2$ are 
closed under inverse.

\noindent (1) Since $G_2 \subseteq G_1$, for every generator 
$\gamma \in \Gamma_2$ there exists a word $w_{\gamma} \in \Gamma_1^{ \ *}$ 
such that $\, \gamma =_{G_1} w_{\gamma}$. Then the total function 

\smallskip

 \ \ \ $\rho_{2,1}: \ x_1 \, \ldots x_n \in \Gamma_2^{ \ *}$ 
$ \ \longmapsto \ $ 
$w_{x_1} \cdot \, \ldots \, \cdot w_{x_n} \in \Gamma_1^{ \ *}$ 

\smallskip

\noindent is a one-to-one linear-time reduction of the word problem of $G_2$
over $\Gamma_2$ to the word problem of $G_1$ over $\Gamma_1$; here, 
``$\cdot$'' denotes concatenation.
The length of $w_{x_1} \cdot \ldots \cdot  w_{x_n}$ is 
$\, \le c \, |x_1 \ldots x_n|$, where 
$\, c = \max\{|w_{\gamma}| : \gamma \in \Gamma_2\}$. Hence, if the word 
problem of $G_1$ over $\Gamma_1$ has time-complexity $\le t_1(.)$, then the
word problem of $G_2$ over $\Gamma_2$ has time-complexity $\le t_1(c n)$ for
inputs of length $n$. 
 
\smallskip

\noindent (2) For every $\gamma \in \Gamma_2$ let 
$w_{\gamma} \in \Gamma_1^{\, *}$ be such that $\gamma = w_{\gamma}$ 
in $G_1$, and let $W(.)$ be the free monoid homomorphism from 
$\Gamma_2^{\, *}$ into $\Gamma_1^{\, *}$ determined by 
$W(\gamma) = w_{\gamma}$ for all $\gamma \in \Gamma_2$. Let 
$\, \rho: \Sigma^* \to m \, \Gamma_2^{\, *}$ be a polynomial-time 
conjunctive reduction of arity $m$, such that for all $v \in \Sigma^*$: 
$v \in L$ iff $\rho(v) = (\varepsilon)^m \,$.
Then 

\smallskip

$\, x \in \Sigma^* \ \longmapsto \ \rho(x) = (y_1, \ldots, y_m) \in $
$m \, \Gamma_2^{\, *}$ \ $\longmapsto$ \ $(W(y_1), \ldots, W(y_m)) \in $
$m \, \Gamma_1^{\, *} \,$ 

\smallskip

\noindent is a polynomial-time conjunctive reduction of arity $m$, from 
$L$ to the he word problem of $G_1$ over $\Gamma_1$. 
 \ \ \ $\Box$

\bigskip

For the word problem of groups, {\em infinite generating sets} cannot always
be avoided, because some groups are not finitely generated, and because 
some finitely generated groups have interesting infinite generating sets. 
In order to apply the concepts of decidability or computational complexity
to groups with infinite generating sets, we encode countable generating sets
over a finite alphabet.
We will use the following.

\begin{defn} \label{DEFencoding} {\bf (encoding).} 

\noindent An {\em encoding} of a countable set $\Gamma$ is an injective 
total function $\, {\sf code}: \Gamma \to \{0,1\}^*$ such that 
$\, {\sf code}(\Gamma) \,$ is a prefix code that is accepted by a 
finite-state automaton. 

For a word of generators $\, w = w_1 \, \ldots \, w_m \in \Gamma^*$ 
we define ${\sf code}(w)$ by the concatenation \ ${\sf code}(w) \ $ $=$ 
$ \ {\sf code}(w_1) \cdot \, \ldots \, \cdot {\sf code}(w_m)$. Hence, 
${\rm Im}({\sf code(.)}) = {\sf code}(\Gamma^*) = ({\sf code}(\Gamma))^*$, 
which is a finite-state language.
\end{defn}
Since the function {\sf code} is injective it has an inverse function,
${\sf code}^{-1}$, whose domain is ${\sf code}(\Gamma^*)$. 
Every countable set admits an encoding of the above type
(e.g., with image set $\, 0^* 1 = \{0^n 1: n \in \omega\}\,$).

\bigskip

\noindent The {\em word problem for a group $G$ with an infinite generating 
set $\Gamma$ and encoding} ${\sf code}: \Gamma \to \{0,1\}^*$ is specified 
as follows.

{\sc Input:} \ $x \in \{0,1\}^*$.

{\sc Precondition:} \ $x \in {\sf code}(\Gamma)^*$. \ (Since 
${\sf code}(\Gamma)^*$ is finite-state, the precondition is easy to check.)

{\sc Question:} \ ${\sf code}^{-1}(x) = {\bf 1}$ in $G\,$?
 \ \ (Here, {\bf 1} denotes the identity element of $G$.)

\medskip

\noindent Equivalently, the word problem is the membership problem of the 
language 

\smallskip

 \ \ \   \ \ \   
$\, {\rm WP}_{G, \Gamma, {\sf code}} \, = \, $
$\{x \in \{0,1\}^* : \ {\sf code}^{-1}(x) = {\bf 1} \ {\rm in} \ G\}$.

\smallskip

\noindent From now on,  by {\em complexity of the word problem} of $G$ over
$\Gamma$ we mean the complexity of ${\rm WP}_{G, \Gamma, {\sf code}}$;
note that the problem depends on $G$, $\Gamma$, and {\sf code}.

\begin{lem} \label{INFINComplexSub}
  \ Let $G_2$ be a subgroup of a countable group $G_1$, let $\Gamma_i \,$ 
$( \subseteq G_i$) be a countable generating set of $G_i$, and let 
${\sf code}_i(.)$ be an encoding of $\Gamma_i$, for $i = 1,2$. We also 
assume that there is a total function $h: \Gamma_2 \to \Gamma_1^{\, *}$ with 
the following properties (that connect the encodings ${\sf code}_1$ and 
${\sf code}_2$):

\smallskip

$\bullet$ \ For all $\gamma \in \Gamma_2$: \ $\gamma = h(\gamma)\,$ 
{\rm in} $G_1$;

\smallskip

$\bullet$ \ the function
 \ \ $h_0: \, {\sf code}_2(\gamma) \in {\sf code}(\Gamma_2)$ 
$ \ \longmapsto \ $
${\sf code}_1(h(\gamma)) \in  {\sf code}(\Gamma_1)^*$ \ is computable in

 \ \ \ linear time.

\smallskip

\noindent The function $h$ is extended to a free-monoid homomorphism
$\Gamma_2^{\, *} \to \Gamma_1^{\, *}\,$ that will also be called $h$; so for 
every $w \in \Gamma_2^{\, *}$ we have: \ $w = h(w)$ {\rm in} $G_1$.

\medskip

\noindent Then the following hold:

\smallskip

\noindent {\bf (1)} \ If $\, {\rm WP}_{G_1, \Gamma_1, {\sf code}_1}$ is in 
${\sf DTime}(t)$ (or in ${\sf NTime}(t)$, or in ${\sf coNTime}(t)$), then 
$\, {\rm WP}_{G_2, \Gamma_2, {\sf code}_2}$ is in ${\sf DTime}(t(O(.)))$
(respectively in ${\sf NTime}(t(O(.)))$, or ${\sf coNTime}(t(O(.)))$).

\smallskip

\noindent {\bf (2)} \ If $L \subseteq \Sigma^*$ (where $\Sigma$ is finite) 
is reducible to $\, {\rm WP}_{G_2, \Gamma_2, {\sf code}_2}$ by a 
polynomial-time conjunctive reduction of arity $m$, then $L$ is also 
reducible to $\, {\rm WP}_{G_1, \Gamma_1, {\sf code}_1}$ by a 
polynomial-time conjunctive reduction of arity $m$. 
Hence, if $\, {\rm WP}_{G_2, \Gamma_2, {\sf code}_2}$ is hard for a 
complexity class (e.g., {\sf coNP}), then 
$\, {\rm WP}_{G_1, \Gamma_1, {\sf code}_1}$ is also hard for that complexity 
class.
\end{lem}
The functions $h$ and $h_0$ in the Lemma have the commuting diagram 

\smallskip

\hspace{1.1in} $h_0 \circ {\sf code}_2(.) = {\sf code}_1 \circ h(.)\,$; 

\smallskip

\noindent equivalently,
$ \, h_0 = {\sf code}_1 \circ h(.) \circ {\sf code}_2^{-1}(.)$, and
$ \, h(.) = {\sf code}_1^{-1} \circ h_0 \circ {\sf code}_2(.)$.
Note that in general $h$ cannot be viewed as a computable function (as 
opposed to $h_0$), since its domain and image are arbitrary countable sets.

\medskip

\noindent
{\sc Proof.}  
(1) For any $w \in \Gamma_2^{ \ *}$ we have: $w = {\bf 1}$ in $G_2$ over 
$\Gamma_2$ \ iff \ $h(w) = {\bf 1}$ in $G_1$ over $\Gamma_1$.
Therefore, $x = {\sf code}_2(w) \in {\rm WP}_{G_2, \Gamma_2, {\sf code}_2}$ 
 \ iff \ $h_0(x) = {\sf code}_1(h(w))$   $\in$  
${\rm WP}_{G_1, \Gamma_1, {\sf code}_1}$. 

Thus we have the following algorithm for the membership problem of
${\rm WP}_{G_2, \Gamma_2, {\sf code}_2}$ on input $x \in \{0,1\}^*$:
First, check whether $x \in {\sf code}_2(\Gamma_2)^*$; this can be checked
in linear time, since ${\sf code}_2(\Gamma_2)^*$ is finite-state.
Second, compute $h_0(x)$ (in linear time). Finally, check whether $h_0(x)$
is in ${\rm WP}_{G_1, \Gamma_1, {\sf code}_1}$, in time 
$\, \le t(|h_0(x)|) \le t(O(|x|))$.

(2) Let $\rho: x \in \{0,1\}^* \longmapsto $
$(y_1, \, \ldots, y_m) \in m \, \{0,1\}^*$ be a polynomial-time conjunctive 
reduction of arity $m$ from $L$ to ${\rm WP}_{G_2, \Gamma_2, {\sf code}_2}$.
Hence, $x \in L$ \ iff \ $\{y_1, \, \ldots, y_m\} \subset$
${\rm WP}_{G_2, \Gamma_2, {\sf code}_2}$.
By the definition of $h$ and $h_0$, the latter holds iff 
 \ $\{h_0(y_1), \, \ldots \, , h_0(y_m)\} \subset$
${\rm WP}_{G_1, \Gamma_1, {\sf code}_1}$. Thus the function
$x \longmapsto (h_0(y_1), \, \ldots \, , h_0(y_m))$, where 
$(y_1, \, \ldots, y_m) = \rho(x)$, is a polynomial-time conjunctive reduction 
of arity $m$ from $L$ to ${\rm WP}_{G_1, \Gamma_1, {\sf code}_1}$.
 \ \ \ $\Box$

\bigskip

\noindent {\bf Some conventions and a fact about the Thompson group 
{\boldmath $V$ over $\Gamma_{\! V} \cup \tau$}:}

\smallskip

\noindent 
We pick a finite generating set $\Gamma_{\! V}$ for $V$, and for notational
convenience we will assume that $\, \Gamma_{\! V} = \Gamma_{\! V}^{\pm 1}$.
We also use the set of bit position transpositions $\tau =$
$\{\tau_{j,j+1} : j \ge 2\}$. 
We assume $\, \Gamma_{\! V} \, \cap \, \tau \,$ $=$ $\varnothing$.

\begin{defn} \label{DEFsizeG21} {\bf (size of a generator).} 
 \ For any generator $\delta \in \Gamma_{\! V} \, \cup \, \tau \,$ we define
the {\em size} $\|\delta\|$ as follows:
For $\delta = \gamma \in \Gamma_{\! V}$ we let $\, \|\gamma\| = 1$, and for
$\delta = \tau_{j,j+1} \in \tau$ we let $\, \|\tau_{j,j+1}\| = j+1$.
For a string of generators $w = w_m \, \ldots \, w_1 \,$ with
$w_i \in \Gamma_{\! V} \, \cup \, \tau \,$ for $i = 1, \, \ldots, m$,
the {\em size} of $w$ is defined by \ $\|w\| = \sum_{i=1}^m \|w_i\|$.

For the word $w$ as above, the {\em length} of $w$ is $|w| = m$.
\end{defn}

\begin{lem} \label{G21incoNP}
 \ The word problem of $V$ over $\Gamma_{\! V} \cup \tau$ belongs to 
{\sf coNP}.
\end{lem}
{\sc Proof.} We have $\ell(\tau_{j,j+1}) = j+1 = \|\tau_{j,j+1}\|\,$ 
(where $\ell(.)$ was defined in Def.\ \ref{DEFlength} based on tables of 
elements of $n \, {\cal RI}^{\sf fin}$).
And there is a positive integer constant $c$ such that for all 
$\gamma \in \Gamma_{\! V}$: $\,\ell(\gamma) \le c$. Hence, for any $w \in$
$(\Gamma_{\! V} \cup \tau)^*$ we have (by Prop.\ \ref{ellAdd}):
  \ $\ell(w) \le c \, \|w\|$.

Now the proof of Lemma \ref{nVinCoNP} can be applied. For any $v \in$
$(\Gamma_{\! V} \cup \tau)^*$, let $f_v \in {\cal RI}^{\sf fin}$ be the
right ideal morphism of $\{0,1\}^*$ generated by $v$. Then for every
$\, w \in (\Gamma_{\! V} \cup \tau)^*$ we have:

\smallskip

 \ \ \  \ \ \ $w = {\bf 1} \,$ in $V$ \ \ \ iff
 \ \ \ $(\forall x \in \{0,1\}^{c \,\|w\|}) \, [ \, f_w(x) = x \, ]$.

\smallskip

\noindent The predicate $R(x,w)$, defined by
$\, [\, x \in \{0,1\}^{c \,\|w\|} \, \Rightarrow \, f_w(x) = x \,]$,
is in {\sf P}. This uses the same proof as Lemma \ref{nVinCoNP}, and the
fact that $w$ is encoded over $\{0,1\}^*$ in such a way that $\ell(w)$,
$\|w\|$, and the length of the encoding, are linearly related.
Hence the above $\forall$-formula is a {\sf coNP}-formula for the word
problem of $V$ over $\Gamma_{\! V} \cup \tau$.
 \ \ \ $\Box$

\bigskip

\bigskip

\noindent {\bf Outline of the proof of {\sf coNP}-hardness of the word
problem of {\boldmath $V$} over {\boldmath $\Gamma_{\! V} \cup \tau$}:}

\smallskip

\noindent In subsections 4.2 - 4.4 we follow (a part of) the strategy of 
\cite{BiCoNP}, where another finitely presented group with 
{\sf coNP}-complete word problem was constructed.

\noindent
1. Every acyclic boolean circuit $C$ is ``simulated'' by a element of $V$,
represented by a word $w_C$ over $\Gamma_{\! V} \cup \tau$, such that the
size of $w_C$ is polynomially bounded by the size of $C$ (subsection 4.2,
Def.\ \ref{DEFsimulat} and Theorem \ref{reduction}).

\noindent
2. The equivalence problem for acyclic boolean circuits is reduced (by a
polynomial-time one-one reduction) to the generalized word problem of the
subgroup ${\rm pFix}_V(0)$ in $V$ (subsection 4.3, Coroll.\ \ref{GenWP}).

\noindent
3. Thanks to the ``commutation test'', the generalized word problem of
${\rm pFix}_V(0)$ in $V$ is reduced to two instances of the word problem of
$V$ over $\Gamma_{\! V} \cup \tau$ (subsection 4.4, Lemma \ref{CommTest01}).
This reduction is a 2-ary conjunctive linear-time reduction (``2'' comes
from the fact that $V$ is 2-generated).

\subsection{Circuits and the Thompson group {\boldmath $V$}}

Our first step in the proof of {\sf coNP}-hardness is to represent acyclic 
boolean circuits by words over the generating set $\Gamma_{\! V} \cup \tau$ 
of $V$. 

An {\em acyclic boolean circuit} is specified by a directed acyclic graph
({\sc dag}) without isolated vertices, together with a vertex labeling. 
This labeling associates 
(1) an input variable with each source vertex, 
(2) an output variable with each sink vertex, and
(3) a gate (of type {\sc not}, {\sc fork}, {\sc and}, or {\sc or}) with 
each interior vertex. 
By definition, a source vertex is a vertex of in-degree 0; a sink vertex  
of out-degree 0; an interior vertex is a vertex whose in-degree and 
out-degree are both non-zero. 
A source vertex is also called {\em input port}, and a sink vertex is also 
called {\em output port}.

A {\em gate} is, by definition, a total function $\{0,1\}^m \to \{0,1\}^n$ 
(for some $m, n \ge 1$). We consider the following four types of gates, 
where $u \in \{0,1\}^{j-1}$, $\,x_j, x_{j+1} \in \{0,1\}$, and 
$v \in \{0,1\}^{n-j} \,\cup\, \{0,1\}^{n-j-1}$.

\smallskip

 \ \ \ {\sc not}$_j: \ u \, x_j \, v \, \mapsto \, u \, \overline{x}_j v$;
 \ here, $m \ge 1$ and $j \le n = m$;

\smallskip

 \ \ \ {\sc and}$_{j,j+1}: \ u \, x_j x_{j+1} \, v \, \mapsto \, $
    $u \, (x_j \, \& \, x_{j+1}) \, v$; \ here, $m \ge 2$ and 
    $j \le m-1 = n$;

\smallskip

 \ \ \ {\sc or}$_{j, j+1}: \ u \, x_j x_{j+1} \, v \, \mapsto$
     $\, u \, (x_j \, {\sf or} \, x_{j+1}) \, v$; \ here, $m \ge 2$ and 
     $j \le m-1 = n$;

\smallskip

 \ \ \ {\sc fork}$_j: \ u \, x_j \, v \, \mapsto \,u \, x_j \, x_j  \, v$;
  \ here, $1 \le j \le m$, and $n = m+1$.

\smallskip

\noindent The operation {\sc fork}$_j$ makes an extra copy of $x_j$. 
In traditional circuit theory, {\sc fork}s are not used separately;
instead, {\sc not}, {\sc and}, and {\sc or} are allowed to produce several 
copies of the output bit. However, using {\sc fork} as a separate gate 
simplifies the conversion of a circuit into a sequence of functions. 
We also use the {\em wire-crossing} operation, which swaps the ``wires'' 
$i$ and $j$ (where $1 \le i < j \le m$); this is the function 

\smallskip

 \ \ \ $\tau_{i,j}: \ u \, x_i \, v \, x_j \, w \, \longmapsto \,$
$u \, x_j \, v \, x_i  \, w$, 

\smallskip

\noindent where, $u \in \{0,1\}^{i-1}$, $v \in \{0,1\}^{j-i}$, 
$w \in \{0,1\}^{m-j-1}$, $m \ge 2$, and $n = m$.  
This operation is not a gate; it is not associated with a vertex, but 
follows from the incidence relation of the graph.

\smallskip

Note that all the gates {\sc not}$_j$ are different functions for different 
values of $j$; the same applies to all {\sc and}$_{j,j+1}$, and all 
{\sc or}$_{j,j+1}$.
However in the presence of the operations $\tau_{i,j}$ it is sufficient to 
use just one set of gates \{{\sc not}, {\sc and}, {\sc or}, {\sc fork}\}, 
applied to bit positions 1, or 1 and 2. E.g., \ {\sc not}$_i$ $=$ 
$\tau_{i,1}$ $\circ$ {\sc not}$_1$ $\circ$ $\tau_{i,1}$. 
Thus, here we view acyclic circuits as expressions over the generating set 
 \ \{{\sc not}, {\sc and}, {\sc or}, {\sc fork}\}
$\cup$ $\{\tau_{i,j} : j > i \ge 1\}$.
Note that $\tau_{i,j} \in V$, with ${\rm domC}(\tau_{i,j})$  $=$  
${\rm imC}(\tau_{i,j})$  $=$  $\{0,1\}^j$.

An acyclic circuit $C$ with sequence of input variables 
$(x_1, \, \ldots \, , x_m)\,$ (with values ranging over $\{0,1\}^m$), and 
sequence of output variables $(y_1, , \ldots \, , y_n)\,$ (with values in 
$\{0,1\}^n$), determines an {\em input-output function} 
$\, f_C: \{0,1\}^m \to \{0,1\}^n$; this is a total function. 
Any total function of the form $\,F: \{0,1\}^m \to \{0,1\}^n$ is called a 
{\em boolean function}. 
In circuit theory it is proved that for every boolean function $F$ there 
exists an acyclic circuit whose input-output function is $F$; see 
e.g.\ \cite{HandbookA, Wegener, JESavage, FredToff}.

\smallskip

Two circuits $C_1$ and $C_2$ are called {\em equivalent} iff
$\, f_{C_1} = f_{C_2}$.

The equivalence problem for acyclic boolean circuits (in short, the {\em
circuit equivalence problem}) is specified as follows:

 \ \ \ {\sc Input:} \ $C_1$, $C_2\,$ (two circuits, described by {\sc dag}s 
  with gate labels on the vertices);

 \ \ \ {\sc Question:} \ $f_{C_1} = f_{C_2}$ ?

\smallskip

In order to consider the complexity of problems about circuits we need to
define the {\it size} of an acyclic boolean circuit $C$, denoted by $|C|$, 
and simply called {\bf circuit size}; it is defined as follows: 
If $C$ has $k_1$ gates of type {\sc not} or {\sc fork}, $k_2$ gates of type 
{\sc and} or {\sc or}, and $n$ output variables, then the size of $C$ is 
defined to be $ \, |C| = k_1 + 2 \cdot k_2 + n$.
Equivalently, $|C|$ {\em is the number of edges} (or {\em wires}) between 
gates, or from an input to a gate, or from a gate to an output (for that 
reason, gates with two input variables are counted twice).

\medskip

\noindent {\bf Remarks concerning circuit definitions:} Acyclic circuits and 
their sizes are defined in a variety of ways in the literature 
\cite{HandbookA, Papadim, JESavage, Wegener, GaryJohns, HU, DuKo, Jordan}; 
however, all these definitions lead to sizes that are {\em polynomially 
equivalent} (i.e., each one is polynomially bounded in terms of every other 
one). In the theory of {\sf NP}- or {\sf coNP}-completeness, polynomial 
differences are not significant.   \\   
(1) In the literature, the circuit size is usually defined as the number 
of vertices. Since we do not use isolated vertices in a circuit, we have 
$n_V \le n_E \le n_V^{ \ 2}\,$ (where $n_V$ and $n_E$ denote the number of
vertices and edges). So $n_V$ and $n_E$ are polynomially equivalent.  \\   
(2) In the literature the input and output variables are usually not called
vertices, but in that case they are nevertheless counted among the vertices 
in the definition of circuit size. \\  
(3) When a circuit is described by a bitstring $s_C$, the length satisfies
$n_E \le |s_C| \le c \, n_E \, \log_2 n_V$, for some constant $c \ge 1$. 
Typically, such a description of $C$ lists all the edges, where each edge is
given as a pair of strings (the names of two vertices, each vertex name 
having length $\le 1 + \log_2 n_V$).  An additional list is given that 
associates a gate or an input variable or an output variable with each 
vertex. An input variable $x_i$ is described by a code word (for $x$) and 
the binary representation of $i$; the output variables $y_j$ are described
similarly. In any case, $|s_C|$ and $n_E$ are polynomially equivalent. \\   
(4) In the literature, the {\sc fork}-gate is usually not used explicitly; 
instead, the {\sc and}-, {\sc or}-, and {\sc not}- gates, as well as the input
variables, are allowed to have a {\em fan-out}. However, even in that case,
every wire goes to a gate or an output variable, so the total of all be 
fan-outs is $\le n_V^{ \ 2}$. A gate with fan-out $k$ can be replaced 
by a gate with fan-out 1 and $k-1$ {\sc fork}-gates. This leads to a
circuit with gates that have fan-out 1, except for {\sc fork}-gates with 
fan-out 2; the size increase is polynomially bounded. \\   
(5) In the literature, {\sc and} and {\sc or}-gates are allowed to have a
{\em fan-in} $\ge 2$. But every fan-in wire comes from a gate of an input
variable, so the total of all be fan-ins is $\le n_V^{ \ 2}$. An
{\sc or}-gate with fan-in $k$ can be replaced by $k-1$ {\sc or}-gates with
fan-in 2 (and similarly for {\sc and}). This leads to a circuit with gates 
that have fan-in $\le 2$; the size increase is polynomially bounded. 

\medskip

\noindent {\bf Remarks on complexity:} The circuit equivalence problem is a 
well-known problem that is {\sf coNP}-complete. 
It is fairly straightforward to prove that the problem is in {\sf coNP}. 
Moreover, the tautology problem for boolean formulas (which is a classical 
{\sf coNP}-complete problem) is a special case of the circuit equivalence 
problem (and is reduced to the circuit equivalence problem by converting a 
boolean formula into a circuit and asking whether a given circuit is 
equivalent to a circuit for the constant-1 function). 
See \cite[Introduction]{Jordan} for comments on the circuit equivalence 
problem, see \cite{Papadim} for a circuit-based proof of 
{\sf NP}-completeness of the satisfiability problem for boolean formulas, 
and see \cite{HandbookA, GaryJohns} for general information.

\smallskip

The following well-known fact implies that every $\tau_{i,j}$ can be 
expressed as a composition of elements of 
$\tau = \{\tau_{k,k+1} : k \ge 1\}$; the expression has linear length in 
terms of $j$.

\begin{lem} \label{classicalFormulas} \
As elements of $V$ the transpositions satisfy 

\smallskip

 \ \ \ $\tau_{i,j} \ = \ $
 $\tau_{i,i+1} \ \tau_{i+1,i+2} \ \ldots \ \tau_{j-2,j-1} \ \tau_{j-1,j}$
 $ \ \tau_{j-2,j-1} \ \ldots \ \tau_{i+1,i+2} \ \tau_{i,i+1}$ , \ \ \
 if $1 \le i < j$.

\smallskip

\noindent The word length of $\tau_{i,j}$ over $\tau$ is therefore 
$ \ \le 2(j-i)-1$.  \ \ \   \ \ \ $\Box$ 
\end{lem}

We want to represent the circuit gates {\sc not}, {\sc or}, {\sc and},
and {\sc fork}, by elements of $V$. For this, the main problem is that
the input-output function of a circuit is not necessarily a permutation.
Therefore we introduce the following notion of ``simulation'' of a circuit 
$C$ by a Thompson group element $\Phi_C$ and by a word $w_C$ over 
$\Gamma_{\! V} \cup \tau$ (Def.\ \ref{DEFsimulat} and Theorem 
\ref{reduction} below). See the discussion in \cite{BiCoNP} for additional
motivation of our definition of simulation.

\begin{defn} \label{DEFsimulat} {\bf (simulation).}
 \ Let $\, f: \{0,1\}^m \to \{0,1\}^n \,$ be a total function. An element 
$\Phi_f \in V$ {\em simulates} $f$ \ iff \ for all 
$\, x \in \{0,1\}^m$:
 \ \ \ $\Phi_f(0 \, x) \ = \ 0 \ f(x) \ x$.

\medskip

\noindent When $\Phi_f$ is represented by a {\em word} $w_f \in$
$(\Gamma_{\! V} \cup \tau)^*$ \ we say that $w_f$ simulates $f$.
\end{defn}
According to this definition, $f$ is faithfully described by the action of 
$\Phi_f$ on $0\,\{0,1\}^*$; but there are no constraints on the values of 
$\Phi_f$ for input strings in $1 \, \{0,1\}^*$.
Since $\Phi_f$ is an element of $V$ it is a bijection between finite
maximal prefix codes, whereas $f$ need not be injective nor surjective.
So there has to be a big difference between $\Phi_f$ and $f$ somewhere.
In subsections 4.3 and 4.4 we show that, nevertheless, the equivalence 
problem of circuits can be reduced to the word problem of $V$ over
$\Gamma_{\! V} \cup \tau$. In the rest of this subsection we construct 
$\Phi_f$.

The next Lemma follows immediately from the definition of simulation.

\begin{lem} \label{simulatVSfix} \
Let $f$ and $g$ be any boolean functions with the same number of input
variables and the same number of output variables.
If $f$ and $g$ are simulated by $\Phi_f$, respectively $\Phi_g$, then
we have

\medskip

 \ \ \ \ \    $f = g$ \ \ \  iff \ \ \
    $(\Phi_f)|_{0 \{0,1\}^*} \ = \ (\Phi_g)|_{0 \{0,1\}^*}$ .
 \ \ \  \ \ \ $\Box$
\end{lem}
We choose the following elements of $V$ to describe the gates {\sc not}, 
{\sc or}, {\sc and}, and {\sc fork}:

\medskip

\( \varphi_{\neg} \ = \ \left[ \begin{array}{ll}
0 & 1 \\
1 & 0
\end{array} \right], \)

\medskip

\( \varphi_{\vee} \ = \ \left[ \begin{array}{ll}
0x_1x_2 \                 & 1x_1x_2   \\
(x_1\vee x_2) \, x_1x_2 \ & ({\ov{x_1\vee x_2}}) \ x_1x_2
\end{array} \right], \hspace{0.3in}
\varphi_{\wedge} \ = \ \left[ \begin{array}{ll}
0x_1x_2 \                   & 1x_1x_2 \\
(x_1\wedge x_2) \, x_1x_2 \ & ({\ov{x_1\wedge x_2}}) \ x_1x_2
\end{array} \right]
\),

\medskip

\noindent
where $x_1$ and $x_2$ range over $\{0,1\}$. Hence, 
${\rm domC}(\varphi_{\neg}) = {\rm imC}(\varphi_{\neg}) = \{0,1\}$, and
${\rm domC}(\varphi_{\vee}) = {\rm imC}(\varphi_{\vee})$ $=$ 
${\rm domC}(\varphi_{\wedge}) = {\rm imC}(\varphi_{\wedge}) = \{0,1\}^3$.
In order to represent the {\sc fork} function we first define 

\medskip

\( \varphi_{\rm 0f} \ = \ \left[ \begin{array}{ccc}
0  \ &\ 10 \ & \ 11   \\
00 \ &\ 01 \ & \ 1
\end{array}        \right] \);

\medskip

\noindent so, ${\rm domC}(\varphi_{\rm 0f}) = \{0, 10, 11\}$, and 
${\rm imC}(\varphi_{\rm 0f}) = \{00, 01, 1\}$.
Then {\sc fork} is simulated by 

\smallskip

$\varphi_{\rm f}$  $=$ 
$\tau_{1,2} \circ \varphi_{\vee} \circ \varphi_{\rm 0f}$.

\smallskip

\noindent Indeed, for all $x_1 \in \{0,1\}$: 
 \ $\tau_{1,2} \circ \varphi_{\vee} \circ \varphi_{\rm 0f}(0x_1)$
$ \,=\, 0x_1 x_1$.

\bigskip

For every acyclic boolean circuit $C$ we want to find a word 
$w_C \in (\Gamma_{\! V} \cup \tau)^*$ that simulates $C$; and we want the 
map $\,C \mapsto w_C$ to be polynomial-time computable (in terms of $|C|$). 

\smallskip

A standard property of {\sc dag}s is that every vertex has a {\em level} 
(or ``layer'') corresponding to its ``depth'' in the {\sc dag}.
The source vertices have level 0.
A gate or an output variable has level 1 iff only input variables of the
circuit feed into it.
A gate or an output variable has level $\ell$ iff it receives input from
levels $< \ell$ only, and at least one of its inputs comes from level
$\ell - 1$. Equivalently, the level of a vertex $v$ is the length of a 
longest path from a source to $v$. The maximum level of any sink vertex 
is called the {\it depth} of the {\sc dag}.

The following theorem is a simplification of \cite[Thm.\ 3.5]{BiCoNP}.
For a word $w \in (\Gamma_{\! V} \cup \tau)^*$ we use the {\em size}, 
denoted by $\|w\|$, as defined in Def.\ \ref{DEFsizeG21}.

\begin{thm} \label{reduction} {\bf (existence of simulation).}
 \ There is an injective function  \ $C  \mapsto w_C$ \ from the set of
acyclic boolean circuits to the set of words over $\Gamma_{\! V} \cup \tau$
with the following properties:

\smallskip

\noindent {\rm (1)} \ \ $w_C$ simulates the input-output function $f_C$ 
of $C$; 

\smallskip

\noindent {\rm (2)} \ \ the size of $w_C$ satisfies 
 \ \ $\|w_C\| \,<\, c \ |C|^6$ \ \ (for some constant $c > 0$); 

\smallskip

\noindent {\rm (3)} \ \ $w_C$ is computable from $C$ in polynomial time, 
in terms of $|C|$.
\end{thm}
{\sc Proof.} \ Item (1) refers to simulation as in Def.\ \ref{DEFsimulat}. 
In the proof we assume that $\varphi_{\neg}$, $\varphi_{\vee}$,
$\varphi_{\wedge}$, $\varphi_{\rm f}$, $\varphi_{\rm 0f}$, and $\tau_{1,2}$,
belong to $\Gamma_{\! V}$.
(If this were not the case, we could express them by fixed words over
$\Gamma_{\! V}$.)

We can assume that our acyclic circuits are {\it strictly layered}, i.e.,
a gate or an output variable at level $\ell$ only receives inputs from level
$\ell-1$. Hence, all the output variables of the circuit are at the same
level $L$, where $L$ is the depth of the circuit.
If the layering of a circuit $C$ is not strict, we can insert
{\it identity gates} to obtain strictness. An identity gate has one input
variable and one output variable, connected by a wire; the two variables 
carry the same boolean value. We will count these identity gates as gates 
in the evaluation of circuit size.
In order to make a circuit $C$ strictly layered, fewer than $|C|^2$ identity
gates need to be introduced. (Indeed, for each gate we add fewer than $|C|$ 
identity gates above it; so, in total we add fewer than $|C|^2$ identity 
gates.)

An acyclic circuit $C$ has input variables $x_1, \ldots, x_m$, output 
variables $y_1, \ldots, y_n$, and {\it internal variables} which correspond 
to the boolean values carried by internal wires (between gates or between a 
gate and an input or an output port). The internal variables at level $\ell$
(for $0 \le \ell \le L$) are denoted by $y_1^{\ell}$, $y_2^{\ell}$, 
$\ldots$, $y_{n_{\ell}}^{\ell}$.
When $\ell = L$ (output level) we have $n_L = n$ and $y_i^L = y_i$; when
$\ell = 0$ (input level) we have $n_0 = m$ and $y_i^0 = x_i$.

For every level $\ell$ (with $1 \le \ell \le L$) there is a circuit
$C_{\ell}$, called the {\it slice of} $C$ at level $\ell$: 
The input variables of the slice $C_{\ell}$ are $y_1^{\ell-1}$, $\ldots$, 
$y_{n_{\ell-1}}^{\ell-1}$; the output variables are $y_1^{\ell}$, $\ldots$, 
$y_{n_{\ell}}^{\ell}$; the gates of $C_{\ell}$ are the gates of $C$ at level 
$\ell$; we use the fact that $C$ is strictly layered.
In addition to gates, a slice $C_{\ell}$ also contains wire-swappings of its
inputs, i.e., a bit-position permutation is applied to the $n_{\ell-1}$ 
input variables. 
Every permutation of $n_{\ell-1}$ wires can be written as the composite of 
$\, \le n_{\ell-1}$ ($< |C_{\ell}|$) transpositions. 
And each $\tau_{i,j}$ has word length $\, \le 2(j-i) -1 \,$ over $\tau$ (by 
Lemma \ref{classicalFormulas}), hence it has size 
$\, \|\tau_{i,j}\| < |C_{\ell}|^2$.  Thus the input-wire permutation of a 
slice $C_{\ell}$ has size $\, < |C_{\ell}|^3$. 
Moreover, every $\tau_{i,j}$ belongs to $V$, so it does not need any 
simulation.

We use the notation \ $Y^{\ell}$  $=$  $y_1^{\ell} y_2^{\ell}$ 
$ \ \ldots \ $  $y_{n_{\ell}}^{\ell}$ \ (i.e., the concatenation of the 
variables $y_i^{\ell}$, for $i = 1, \ldots, n_{\ell}$, and 
$\ell = 0, \ldots, L$).

\medskip

\noindent {\sc Simulation of one slice}

\smallskip

In order to construct $w_C$ we first consider the special case where the
circuit $C$ consists of just one slice, hence $C$ has depth 2 (the gates of 
the slice have depth 1, the output variables have depth 2).  
Identity gates are allowed. We number the gates of $C$ from left to right.  

For $k \geq 0$, let $K$ consist of the first $k$ gates of a slice; so, $K$ 
is a one-slice circuit that has $k$ gates. When $k=0$, $K$ is empty, $w_K$
is the empty string, and its input-output function is the identity function 
($\in V$).
Inductively, let $C$ be a slice obtained from $K$ by adding one gate 
({\sc and}, {\sc or}, {\sc not}, identity, or {\sc fork}) on the right of 
$K$ (with number $k+1$). Inductively we assume that $K$ satisfies the 
Theorem and that $w_K$ has been constructed.  
Let $x_1, \ldots, x_m$ be the input variables and let $y_1, \ldots, y_n$ be 
the output variables of $K$. We now construct $w_C$ from $w_K$ and the gate
being  added. 

\medskip

{\sc Case 1:} \  Suppose the slice $C$ is obtained from $K$ by adding, on 
the right of $K$, an identity gate or a {\sc not} gate, with new input 
variable $x_{m+1}$ and new output variable $y_{n+1}$.
If a {\sc not} gate is added, the input-output function of $C$ is
$\, f_C(x_1, \ldots, x_m, x_{m+1}) = $
$(y_1, \ldots, y_n, \ov{x_{m+1}})$,
where $\,f_K(x_1, \ldots, x_m) = (y_1, \ldots, y_n)$.
The boolean function $f_C$ is to be simulated by a Thompson group element 
$\Phi_{f_C}$ such that 

\smallskip

$\Phi_{f_C}(0 \, x_1 \ldots x_m, x_{m+1}) \ = \ $
$0 \, y_1 \ldots y_n \ \ov{x_{m+1}} \ x_1 \ldots x_m x_{m+1}$.

\smallskip

\noindent We have $w_K \in (\Gamma_{\! V} \cup \tau)^*$, where 
$\Phi_{f_K} \in V$ is the simulation of $f_K$, which exists by induction. 
We find $w_C$ as follows:

\smallskip

$0\, x_1 \ldots x_m \ x_{m+1} \ $
$\stackrel{\Phi_{f_K} }{\longmapsto} \ $
$0 \ y_1 y_2 \ldots y_n \ x_1 \ldots x_m \ x_{m+1} \ $

\smallskip

$\stackrel{\tau_{2,n+m+2}}{\longmapsto} \ $
$0 \ x_{m+1} \ y_2 \ldots y_n \ x_1 \ldots x_m \ y_1$
$\stackrel{\varphi_{\rm f}}{\longmapsto} \ $
$0 \ x_{m+1} \ x_{m+1} \ y_2 \ldots y_n \ x_1 \ldots x_m \ y_1$

\smallskip

$\stackrel{\tau_{1,2}}{\longmapsto} \ $
$x_{m+1} \ 0 \ x_{m+1} \ y_2 \ldots y_n \ x_1 \ldots x_m \ y_1$

\smallskip

$\stackrel{\varphi_{\neg}}{\longmapsto} \ $
$\ov{x_{m+1}} \ 0 \ x_{m+1} \ y_2 \ldots y_n \ x_1 \ldots x_m \ y_1$
$ \ \stackrel{\tau_{3,n+m+3}}{\longmapsto}$  
$\stackrel{\tau_{1,2}}{\longmapsto} \ $
$0 \ \ov{x_{m+1}} \ y_1 y_2 \ldots y_n \ x_1 \ldots x_m x_{m+1}$ ;

\smallskip

applying 
 \ $\tau_{n+1,n+2} \ \tau_{n,n+1} \ \dots \ \tau_{3,4} \ \tau_{2,3}(.)$
 \ then yields

\smallskip

$0 \ y_1 \ldots y_n \ \ov{x_{m+1}} \ x_1 \ldots x_m x_{m+1}$.

\smallskip

\noindent So, \ $w_C$ $\,=\,$
$\tau_{n+1,n+2} \ \tau_{n,n+1} \ \dots \ \tau_{3,4} \ \tau_{2,3}$
$\tau_{1,2}$  $\tau_{3,n+m+3}$
$\varphi_{\neg}$  $\tau_{1,2}$   $\varphi_{\rm f}$   $\tau_{2,n+m+2}$
$w_K$.

\smallskip

\noindent The case where, instead of a {\sc not} gate, an identity gate
is added is similar (except that we simply omit $\varphi_{\neg}$).
By Lemma \ref{classicalFormulas}, we can express $\tau_{2,n+m+2}$ and 
$\tau_{3,n+m+3}$ over $\tau = \{\tau_{k,k+1} : k \ge 1\}$.  
Then the size of $w_C$ is 

\smallskip

$\|w_C\|$ $ \ \le \ $  $\|w_K\|$ $+$ $\|\tau_{2,n+m+2}\|$ $+$  
   $\|\tau_{3,n+m+3}\| + 4$ $+$ $\, \sum_{k=2}^{n+1} \|\tau_{k,k+1}\|$ 

\smallskip

\hspace{0.37in}  $\le \, \|w_K\| \, + \, c \, (n+m)^2 + c$,
  \ \ \ for some constant $c > 1$. 

\bigskip

{\sc Case 2:} \  Suppose our slice $C$ is obtained by adding an {\sc and} 
gate or an {\sc or} gate to $K$ on the right, with new output variable 
$y_{n+1}$ and new input variables $x_{m+1}, \, x_{m+2}$. 
We only analyze the {\sc or} case, the {\sc and} case being almost the same.
The input-output function of $C$ is

\smallskip

$f_C(x_1, \ldots, x_m, x_{m+1}, x_{m+2}) \ = \ $
$(y_1, \ldots, y_n, \ x_{m+1} \vee x_{m+2})$,

\smallskip

\noindent where $f_K(x_1, \ldots, x_m) = (y_1, \ldots, y_n)$.
The function $f_C$ is to be simulated by a Thompson group element 
$\Phi_{f_C}$ such that

\smallskip

$\Phi_{f_C}(0 \, x_1 \ldots x_m \ x_{m+1} x_{m+2}) \ = \ $
$0 \ y_1 \ldots y_n \ (x_{m+1} \vee x_{m+2}) \ x_1 \ldots x_m x_{m+1}x_{m+2}$

\smallskip

\noindent Let $w_K \in (\Gamma_2 \cup \tau)^*$ be such that 
$\Phi_{f_K} \in V$ simulates $f_K$. Then we construct $w_C$ as follows: 

\smallskip

$0\,x_1 \ldots x_m \ x_{m+1} x_{m+2} \ $
$\stackrel{\Phi_{f_K} }{\longmapsto} \ $
$0\,y_1 \ldots y_n \ x_1 \ldots x_m \, x_{m+1}\, x_{m+2}$

\smallskip

$\stackrel{\varphi_{\rm 0f}}{\longmapsto}$ 
$ \ 0 0 \ y_1 \ldots y_n \ x_1 \ldots x_m \, x_{m+1}\, x_{m+2}$

\smallskip

$\stackrel{\tau_{2,n+m+3}}{\longmapsto}$ 
$\stackrel{\tau_{3,n+m+4}}{\longmapsto} \ $
$0 \, x_{m+1} x_{m+2} \ y_2 \ldots y_n \ x_1 \ldots x_m \ 0 y_1$

\smallskip

$\stackrel{\varphi_{\vee}}{\longmapsto} $
$(x_{m+1} \vee x_{m+2})\ x_{m+1}x_{m+2} \ y_2 \ldots y_n \ $
              $x_1 \ldots x_m \ 0 y_1$

\smallskip

$\stackrel{\tau_{2,n+m+3}}{\longmapsto}$
$\stackrel{\tau_{3,n+m+4}}{\longmapsto} \ $
$(x_{m+1} \vee x_{m+2}) \ 0 \ y_1 y_2 \ldots y_n \ $
              $x_1 \ldots x_m  \, x_{m+1} \, x_{m+2}$ ;

\smallskip

applying \ $\tau_{n+1,n+2} \ \ldots \ \tau_{2,3}\, \tau_{1,2}$ \ then yields

\smallskip

$0 \ y_1 y_2 \ldots y_n \ (x_{m+1} \vee x_{m+2}) \ x_1 \ldots x_m $
$x_{m+1} x_{m+2}$.

\smallskip

\noindent Thus $C$ is simulated by the word 

$w_C = $ 
$\tau_{n+1,n+2} \ \ldots \ \tau_{2,3}\, \tau_{1,2}$
$\tau_{3,n+m+4} \, \tau_{2,n+m+3}$
$\varphi_{\vee}$ $\tau_{3,n+m+4}$  $\tau_{2,n+m+3}$
$\varphi_{\rm 0f}$  $w_K$

\noindent of size 

\smallskip 

$\|w_C\|$ $ \ \le \ $  $\|w_K\|$ $+$ $2\,\|\tau_{2,n+m+3}\|$ $+$ 
$2\, \|\tau_{2,n+m+4}\| +2$ $+$ $\, \sum_{k=1}^{n+1} \|\tau_{k,k+1}\|$

\smallskip

\hspace{0.37in} $\le \, \|w_K\| \, + \, c \, (n+m)^2 + c$, 
 \ \ \ for some constant $c > 1$.

\bigskip

{\sc Case 3:} \ Suppose our slice $C$ is obtained by adding a {\sc fork} 
gate on the right of $K$, with a new input variable $x_{m+1}$ and 
two new output variables $y_{n+1}$ and $y_{n+2}$. The input-output function 
of $C$ is

\smallskip

$f_C(x_1, \ldots, x_m, x_{m+1}) \ = \ (y_1, \ldots, y_n, x_{m+1}, x_{m+1})$,

\smallskip

\noindent where $f_K(x_1, \ldots, x_m) = (y_1, \ldots, y_n)$.
The boolean function $f_C$ is to be simulated by a Thompson group element
$\Phi_f$ such that

\smallskip

$\Phi_f(0 \, x_1 \ldots x_m x_{m+1}) \ = \ $
$ 0 \ y_1 \ldots y_n\ x_{m+1}x_{m+1} \ x_1 \ldots x_m x_{m+1}$.

\smallskip

\noindent
Let $w_K \in (\Gamma_{\! V} \cup \tau)^*$ and $\Phi_{f_K} \in V$ be the 
simulation of $f_K$, which exists by induction. Then

\smallskip

$0\, x_1 \ldots x_m \ x_{m+1} \ $
$\stackrel{\Phi_{f_K} }{\longmapsto} \ $
$0 \ y_1 y_2 \ldots y_n \ x_1 \ldots x_m \ x_{m+1} \ $

\smallskip

$\stackrel{\tau_{2,n+m+2}}{\longmapsto} \ $
$0 \ x_{m+1} \ y_2 \ldots y_n \ x_1 \ldots x_m \ y_1 \ $
$\stackrel{\varphi_{\rm f}}{\longmapsto}$ 
$\stackrel{\varphi_{\rm f}}{\longmapsto} \ $
$0 \, x_{m+1} x_{m+1}x_{m+1} \ y_2 \ldots y_n \ x_1 \ldots x_m \ y_1 \ $

\smallskip

$\stackrel{\tau_{4,n+m+4}}{\longmapsto} \ $ 
$0 \, x_{m+1} \, x_{m+1} \ y_1 \ldots y_n \ x_1 \ldots x_m \ x_{m+1}$ ;

\smallskip
 
applying \ $\tau_{n+2,n+3} \ \ldots \ \tau_{3,4}$ \ and then 
 \ $\tau_{n+1,n+2} \ \ldots \ \tau_{2,3}$ \ yields

\smallskip

$0 \ y_1 y_2 \ldots y_n \ x_{m+1} x_{m+1} \ x_1 \ldots x_m x_{m+1}$.

\smallskip

\noindent This simulates $f_C$ by a word 

$w_C = $  $\tau_{n+1,n+2} \ \ldots \ \tau_{2,3}$  
$\tau_{n+2,n+3} \ \ldots \ \tau_{3,4}$  $\tau_{4,n+m+4}$  
$\varphi_{\rm f} \ \varphi_{\rm f} \ \tau_{2,n+m+2}$ $w_K$

\noindent of size

$\|w_C\| \ \le \ \|w_K\| + 2 + \|\tau_{2,n+m+2}\| + \|\tau_{4,n+m+4}\|$ $+$
$ \ \sum_{k=3}^{n+2} \|\tau_{k,k+1}\|$  $ \ + \ $
$\sum_{k=2}^{n+1} \|\tau_{k,k+1}\|$

\hspace{0.37in} $\le \, $ $\|w_K\| + c \, (m+n)^2 + c$, 
 \ \ \ for some constant $c > 1$.

\bigskip

\noindent In all three cases the slice $C$ is simulated by a word $w_C \in$
$(\Gamma_{\! V} \cup \tau)^*$ of size 

\smallskip

$\|w_C\| \le \|w_K\| + c \, (m+n)^2 + c$. 

\smallskip

\noindent Let $S$ now be any slice, and let $n_i$ be the number of interior 
vertices of $S$ (i.e., the vertices labeled by gates). Then if $w_S$ is 
constructed by adding $n_i$ ($< |S|)$ gates to slices (starting with $K$ 
being the empty slice, and ending with $K$ being the desired slice $S$), 
the size of $w_S$ is

\smallskip
 
$\|w_S\| \ \le \ n_i \, (c \, (m+n)^2 + c) \ \le \ c_0 \, |S|^3$, 

\smallskip

\noindent  for some constant $c_0 > 1$ (that does not depend on $S$).

Moreover, as we saw when we introduced the notion of slice, 
in all three cases a bit-position permutation of the input wires of the 
slice $S$ is attached at the beginning of $w_S$.
This permutation belongs to $V$ and has size $< |S|^3$.

The above construction of each word $w_S$ from $S$ is a polynomial-time 
algorithm (in terms of $|S|$).

\bigskip

\noindent {\sc Simulation of a multi-slice circuit} 

\smallskip

Assume that $C$ is a circuit of depth $L > 2$; the depth is the number of
slices.  In order to define $w_C$ we use the fact that we have already 
defined the word $w_{C_{\ell}}$ that simulates the slice $C_{\ell}$ of $C$ 
(for every $\ell$, $1 \le \ell \le L$).
Each word $w_{C_{\ell}}$ has all the properties claimed in Theorem
\ref{reduction}; in particular, $w_{C_{\ell}}$ represents the function

\smallskip

$\Phi_{C_{\ell}} : \ \ 0 \ Y^{\ell-1} \ \longmapsto \ $
$0 \ Y^{\ell} \ Y^{\ell-1}$.

\smallskip

\noindent
Hence, since $\Phi_{C_{\ell}}$ is a right ideal isomorphism, we also have

\smallskip

$0 \ Y^{\ell-1} \ Y^{\ell-2} \ \ldots \ Y^1 \ x_1 \ldots x_m$
$ \ \ \stackrel{\Phi_{C_{\ell}}}{\longmapsto} \ \ $
$0 \ Y^{\ell} \ Y^{\ell-1} \ Y^{\ell-2} \ \ldots \ Y^1 \ x_1 \ldots x_m$. 

\smallskip

\noindent Therefore,
 \ $w_{C_L} \, w_{C_{L-1}} \ \ldots \ w_{C_{\ell}} \ \ldots \ w_{C_1}$ 
 \ represents the function 

\smallskip

$\Phi_{C_L C_{L-1} \ldots C_1} : $
$ \ \ 0 \, x_1 \ldots x_m \ \ \longmapsto \ \ $
$0 \ y_1 \ldots y_n \, Y^{L-1} \, \ldots \, Y^{\ell} \, \ldots \, $
$Y^2 \, Y^1 \, x_1 \ldots x_m \ \ (=_{\rm def} \ Z)$, 

\smallskip

\noindent where $\,y_1 \ldots y_n = Y^L \,$ is the output of $C$, and 
$x_1 \ldots x_m$ is the input of $C$.

\medskip

\noindent The length of the word $Z$ ($\in \{0,1\}^*$) is
 \ $|Z| \,\le\, 1 + |C|$. Indeed, the total number of variables in the 
circuit (i.e., $n_L + \,\ldots\, + n_1 + m$) is equal to the total number 
of wires (i.e., $|C|$); the ``$+1$'' comes from the leading bit $0$.

\smallskip

Let $\,\sigma_{i,j} = \tau_{j-1,j} \, \tau_{j-2,j-1} \ \ldots \ $
$\tau_{i+1,i+2} \, \tau_{i,i+1}(.)$ \ (for $1 \le i < j$).
Then $ \ \pi_1 \ = \ (\sigma_{1,|Z|})^n \ $ transforms the word $Z$ into

\smallskip

$0 \ Y^{L-1} \ \ldots \ Y^{\ell} \ \ldots \ Y^2 \ Y^1 \ x_1 \ldots x_m$ \
$y_1 \ldots y_n$.

\smallskip

\noindent Next (and this is a fundamental and crucial idea from 
{\em reversible computing}, see e.g., \cite{Ben89, Ben73, FredToff}), to 
the latter string we apply

\medskip

 \ \ \  \ \ \   
$(w_{C_{L-1}} \ \ldots \ w_{C_{\ell}} \ \ldots \ w_{C_2}\, w_{C_1})^{-1}$  

\medskip

\noindent in order to clear away intermediate outputs of all the internal 
slices. This yields  

\smallskip

 \ \ \  \ \ \ $0 \ x_1 \ldots x_m \ y_1 \ldots y_n$.

\smallskip

\noindent Finally, applying the permutation 
$\, \pi_2 \ = \ (\sigma_{1,n+m})^m \,$ produces the desired final output

\smallskip

 \ \ \  \ \ \ $0 \ y_1 \ldots y_n \ x_1 \ldots x_m$.

\smallskip

\noindent Therefore we define $w_C \in (\Gamma_{\! V} \cup \tau)^*$) by

\medskip

 \ \ \  \ \ \   
$w_C \ = \ \pi_2 \ (w_{C_{L-1}} \ \ldots \ w_{C_1})^{-1} \  \pi_1 \ $
$ w_{C_L} \ w_{C_{L-1}} \ \ldots \ w_{C_1}\,$.

\medskip

\noindent The word length of $\pi_1$ over $\tau$ is less than $\,n \ |Z|$. 
Since all subscripts in $\sigma_{1,|Z|}$ are $\le |Z|$, the size of $\pi_1$
is $\, \|\pi_1\| \,<\, |Z| \ n \ |Z| \, \le\,  (|C| + 1)^3$.
Since $m+n \le |C|$, the size of $\pi_2$ is also less than $\,(|C| + 1)^3$.

\smallskip

\noindent For the size of $w_C$ we have 

\smallskip

$\|w_C\| \le \|\pi_2\| + \|\pi_1\| + \|w_{C_L}\|$ 
        $ + \ 2 \, \sum_{\ell=1}^{L-1} \|w_{C_{\ell}}\|$.

\smallskip

\noindent We saw that $ \ \|w_{C_{\ell}}\| \le c_0 \, |C_{\ell}|^3 \ $ 
(for $1 \le \ell \le L$); and 
$\, \sum_{\ell = 1}^L |C_{\ell}| = |C| \,$ implies
$\, \sum_{\ell=1}^L |C_{\ell}|^3 \le |C|^3$.
Thus $\, \|w_C\| \le c \cdot |C|^3$, for some positive constant $c$.

\smallskip

\noindent Since $|C|$ was at most squared in order to obtain strict layering, 
the above bound becomes

\medskip

$\|w_C\| \, \le \,  c \ |C|^6$ ,

\medskip

\noindent in terms ot the original (not necessarily strictly layered) 
circuit $C$. 

The word $w_C$ can be written down in linear time, based on the words
$w_{C_{\ell}}$ ($1\le \ell \le L$), and we saw that each $w_{C_{\ell}}$
can be computed in polynomial time from $C_{\ell}$.
 \ \ \ $\Box$

\subsection{Reduction to a generalized word problem of {\boldmath $V$} \\
  (over an infinite generating set)}

We first extend the classical concepts of stabilizer and fixator to the 
case of partial injections.

\begin{defn} \label{DEFpStabpFix}
A function $g$ {\em partially stabilizes} a set $S \subseteq \{0,1\}^*$ 
 \ iff \ $g(S) \cup g^{-1}(S) \subseteq S$.
For a subgroup $G \subseteq V$, the {\em partial stabilizer} of $S$ (in $G$) 
is

\medskip

 \ \ \ \ \  ${\rm pStab}_G(S) \ = \ $
  $\{g \in G : \ g(S) \,\cup\, g^{-1}(S) \,\subseteq\, S \}$.

\medskip

\noindent A function $g$ {\em partially fixes} a set $S$ \ iff 
 \ $g(x) = x\,$ for every $\, x \in $
$S \,\cap\, {\rm Dom}(g) \,\cap\, {\rm Im}(g)$.
This is also called partial pointwise stabilization.
For a subgroup $G \subseteq V$, the {\em partial fixator} of $S$ (in $G$) is

\medskip

 \ \ \ \ \  ${\rm pFix}_G(S) \ = \ \{g \in G : \ $
$(\forall x \in S \,\cap\, {\rm Dom}(g) \,\cap\, {\rm Im}(g))$
$[\,g(x) = x\,]\,\}$.
\end{defn}
We will only use partial stabilizers and fixators for sets $S$ that are
right ideals; then ${\rm pStab}_G(S)$ and ${\rm pFix}_G(S)$ are
groups \cite[Lemma 4.1]{BiCoNP}. 
When $S = P \{0,1\}^*$ is a right ideal, where $P$ is a prefix 
code, we will abbreviate ${\rm pFix}_G(P \, \{0,1\}^*)$ and
${\rm pStab}_G(P \, \{0,1\}^*)$ by ${\rm pFix}_G(P)$, respectively 
${\rm pStab}_G(P)$.  In particular, we abbreviate 
${\rm pFix}_V(0 \, \{0,1\}^*)$ to ${\rm pFix}_V(0)$.

\begin{lem} \label{pFixINpStab}
 \ We have: \ \ ${\rm pFix}_V(0) \ \subset \ $
${\rm pStab}_V(0 \, \{0,1\}^*) \ \cap \ {\rm pStab}_V(1 \, \{0,1\}^*)$.
\end{lem}
{\sc Proof.} Obviously, ${\rm pFix}_V(0) \subset$
${\rm pStab}_V(0 \, \{0,1\}^*)$. Moreover, if we had $g(1x) = 0y$ for any
$g \in {\rm pFix}_V(0)$ and $x, y \in \{0,1\}^*$, then
$0y = g^{-1}(0y) = g^{-1}g(1x) = 1x$; the first equality holds since 
$g^{-1} \in {\rm pFix}_V(0)$.  But $0y = 1x$ is false since a string does 
not start with both 0 and 1. 
 \ \ \ $\Box$

\bigskip

\noindent The following is little more than a reformulation of the 
definition of simulation and Lemma \ref{simulatVSfix}.

\begin{lem} \label{simulatVSfixators} \
Let $f$ and $g$ be any boolean functions such that $f$ and $g$ have the same
number of input variables, and $f$ and $g$ have the same number of output
variables. Suppose $f$ and $g$ are simulated by $\Phi_f$, respectively
$\Phi_g$ \ ($\Phi_f, \Phi_g \in V$). Then,

\medskip

\hspace{1.in}  $f = g$ \ \ \ iff
 \ \ \ $\Phi_f^{-1} \, \Phi_g \,\in\, {\rm pFix}_V(0)$.
\end{lem}
{\sc Proof.} Let $\{0,1\}^m$ be the common domain of $f$ and $g$. Then by
Lemma \ref{simulatVSfix}, $f = g$ iff for all $x \in \{0,1\}^m$: 
$\, \Phi_f(0x) = \Phi_g(0x)$. Then for all $x \in \{0,1\}^m$: 
$\, 0x = \Phi_f^{-1} \, \Phi_g(0x) = \Phi_g^{-1} \, \Phi_f(0x)\,$
(and $\Phi_g^{-1} \, \Phi_f = (\Phi_f^{-1} \, \Phi_g)^{-1}$). Hence,
$f = g\,$ iff $\,\Phi_f^{-1} \, \Phi_g \in {\rm pFix}_V(0)$.
 \ \ \   $\Box$

\bigskip

\noindent Theorem \ref{reduction} and Lemma \ref{simulatVSfixators} give a
polynomial-time one-one reduction from the circuit equivalence problem 
to the generalized word problem of
$\, {\rm pFix}_V(0)$ in $V$, where the elements of $V$ written 
over $\Gamma_{\! V} \cup \tau$.
Since the circuit equivalence problem is {\sf coNP}-complete, it follows 
that this generalized word problem is {\sf coNP}-hard.  Hence we have:

\begin{cor} \label{GenWP} {\bf (coNP-hard generalized word problem).} 
The generalized word problem of $\,{\rm pFix}_V(0)$ in $V$ over 
$\Gamma_{\! V} \cup \tau$ is {\sf coNP}-hard.
 \ \ \  \ \ \  \ \ \  $\Box$
\end{cor}

\subsection{Reduction to the word problem of {\boldmath $V$}} 

We will give a linear-time $2$-ary conjunctive reduction from the 
generalized word problem of ${\rm pFix}_V(0)$ to the word problem of $V$ 
over the infinite generating set $\Gamma_{\! V} \cup \tau$.
This reduction is based on a ``commutation test'', that was studied in 
greater generality in \cite[Section 5]{BiCoNP}; here we just use $V$,
based on an alphabet of size 2, which makes everything simpler.

We first need a few lemmas. Recall the notation $u \para_{\rm pref} v\,$ 
($u$ and $v$ are prefix-comparable) and its negation $\nparallel_{\rm pref}$. 
For $x \in \{0,1\}^*$ and $L \subseteq \{0,1\}^*$, we define
$\,x^{-1} L \,=\, \{v \in \{0,1\}^*: xv \in L\}$.

\begin{lem} \label{LEMg0xVS0x} 
 \ If $g \not\in {\rm pFix}_V(0)$ but $g \in {\rm pStab}_V(0)$, then there 
exists $0x \in {\rm domC}(g)$ such that 

\medskip

\hspace{0.3in}  $0x \nparallel_{\rm pref} g(0x)$.
\end{lem}
Hence, $\, 0x u \nparallel_{\rm pref} g(0x u) \,$ \big($= g(0x) \, u$\big), 
for all $u \in \{0,1\}^*$.

\medskip

\noindent {\sc Proof.} Lemma \ref{LEMg0xVS0x} is a special case of 
\cite[Lemma 9.6]{BiCoNP}, with a simpler proof.
(Note that in \cite{BiCoNP} the notation $\le_{\rm pref}$ for the prefix 
order was reversed; here, ``$p \le_{\rm pref} w$'' always means $p$ is a 
prefix of $w$.)

\medskip

We prove the contrapositive, i.e., if for all $0x \in {\rm domC}(g)$ we 
have $0x \para_{\rm pref} g(0x)$, then $g \in {\rm pFix}_V(0)$.

\smallskip

\noindent {\sc Case} 1: \ $0x <_{\rm pref} g(0x)$.

\smallskip

Then $g(0x) = 0x \, v$, for some $v \in \{0,1\}^+$, so 
$v \in (0x)^{-1} {\rm imC}(g)$. 
Now, $(0x)^{-1}({\rm imC}(g))$ is a maximal finite prefix code (by 
\cite[Lemma 9.4]{BiCoNP}), which contains the non-empty string $v$. Hence 
$(0x)^{-1} {\rm imC}(g)$ contains at least one other non-empty string (by 
\cite[Lemma 9.5]{BiCoNP}), i.e., ${\rm imC}(g)$ contains $0x w$ ($\ne 0xv$), 
for some $w \in \{0,1\}^+$. 
Hence (since $g^{-1}$ stabilizes $0 \{0,1\}^*$), there exists 
$0x' \in {\rm domC}(g)$ such that $0x' \ne 0x$, and
$g(0x') \in {\rm imC}(g)$ and $g(0x') = 0x \, w >_{\rm pref} 0x$. 
Since ${\rm imC}(g)$ is a prefix code, $g(0x) \nparallel_{\rm pref} g(0x')$. 

By the (contrapositive) assumption, $0x' \para_{\rm pref} g(0x')$.
Hence there are two possibilities:

\smallskip

\noindent (1) \ $0x' \ge_{\rm pref} g(0x')$: \ Then 
$0x' \ge_{\rm pref} g(0x') >_{\rm pref} 0x$. So $0x' >_{\rm pref} 0x$,
which contradicts the fact that ${\rm domC}(g)$ is a prefix code.

\smallskip

\noindent (2) \ $0x' <_{\rm pref} g(0x')$: \ Then
$0x' <_{\rm pref} g(0x') = 0x' \, z$, for some $z \in \{0,1\}^+$; and we saw
that also $g(0x') = 0x \, w$. This implies that $0x \|_{\rm pref} 0x'$. 
Again, this contradicts that ${\rm domC}(g)$ is a prefix code.

Thus, case 1 is impossible.

\smallskip

\noindent {\sc Case} 2: \ $0x >_{\rm pref} g(0x)$.

\smallskip

Then $0x = g(0x) \, u$, for some $u \in \{0,1\}^+$, hence 
$u \in (g(0x))^{-1}{\rm domC}(g)$. Now $(g(0x))^{-1}{\rm domC}(g)$ is a 
finite maximal prefix code, containing the non-empty string $u$, hence it 
contains some other non-empty string. So there exists $0x'$ ($\ne 0x$) with 
$0x' \in {\rm domC}(g) \,\cap\, g(0x) \, \{0,1\}^+$.

By the (contrapositive) assumption, $0x' \para_{\rm pref} g(0x')$.
Again, we have two possibilities:

\smallskip

\noindent (1) \ $0x' \le_{\rm pref} g(0x')$:  \ Then
$g(0x') \ge_{\rm pref} 0x'$, and  $0x' >_{\rm pref} g(0x)$ (since 
$0x' \in g(0x) \, \{0,1\}^+$).  Thus, $g(0x') >_{\rm pref} g(0x)$, 
which contradicts the fact that ${\rm imC}(g)$ is a prefix code.

\smallskip

\noindent (2) \ $0x' >_{\rm pref} g(0x')$: \ Then
$0x' = g(0x') \, z$, for some $z \in \{0,1\}^+$; and
$0x' = g(0x) \, w$, for some $w \in \{0,1\}^+$ (since 
$0x' \in g(0x) \, \{0,1\}^+$).
Thus, $0x' = g(0x') \, z = g(0x) \, w$, which implies 
$g(0x') \para_{\rm pref} g(0x)$. Again, this contradicts the fact that 
${\rm imC}(g)$ is a prefix code.

We conclude that case 2 is impossible.

\smallskip

Now, having ruled out cases 1 and 2, the only remaining possibility is that 
$0x = g(0x)$, for all $0x \in {\rm domC}(g)$. This means that 
$g \in {\rm pFix}_V(0)$.  
 \ \ \ $\Box$

\begin{lem} \label{LEMf0xVS0y} 
 \ For every $\, 0x, 0y \in 0\,\{0,1\}^*$ such that 
$0x \nparallel_{\rm pref} 0y$, there exists $f_0 \in {\rm pFix}_V(1)$ and 
$u \in \{0,1\}^*$ such that 

\medskip

\hspace{1.1in} $f_0(0xu) = 0xu$ \ \ and \ \ $f_0(0yu) \ne 0yu$.
\end{lem}
{\sc Proof.} This Lemma is a simplification of 
\cite[Prop.\ 9.14(1)]{BiCoNP}, and we adapt that proof.  

Let $0x, 0y \in 0\,\{0,1\}^*$ be two prefix-incomparable strings, and let 
$a, b \in \{0,1\}$ with $a \ne b$.  Then $0x, 0ya, 0yb$ are 
prefix-incomparable two-by-two (as is easy to check).
We now use \cite[Lemma 9.7]{BiCoNP} to construct a finite maximal prefix 
code $\,Q \,\cup\, \{0x, 0ya, 0yb, 1\}$, with $Q \subset 0\,\{0,1\}^*$.
We define $f_0 \in V$ by 

\smallskip

 \ \ \ $f_0(0ya) = 0yb, \ \ f_0(0yb) = 0ya, \ \ f_0(0x) = 0x$, \ \ and

\smallskip

 \ \ \ $f_0$ is the identity on $Q \,\cup\, \{1\}$.

\smallskip

\noindent So, $Q \,\cup\, \{0x, 0ya, 0yb, 1\}$ is the domain code and image 
code of $f_0$. Then $f_0 \in {\rm pFix}_V(1)$,
$f_0(0ya) \ne 0ya$, and $f_0(0xa) = 0xa$ (since $f_0(0x) = 0x$). So here, 
$a$ plays the role of $u$. 
 \ \ \ $\Box$

\begin{lem} \label{CommTest01} {\bf (commutation test).} 
 \ For all $g \in V$ we have: 

\smallskip

\hspace{0.6in} $g \in {\rm pFix}_V(0)$ \ \ \ iff 
 \ \ \ $\big(\forall f \in {\rm pFix}_V(1)\big) \,[\,fg = gf\,]$.
\end{lem}
In words: \ An element $g \in V$ belongs to the subgroup ${\rm pFix}_V(0)$ 
iff $g$ commutes with all the elements of the subgroup ${\rm pFix}_V(1)$.

\medskip

\noindent
{\sc Proof.} ${\boldmath [\Leftarrow]}$ \ Suppose $fg = gf$, for all 
$f \in {\rm pFix}_V(1)$, and hence also $\, g^{-1} f = f g^{-1}$.

\smallskip

\noindent (1) We first prove that $g \in {\rm pStab}_V(0)$. 

If $g(0x) = 1y$ for some $x, y \in \{0,1\}^*$,  then
$fg(0x) = f(1y) = 1y$ for all $f \in {\rm pFix}_V(1)$.
And $1y = fg(0x) = gf(0x)$.
Hence, $g(0x) = 1y = g(f(0x))$, hence by injectiveness, $0x = f(0x)$ 
for all $f \in {\rm pFix}_V(1)$. 
So, $f(0x0) = 0x0$ and $f(0x1) = 0x1$, and $0x0 \nparallel_{\rm pref} 0x1$. 
Hence by Lemma \ref{LEMf0xVS0y}, there exists $f_o \in {\rm pFix}_V(1)$ 
such that $f_o(0x0 u) = 0x0 u$, and $f_o(0x1 u) \ne 0x1 u$ (for some 
$u \in \{0,1\}^*$).   The latter inequality contradicts the fact that 
$f(0x) = 0x$ for all $f \in {\rm pFix}_V(1)$.

In a similar way one obtains a contradiction if $g^{-1}(0x) = 1y$ for some
$x, y \in \{0,1\}^*$.
 
\smallskip

\noindent (2) We prove next that $g \in {\rm pFix}_V(0)$.

Suppose $fg = gf$ for all $f \in {\rm pFix}_V(1)$; we saw that then 
$g \in {\rm pStab}_V(0)$.
If, by contradiction, $g \not\in {\rm pFix}_V(0)$, then by Lemma 
\ref{LEMg0xVS0x}, there exists $0x \in {\rm domC}(g)$ such that 
$\,0x \nparallel_{\rm pref} g(0x) = 0y$. 

Then, $fg(0x) = f(0y) = gf(0x)$. By Lemma \ref{LEMf0xVS0y} there exists
$f_o \in {\rm pFix}_V(1)$ such that $f_o(0x u) = 0x u$, and 
$f_o(0y u) \ne 0y u$ \ (for some $u \in \{0,1\}^*$).
Then $f_o(0y u) = f_o g(0x u) = g f_o(0xu) = g(0x u) = 0y u$. So, 
$f_o(0y u) = 0y u$, which contradicts $f_o(0y u) \ne 0y u$.

\smallskip

{\boldmath $[\Rightarrow]$} \ Let $g \in {\rm pFix}_V(0)$ and 
$f \in {\rm pFix}_V(1))$. Then ${\rm domC}(f) = \{1\} \,\cup\, 0 P$, and
${\rm domC}(g) = \{0\} \,\cup\, 1 Q$, where $P, Q \subset \{0,1\}^*\,$ are 
finite maximal prefix codes. So, $0 P \,\cup\, 1 Q$ is a finite maximal 
prefix code.

Then for every $0x \in 0 P$: $\, fg(0x) = f(0x)$, since 
$g \in {\rm pFix}_V(0)$; and
$\, gf(0x) = f(0x)$, since $f(0x) \in 0\{0,1\}^*$ and $g \in {\rm pFix}_V(0)$.
So, $fg(0x) = gf(0x)$.

Similarly, for all $1x \in 1 Q$: $\, gf(1x) = g(1x)$, since 
$f \in {\rm pFix}_V(1))$; and 
$\, fg(1x) = g(1x)$, since $g(1x) \in 1\{0,1\}^*$ and $f \in {\rm pFix}_V(1))$.
So, $fg(1x) = gf(1x)$.

Hence, $fg = gf$ on the finite maximal prefix code 
$\,0 P \,\cup\, 1 Q$. Hence $fg = gf$ in $V$.
 \ \ \ $\Box$

\begin{lem} \label{Fix_iso_V_endm}
 \ The subgroups ${\rm pFix}_V(1)$ and ${\rm pFix}_V(0)$ are isomorphic 
to $V$.
\end{lem}
{\bf Proof.} Every element of $V$ has a table of the form 
$\,\{(x_i, y_i) : 1 \le i \le n\}$, where $\{x_1, \ldots, x_n\}$ and 
$\{y_1, \ldots, y_n\}$ are finite maximal prefix codes over $\{0,1\}$.
An isomorphism $\,V \to {\rm pFix}_V(1)\,$ is given by 
\[ g \ = \ \left[ \begin{array}{ccc}
x_1 & \ldots & x_n \\
y_1 & \ldots & y_n
\end{array}        \right] \ \ \longmapsto
 \ \ \theta(g) \ = \ \left[ \begin{array}{c ccc}
1 & 0x_1 & \ldots & 0x_n \\
1 & 0y_1 & \ldots & 0y_n
\end{array}        \right]. 
\]
The map $\theta$ is obviously a bijection from $V$ onto ${\rm pFix}_V(1)$, 
and it is easy to check that it is a homomorphism.
 \ \ \ $\Box$

\bigskip

\noindent {\bf {\sf coNP}-hardness of the word problem of {\boldmath $V$} 
over {\boldmath $\Gamma_{\! V} \cup \tau$}:}

\smallskip

The commutation test of Lemma \ref{CommTest01} reduces the
generalized word problem of ${\rm pFix}_V(0)$ in $V$ (over
$\Gamma_{\! V} \cup \tau$) to an infinite set of word problems of $V$, namely
$\, \{fg = gf : f \in {\rm pFix}_V(1)\}$.

However, ${\rm pFix}_V(1)$ is 2-generated; this follows from Lemma 
\ref{Fix_iso_V_endm} and the fact that $V$ is 2-generated 
\cite{Th, DMason, BleakQuick}.
Obviously, $g$ commutes with all of ${\rm pFix}_V(1)$ iff $g$ commutes 
with the two generators of ${\rm pFix}_V(1)$. 
This reduces the generalized word problem of ${\rm pFix}_V(0)$ in $V$ (over
$\Gamma_{\! V} \cup \tau$) to the conjunction of two instances of the word 
problem of $V$ over $\Gamma_{\! V} \cup \tau$. Hence, the word problem of $V$ 
over $\Gamma_{\! V} \cup \tau$ is {\sf coNP}-hard with respect to 2-ary 
conjunctive polynomial-time reduction.

\begin{thm} \label{WPofVcoNPhard} {\bf (coNP-complete word problem).}
 \ The word problem of $V$ over the generating set $\Gamma_{\! V} \cup \tau$
is {\sf coNP}-complete.
\end{thm}
{\sc Proof.}  By Lemma \ref{G21incoNP}, this word problem belongs to
{\sf coNP}. 
By the reasoning in the above few lines, the word problem is 
{\sf coNP}-hard.  
 \ \ \  $\Box$

\bigskip

\subsection{Alternative proof of {\sf coNP}-completeness 
of the word problem of {\boldmath $V$ over $\Gamma_{\! V} \cup \tau$} } 

The above proof of {\sf coNP}-completeness of the word problem of $V$ 
over $\Gamma_{\! V} \cup \tau$ was derived from a similar proof for 
$G_{3,1}$ \cite{BiCoNP} (in 2003). 
Since then, Stephen Jordan \cite{Jordan} (in 2013) proved that the 
equivalence problem for bijective circuits built from copies of the Fredkin
gate is {\sf coNP}-complete. A {\em bijective circuit} is an acyclic circuit 
in which every gate has a permutation of $\{0,1\}^j$ as its input-output 
function (for some $j > 0$, depending on the gate). The input-output 
function of such a circuit is a permutation of $\{0,1\}^n$ for some $n > 0$ 
(see e.g.\ \cite{ShPMH}).  The {\em Fredkin gate}, on an input 
$x_1 x_2 x_3 \in \{0,1\}^3$, is defined by 

\smallskip

 \ \ \  \ \ \  ${\small \sf F}(0 \, x_2 x_3) = 0 \, x_2 x_3\,$, 

 \ \ \  \ \ \  ${\small \sf F}(1 \, x_2 x_3) = 1 \, x_3 x_2\,$;

\smallskip

\noindent see e.g.\ \cite{FredToff}. This gate is also called the
``controlled transposition'' (of $x_2$ and $x_3$).
Clearly, ${\small \sf F}$ is the table of an element of $V$; moreover, 
with $\{{\small \sf F}\} \cup \tau$ we can compute 
${\small \sf F}(x_i x_j x_k)$ for any three different variables in an 
input $x_1 \ldots x_n$ with $i, j, k \in \{1,\ldots,n\}$. Hence 
Jordan's result can be recast as follows:

\begin{thm} \label{THMVversJordan} {\bf (Thompson group form of Jordan's
theorem).}
 \ The subgroup of $V$ generated by $\,\{{\small \sf F}\} \cup \tau$ has a
{\sf coNP}-complete word problem, with respect to many-one polynomial-time 
reduction. \ \ \  \ \ \  $\Box$
\end{thm}
See \cite{BiFact} (and \cite{ShPMH}) for further connections between 
bijective (``reversible'') circuits and the Thompson groups.

\medskip

Theorem \ref{THMVversJordan} immediately implies Theorem 
\ref{WPofVcoNPhard}, as the word problem of the subgroup 
$\langle \{{\small \sf F}\} \cup \tau \rangle_{_V}$ reduces to the word 
problem of $V$ over $\Gamma_{\! V} \cup \tau$ by the inclusion map.  
(Here we assume that ${\small \sf F} \in \Gamma_{\! V}$; if that is not the 
case we can represent ${\small \sf F}$ by a fixed word over $\Gamma_{\! V}$ 
for the reduction; see Lemma \ref{INFINComplexSub}(2).)  

\smallskip

An advantage of our method of subsections 4.2 - 4.4 is that it is 
direct, whereas Jordan's theorem is based on Barrington's theorem 
\cite{Barrington}, which is itself a deep result.   However, using Jordan's 
theorem has the advantage that it yields the following:
The word problem of $V$ over $\Gamma_{\! V} \cup \tau$ is 
{\sf coNP}-complete with respect to polynomial-time {\em many-one reduction}. 
The earlier proof only yields polynomial-time binary conjunctive reduction.

\subsection{The shift, and the word problem of {\boldmath $nV$} }

For all $\tau_{j, j+1} \in \tau \subset V$ with $j \ge 1$, we define 
 \ $\tau_{j,j+1} \times {\bf 1}$: $\{0,1\}^* \times \{0,1\}^*$  
$\,\longrightarrow\,$  $\{0,1\}^* \times \{0,1\}^*$ \ by
 
\medskip

 \ \ \   \ \ \ $\tau_{j,j+1} \times {\bf 1}:$
$ \ \ (x,y) \ \longmapsto \ (\tau_{j,j+1}(x), \ y)$.

\medskip

\noindent So, ${\rm domC}(\tau_{j,j+1} \times {\bf 1}) = $
$\{0,1\}^{j+1} \times \{\varepsilon\}$. 

\medskip

\noindent The {\em shift} $\,\sigma \in 2V$ is defined by 
$\, {\rm domC}(\sigma) = \{\varepsilon\} \times \{0,1\}$, 
$\, {\rm imC}(\sigma) = \{0,1\} \times \{\varepsilon\}$, and

\smallskip

 \ \ \   \ \ \  $\sigma(\varepsilon, a) = (a, \varepsilon)$, 

\smallskip

\noindent for all $a \in \{0,1\}$.
Hence, $\sigma(x, \, ay) = (ax, \, y)$, for all $a \in \{0,1\}$, and
$\, (x, y) \in \{0,1\}^* \times \{0,1\}^*$.

\begin{lem} \label{tauVSsigma}
 \ For all $j \ge 1$: \ \ \ $\tau_{j,j+1} \times {\bf 1}(.) \ = \ $
$\sigma^{j-1} \circ (\tau_{1,2} \times {\bf 1}) \circ \sigma^{-j+1}(.)$ .
\end{lem}
{\sc Proof.} For any $(x,y) \in \{0,1\}^* \times \{0,1\}^*$, where 
$\, x = u \, x_j x_{j+1} \, v$ with $|u| = j-1 \ge 0$, and 
$v \in \{0,1\}^*$, we have:

\smallskip

$(u \, x_j x_{j+1} v, \, y) \ \stackrel{\sigma^{-j+1}}{\longmapsto} \ $
$(x_j x_{j+1} \, v, \   u^{\rm rev} \, y) \ $
$\stackrel{\tau_{1,2} \times {\bf 1}}{\longmapsto} \ $
$(x_{j+1} x_j \, v, \   u^{\rm rev} \, y) \ $
$\stackrel{\sigma^{j-1}}{\longmapsto} \ $
$(u \, x_{j+1} x_j \, v, \   y)$.

\smallskip

\noindent Here, $u^{\rm rev}$ denotes the reverse of $u$.
 \ \ \ $\Box$
 
\bigskip

\noindent {\bf Proof of Theorem \ref{2g21wpCoNP}:}

\smallskip

By Lemma \ref{nVinCoNP} the word problem of $nV$ belongs to {\sf coNP}.

\smallskip

By Theorem \ref{WPofVcoNPhard} the word problem of $V$ over
$\Gamma_{\! V} \cup \tau$ is {\sf coNP}-hard.
By Lemma \ref{tauVSsigma}, the word problem of $V$ over 
$\Gamma_{\! V} \cup \tau$, reduces to the word problem of $2V$ over a finite
generating set; this reduction is the one-one reduction that replaces every 
generator $\gamma \in \Gamma_{\! V}$ by $\gamma \times {\bf 1}$, and replaces 
$\tau_{j,j+1}$ by 
$\,\sigma^{j-1} \circ (\tau_{1,2} \times {\bf 1}) \circ \sigma^{-j+1}$, 
as in Lemma \ref{tauVSsigma}.  We can include the set
$\{\gamma \times {\bf 1} : \gamma \in \Gamma_{\! V}\} \cup \{\sigma\}$ into 
the finite generating set of $2V$, or we can express all the elements of this
set by a finite set of strings over some other finite generating set of 
$2V$.  
Thus the word problem of $2V$ over a finite generating set is 
{\sf coNP}-hard.

\smallskip

To show that word problem of $nV$ (for $n > 2$) over a finite generating 
set is {\sf coNP}-hard, we use the fact that $2V$ is a finitely generated 
subgroup of $nV$, and apply Lemma \ref{ComplexSub}(2). 
 \ \ \ $\Box$

\bigskip

\noindent {\bf Remark on the distortion of {\boldmath $V$} in 
{\boldmath $2V$}:} Burillo and Cleary \cite{BurCleary} show that $V$ is 
exponentially distorted in $2V$ (when both $V$ and $2V$ are over finite 
generating sets). 
In Lemma \ref{tauVSsigma} we proved that $\tau_{j-1,j}$ has linear word
length in $2V$ (as a function of $j$); but $\tau_{j-1,j}$ has exponential 
word length in $V$ over any finite generating set.  This, again, shows that 
the distortion of $V$ in $2V$ is at least exponential. 

Indeed, for all $j \ge 2$, $\,\tau_{j-1,j}$ has a table 
 \ $u x_{j-1} x_j \in \{0,1\}^j \longmapsto u x_j x_{j-1} \in \{0,1\}^j$ 
 \ (see the beginning of subsection 4.4).
It follows from Lemma \ref{extStep} that the table of $\tau_{j-1,j}$
is maximally extended. So, $\tau_{j-1,j}$ has table-size $2^j$. Therefore, 
by \cite[Thm.\ 3.8]{BiThomps}, the word length $|\tau_{j-1,j}|_{_V}$ of 
$\tau_{j-1,j}$ in $V$ (over any finite generating set) satisfies  
 \ $\alpha\, 2^j \,\le\, |\tau_{j-1,j}|_{_V} \,\le\, \beta\, j \, 2^j$
 \ (for some constants $\alpha, \beta > 0$). 
On the other hand, the embedding of $V$ into $2V$, used in Lemma
\ref{tauVSsigma}, represents $\tau_{j-1,j}$ by a word of length $2j-3$.

\bigskip

\noindent {\bf ``Why'' is the word problem of {\boldmath $2V$} 
{\sf coNP}-complete?}
The table-size of an element $f \in 2V$ can be exponentially larger than the 
word length of $f$ (over a finite generating set); 
hence, the polynomial-time algorithm for the word problem of $V$ (consisting 
of simply composing the tables of the generators) turns into an 
exponential-time algorithm in $2V$. 
In $V$ we have the table-size formula 
$\,|{\rm domC}(f_2 \circ f_1)| \le |{\rm domC}(f_2)| + |{\rm domC}(f_1)|$; 
in $2V$ there is no such formula.
However, the length-formula of Lemma \ref{ellAdd} implies rather directly 
that the word problem of $2V$ belongs to {\sf coNP}. 
 
The {\sf coNP}-hardness is less intuitive. The proof that $V$ (over 
$\Gamma_{\! V} \cup \tau$) can simulate circuits is intuitive (if tedious).
The commutation test, reducing a generalized word problem to a word problem, 
is less intuitive, and it is a priori not related to computing. 
The alternative proof of {\sf coNP}-hardness of the word problem of $V$ over
$\Gamma \cup \tau$ is derived from Jordan's theorem, which is itself based 
on Barrington's theorem; the latter has always been considered a surprising 
result. 

Using the shift to represent the infinite set $\tau$ by a finite set is 
easy.  But the shift is not a circuit element (although it has a 
computational meaning, namely, as an operation in multi-stack machines).

\bigskip

\noindent {\bf Acknowledgement:} I would like to thank the referee for a 
thorough reading of the paper.


\bigskip



\bigskip 

{\small 
\noindent
birget@camden.rutgers.edu
}

\end{document}